\documentclass[12pt,twoside, leqno]{article} 

\usepackage{ifthen}

\newboolean{tab}
\setboolean{tab}{false}

\usepackage[latin1]{inputenc}
\usepackage{graphicx}
\usepackage{latexsym}
\usepackage{amssymb}
\usepackage{epsfig}
\usepackage{enumitem}

\newcommand{\Var}{\mbox{Var}}
\def\eps{\varepsilon}
\def \R{\mathbb{R}}
\def \N{\mathbb{N}}

\def \E{\mathbb{E}}
\def \Cov{\mbox{Cov}}

\newcommand{\nb}{\nonumber}

\newcommand{\ef}{\mathcal{F}}

\newcommand{\Y}{\mathcal{Y}}
\newcommand{\G}{\mathcal{G}}

\newcommand{\XX}{\mathbf{X}}
\newcommand{\YY}{\mathbf{Y}}
\newcommand{\WW}{\mathbf{W}}

\newcommand{\bean}{\begin{eqnarray*}}
\newcommand{\eean}{\end{eqnarray*}}
\newcommand{\bea}{\begin{eqnarray}}
\newcommand{\eea}{\end{eqnarray}}
\newcommand{\be}{\begin{eqnarray}}
\newcommand{\ee}{\end{eqnarray}}
\newcommand{\beq}{\begin{equation}}
\newcommand{\eeq}{\end{equation}}

\begin{document}

\title{\bf The independence process in conditional quantile location-scale models and an application to testing for monotonicity}

\author{Melanie Birke\\Universitõt Bayreuth\\ 
\and Natalie Neumeyer\\
Universitõt Hamburg\\ 
\and Stanislav Volgushev\footnote{ The authors would like to thank two anonymous referees and the associate editor for careful reading and for very constructive suggestions to improve the paper. Our special thanks go to one of the referees for several very careful readings of the manuscript and insightful comments. Part of this work was conducted while Stanislav Volgushev was postdoctoral fellow at the Ruhr University Bochum, Germany. During that time Stanislav Volgushev was supported by the Sonderforschungsbereich ``Statistical modelling of nonlinear dynamic processes" (SFB~823), Teilprojekt (C1), of the Deutsche Forschungsgemeinschaft.}
\\University of Toronto\\
}

\maketitle

\newtheorem{theo}{Theorem}[section]
\newtheorem{lemma}[theo]{Lemma}
\newtheorem{cor}[theo]{Corollary}
\newtheorem{rem}[theo]{Remark}
\newtheorem{prop}[theo]{Proposition}
\newtheorem{defin}[theo]{Definition}
\newtheorem{example}[theo]{Example}

\begin{abstract} 
In this paper the nonparametric quantile regression model is considered in a location-scale context. The asymptotic properties of the empirical independence process based on covariates and estimated residuals are investigated. In particular an asymptotic expansion and weak convergence to a Gaussian process are proved. 
The results can, on the one hand, be applied to test for validity of the location-scale model. On the other hand, they allow to derive various specification tests in conditional quantile location-scale models. 
In detail a test for monotonicity of the conditional quantile curve is investigated. 
For the test for validity of the location-scale model as well as for the monotonicity test smooth residual bootstrap versions of Kolmogorov-Smirnov and CramÚr-von Mises type test statistics are suggested. We give rigorous proofs for bootstrap versions of the weak convergence results. The performance of the tests is demonstrated in a simulation study. 
\end{abstract}



\newpage

\section{Introduction}
\def\theequation{1.\arabic{equation}}
\setcounter{equation}{0}

Quantile regression was introduced by Koenker and Bassett (1978) as an extension of least squares methods focusing on the estimation of the conditional mean function. Due to its many attractive features as robustness with respect to outliers and equivariance under monotonic transformations that are not shared by the mean regression, it has since then become increasingly popular in many important fields such as medicine, economics and environment modelling [see Yu et al.\ (2003) or Koenker (2005)]. Another important feature of quantile regression is its great flexibility. While mean regression aims at modelling the average behaviour of a variable $Y$ given a covariate $X=x$, quantile regression allows to analyse the impact of $X$ in different regions of the distribution of $Y$ by estimating several quantile curves simultaneously. See for example Fitzenberger et al.\ (2008), who demonstrate that the presence of certain structures in a company can have different effects on upper and lower wages. For a more detailed discussion, we refer the interested reader to the recent monograph by Koenker (2005).

The paper at hand has a twofold aim. On the one hand it proves a weak convergence result for the empirical independence process of covariates and estimated errors in a nonparametric location-scale conditional quantile model. On the other hand it suggests a test for monotonicity of the conditional quantile curve. To the authors' best knowledge this is the first time that those problems are treated for the general nonparametric quantile regression model. 

The empirical independence process results from the distance of a joint empirical distribution function and the product of the marginal empirical distribution functions. It can be used to test for independence; see Hoeffding (1948), Blum et al.\ (1961) and ch.\ 3.8 in van der Vaart and Wellner (1996). When applied to covariates $X$ and estimators of error terms $\eps=(Y-q(X))/s(X)$ it can be used to test for validity of a location-scale model $Y=q(X)+s(X)\eps$ with $X$ and $\eps$ independent. Here the conditional distribution of $Y$, given $X=x$, allows for a location-scale representation $P(Y\leq y\mid X=x)=F_\eps((y-q(x))/s(x))$, where $F_\eps$ denotes the error distribution function. To the best of our knowledge, Einmahl and Van Keilegom (2008a) is the only paper that considers such tests for location-scale models in a very general setting (mean regression, trimmed mean regression,\dots). However, the assumptions made there rule out the quantile regression case, where $q$ is defined via $P(Y\leq q(x)\mid X=x)=\tau$ for some $\tau\in(0,1)$, $\forall x$. The first part of our paper can hence be seen as extension and completion of the results by Einmahl and Van Keilegom (2008a). Plenty of technical effort was necessary to obtain the weak convergence result in the quantile context (see the proof of Theorem \ref{theo1} below). Validity of a location-scale model means that the covariates have influence on the trend and on the dispersion of the conditional distribution of $Y$, but otherwise do not affect the shape of the conditional distribution (such models are frequently used, see Shim et al., 2009, and Chen et al., 2005). Contrariwise if the test rejects independence of covariates and errors then there is evidence that the influence of the covariates on the response goes beyond location and scale effects. Note that our results easily can be adapted to test the validity of location models $P(Y\leq y\mid X=x)=F_\eps(y-q(x))$; see also Einmahl and Van Keilegom (2008b) and Neumeyer (2009b) in the mean regression context.
\\
Further if there is some evidence that certain quantile curves might be monotone one should check by a statistical test, that this assumption is reasonable. Such evidence can e.g.\ come from an economic, physical or biological background. In classical mean regression there are various methods for testing monotonicity. It has already been considered e.g.\ in Bowman et al.\ (1998), Gijbels et al.\ (2000), Hall and Heckman (2001), Goshal et al.\ (2000), Durot (2003), Baraud et al.\ (2003) or Dom\'inguez-Menchero et al.\ (2005) and Birke and Dette (2007). More recent work on testing monotonicity is given in Wang and Meyer (2011) who use regression splines and use the minimum slope in the knots as test criterion, and Birke and Neumeyer (2013) who use empirical process techniques for residuals built from isotonized estimators. While most of the tests are very conservative and not powerful against alternatives with only a small deviation from monotonicity the method proposed by Birke and Neumeyer (2013) has in some situations better power than the other tests and can also detect local alternatives of order $n^{-1/2}$.
While there are several proposals for monotone estimators of a quantile function (see e.g.\ Cryer et al. (1972) or Robertson and Wright (1973) for median regression and Casady and Cryer (1976) or Abrevaya (2005) for general quantile regression), the problem of testing whether a given quantile curve is increasing (decreasing) has received nearly no attention in the literature. Aside from the paper by Duembgen (2002) which deals with the rather special case of median regression in a location model, the authors - to the best of their knowledge - are not aware of any tests for monotonicity of conditional quantile curves. The method, which is introduced here is based on the independence process considered before. 
Note that the test is not the same as the one considered by  Birke and Neumeyer (2013) for mean regression adapted to the quantile case. It turned out that in quantile regression the corresponding statistic would not be suitable for constructing a statistical test (see also Section \ref{sec-mon}).\\ 
The paper is organized as follows. In Section 2 we present the location-scale model, give necessary assumptions and define the estimators. In Section 3 we introduce the independence process, derive asymptotical results and construct a test for validity of the model. Bootstrap data generation and asymptotic results for a bootstrap version of the independence process are discussed as well. The results derived there are modified in Section 4 to construct a test for monotonicity of the quantile function. In Section 5 we present a small simulation study while we conclude in Section 6. All proofs are deferred to an appendix and supplementary material.


\section{The location-scale model, estimators and assumptions}\label{sec-est}
\def\theequation{2.\arabic{equation}}
\setcounter{equation}{0}

For some fixed $\tau\in (0,1)$, consider the nonparametric quantile regression model of location-scale type [see e.g. He (1997)], 
\begin{eqnarray}\label{mod-het}
Y_i &=& q_\tau(X_i) + s(X_i)\eps_i, \quad\quad i = 1, \ldots , n, 
\end{eqnarray}
where $q_\tau(x)=F_Y^{-1}(\tau|x)$ is the $\tau$-th conditional quantile function, $(X_i, Y_i)$, $i=1,\ldots,n$, is a bivariate sample of i.i.d.\ observations and $F_Y(\cdot|x)=P(Y_i\leq\cdot|X_i=x)$ denotes the conditional distribution function of $Y_i$ given $X_i=x$. Further, $s(x)$ denotes the median of $|Y_i-q_\tau(X_i)|$, given $X_i=x$. 
We assume that $\eps_i$ and $X_i$ are independent and, hence, that $\eps_i$ has $\tau$-quantile zero and $|\eps_i|$ has median one, because
\begin{eqnarray*}
\tau &=& P\Big(Y_i\leq q_\tau(X_i)\,\Big|\; X_i=x\Big)\;=\; P(\eps_i\leq 0)\\
\frac{1}{2} &=& P\Big(|Y_i-q_\tau(X_i)|\leq s(X_i)\,\Big|\; X_i=x\Big)\;=\; P(|\eps_i|\leq 1).
\end{eqnarray*}
Denote by $F_\eps$ the distribution function of $\eps_i$. Then for the conditional distribution we obtain a location-scale representation as $F_Y(y|x)=F_\eps((y-q_\tau(x))/s(x))$, where $F_\eps$ as well as $q_\tau$ and $s$ are unknown. 

For example, consider the case $\tau=\frac{1}{2}$. Then we have a median regression model, which allows for heteroscedasticity in the sense, that the conditional median absolute deviation $s(X_i)$ of $Y_i$, given $X_i$, may depend on the covariate $X_i$. Here the median absolute deviation of a random variable $Z$ is defined as 
$\mbox{MAD}(Z) = \mbox{median} ( |Z-\mbox{median}(Z)|)$
and is the typical measure of scale (or dispersion), when the median is used as location measure. 
This heteroscedastic median regression model is analogous to the popular heteroscedastic mean regression model 
$Y_i = m(X_i) + \sigma(X_i)\eps_i$, $i = 1, \ldots , n$,
where $X_i$ and $\eps_i$ are assumed to be independent, $E[\eps_i]=0$, $\mbox{sd}(\eps_i)=1$, and hence, $m(x)=E[Y_i\mid X_i=x]$, $\sigma(x)=\mbox{sd}(Y_i\mid X_i=x)$ (see among many others e.g.\ Efromovich (1999), chapter 4.2 for further details).

\begin{rem}\rm \label{reskal}
Note that assuming $|\eps_i|$ to have median one is not restrictive. More precisely, if the model $Y_i = q_\tau(X_i) + \tilde s(X_i)\eta_i$ with $\eta_i$ i.i.d.\ and independent of $X_i$ and some positive function $\tilde s$ holds, the model $Y_i = q_\tau(X_i) + s(X_i)\eps_i$ with $s(X_i) := \tilde s(X_i)F_{|\eta|}^{-1}(1/2)$, $\eps_i := \eta_i/F_{|\eta|}^{-1}(1/2)$ will also be true, where $F_{|\eta|}$ denotes the distribution function of $|\eta_i|$. Then in particular $P(|\eps_i| \leq 1) = P(|\eta_i| \leq F_{|\eta|}^{-1}(1/2)) = 1/2$.
$\blacksquare$
\end{rem}

In the literature, several non-parametric quantile estimators have been proposed [see e.g.\ Yu and Jones (1997, 1998), Takeuchi et al.\ (2006) or Dette and Volgushev (2008), among others]. In this paper we follow the last-named authors who proposed non-crossing estimates of quantile curves using a simultaneous inversion and isotonization of an estimate of the conditional distribution function. To be precise, let
\be \label{def-Fdach}
\hat F_Y(y|x) := (\XX^t\WW\XX)^{-1}\XX^t\WW\YY
\ee
with
\bean
\XX &=& \left(
\begin{array}{cccc}
1 & (x-X_1)& ... & (x-X_1)^p\\
\vdots & \vdots & ... & \vdots\\
1 & (x-X_n)& ... & (x-X_n)^p
\end{array}
\right), \qquad
\YY := \Big(\Omega \Big(\frac{y-Y_1}{d_n}\Big),\dots,\Omega \Big(\frac{y-Y_n}{d_n}\Big) \Big)^t 
\\
\WW &=& \mbox{Diag}\Big(K_{h_n,0}(x-X_1),...,K_{h_n,0}(x-X_n)\Big),
\eean
denote a smoothed local polynomial estimate (of order $p\geq 2$) of the conditional distribution function $F_Y(y|x)$ where $\Omega(\cdot)$ is a smoothed version of the indicator function and we used the notation
$K_{h_n,k}(x) := K(x/h_n)(x/h_n)^k$. Here $K$ denotes a nonnegative kernel and $d_n,h_n$ are bandwidths converging to 0 with increasing sample size. Note that the estimator $\hat F_Y(y|x)$ can be represented as weighted average 
\be\label{def-LL}
\hat F_Y(y|x) = \sum_{i=1}^n W_i(x)\Omega \Big(\frac{y-Y_i}{d_n}\Big).
\ee
Following Dette and Volgushev (2008) we consider a strictly increasing distribution function $G: \mathbb{R} \to (0,1)$, a nonnegative kernel
$\kappa$ and a bandwidth $b_n$, and define the functional
\[
H_{G,\kappa,\tau,b_n}(F) := \frac{1}{b_n}\int_0^1 \int_{-\infty}^\tau \kappa\Big(\frac{F(G^{-1}(u)) - v}{b_n} \Big) dvdu. 
\] 
Note that it is intuitively clear that $H_{G,\kappa,\tau,b_n}(\hat F_{Y}(\cdot|x))$, where $\hat F_{Y}$ is the estimator of the conditional distribution function defined in (\ref{def-Fdach}), is a consistent estimate of $H_{G,\kappa,\tau,b_n}(F_{Y}(\cdot|x))$. If $b_n \to 0$, this quantity 
can be approximated as follows
\begin{eqnarray*} 
H_{G,\kappa,\tau,b_n}(F_{Y}(\cdot|x)) &\approx&
\int_{\R}I\{F_{Y}(y|x)\leq \tau\}dG(y) \\
&=& \nonumber  \int_0^1I\{F_{Y}(G^{-1}(v)|x)\leq \tau\}dv
 \: = \: G \circ F_{Y}^{-1} (\tau|x),
\end{eqnarray*}
and as a consequence an estimate of the conditional quantile function $q_\tau(x)=F_{Y}^{-1}(\tau|x)$ can be defined by
$$ 
\hat q_\tau(x) :=  G^{-1}(H_{G,\kappa,\tau,b_n}(\hat F_{Y}(\cdot|x))).
$$
Finally, note that the scale function $s$ is the conditional median of the distribution of $|e_i|$, given the covariate $X_i$, where 
$e_i=Y_i- q_\tau(X_i)=s(X_i)\eps_i$, $i=1,\dots,n$.
Hence, we apply the quantile-regression approach to $|\hat e_i|=|Y_i- \hat q_\tau(X_i)|$, $i=1,\dots,n$, and obtain the estimator
\begin{equation} \label{def-szeta}
\hat s(x) = G_s^{-1}(H_{G_s,\kappa,1/2,b_n}(\hat F_{|e|}(\cdot|x)))
\: .
\end{equation}
Here $G_s: \R \rightarrow (0,1)$ is a strictly increasing distribution function and $\hat F_{|e|}(\cdot|x)$ denotes the estimator of the conditional distribution function $F_{|e|}(\cdot|x)=P(|e_i|\leq\cdot|X_i=x)$ of $|e_i|$, $i=1,\dots,n$, i.\,e.\
\begin{equation} \label{def-Fedach}
\hat F_{|e|}(y|x) = \sum_{i=1}^n W_i(x)I\{|\hat e_i| \leq y\}
\end{equation}
with the same weights $W_i$ as in (\ref{def-LL}). We further use the notation $F_{e}(\cdot|x)=P(e_1\leq\cdot|X_1=x)$. \\
\\
For a better overview and for later reference, below we collect all the technical assumptions concerning the estimators needed throughout the rest of the paper. First, we collect the assumptions needed for the kernel functions and functions $G,G_s$ used in the construction of the estimators.

\begin{enumerate}[label=(\textbf{K\arabic{*}})]
\item \label{as:k1} The function $K$ is a symmetric, positive, Lipschitz-continuous density with support $[-1,1]$. Moreover, the matrix $\mathcal{M}(K)$ with entries \[
(\mathcal{M}(K))_{k,l} = \mu_{k+l-2}(K) := \int u^{k+l-2}K(u)du
\] 
is invertible.
\item \label{as:k2} The function $K$ is two times continuously differentiable, $K^{(2)}$ is Lipschitz continuous, and for $m=0,1,2$ the set $\{x|K^{(m)}(x)>0\}$ is a union of finitely many intervals.
\item \label{as:om} The function $\Omega$ has derivative $\omega$ which has support $[-1,1]$, is a kernel of order $p_\omega$, and is two times continuously differentiable with uniformly bounded derivatives.
\item \label{as:kap} The function $\kappa$ is a symmetric, uniformly bounded density, and has one Lipschitz-continuous derivative.
\item \label{as:G} The function $G: \R \rightarrow [0,1]$ is strictly increasing. Moreover, it is two times continuously differentiable in a neighborhood of the set $Q := \{q_\tau(x)|x \in [0,1]\}$ and its first derivative is uniformly bounded away from zero on $Q$.
\item \label{as:Gs} The function $G_s: \R \rightarrow (0,1)$ is strictly increasing. Moreover, it is two times continuously differentiable in a neighborhood of the set $S := \{s(x)|x \in [0,1]\}$ and its first derivative is uniformly bounded away from zero on $S$.
\end{enumerate}

The data-generating process needs to satisfy the following conditions.

\begin{enumerate}[label=(\textbf{A\arabic{*}})]
\item \label{as:fx} $X_1,\dots,X_n$ are independent and identically distributed with distribution function $F_X$ and Lipschitz-continuous density $f_X$ with support $[0,1]$ that is uniformly bounded away from zero and infinity.
\item \label{as:s} The function $s$ is uniformly bounded and $\inf_{x\in [0,1]} s(x) = c_s > 0$.
\item \label{as:fdiff} The partial derivatives $\partial_x^k\partial_y^l F_{Y}(y|x), \partial_x^k\partial_y^l F_{e}(y|x)$ exist and are continuous and uniformly bounded on $\R\times [0,1]$ for $k\vee l\leq 2$ or $k+l \leq d$ for some $d\geq 3$.
\item \label{as:eps} The errors $\eps_1,\dots,\eps_n$ are independent and identically distributed with strictly increasing distribution function $F_\eps$ (independent of $X_i$) and density $f_\eps$, which is positive everywhere and continuously differentiable such that $\sup_{y\in\R}|yf_\eps(y)| <\infty$ and $\sup_{y\in\R}|y^2f_\eps'(y)|<\infty$. The $\eps_i$ have $\tau$-quantile zero and $F_{|\eps|}(1) = 1/2$, that is $|\eps_1|$ has median one.
\item \label{as:fbound} For some $\alpha>0$ we have $\sup_{u,y}|y|^\alpha (F_{Y}(y|u)\wedge(1-F_{Y}(y|u))) < \infty$.
\end{enumerate}

Finally, we assume that the bandwidth parameters satisfy 
\begin{enumerate}[label=(\textbf{BW})]
\item \label{as:bw} 
$\displaystyle 
\frac{\log n}{nh_n(h_n\wedge d_n)^4} = o(1),\quad \frac{\log n}{nh_n^2b_n^2} = o(1),\quad d_n^{2(p_\omega \wedge d)} + h_n^{2((p+1)\wedge d)} + b_n^4= o(n^{-1}),
$

with $p_\omega$ from (K3), $d$ from (A3) and $p$ the order of the local polynomial estimator in (\ref{def-Fdach}).
\end{enumerate}

\begin{rem}\rm \label{rem:qdiff} 
Assumptions \ref{as:fx} and \ref{as:s} are mild regularity assumptions on the data-generating process. Assumption \ref{as:fbound} places a very mild condition on the tails of the error distribution, and is satisfied even for distribution functions that don't have finite first moments. Assumptions \ref{as:fdiff} and \ref{as:eps} are probably the strongest ones. Note that by the implicit function theorem they imply that $x\mapsto q_\tau(x)$ and $x\mapsto s(x)$ are $2$ times continuously differentiable with uniformly bounded derivatives. 
Those assumptions play a crucial role throughout the proofs. In principle, this kind of condition is quite standard in the non-parametric estimation and testing literature. Note that due to the additional smoothing of $\hat F_Y(y|x)$ in $y$-direction, we require more than the existence of just all the second-order partial derivatives of $F_Y(y|x)$. The smoothing is necessary for the proofs, and it leads to a slightly better finite-sample performance of the testing procedures. Regarding the bandwidth assumption \ref{as:bw}, observe that if for example $d=p_\omega = p =3$ and we set $d_n = h_n = n^{-1/6-\beta}$ for some $\beta \in (0,1/30)$, $b_n = h_n^{-1/4-\alpha}$ such that $\alpha + \beta \in (0,1/12)$, condition \ref{as:bw} holds. 
$\blacksquare$ \end{rem}


\section{The independence process, asymptotic results and testing for model validity}\label{sec-asy}
\def\theequation{3.\arabic{equation}}
\setcounter{equation}{0}

As estimators for the errors we build residuals 
\begin{eqnarray}\label{pse-neu}
\hat\eps_{i}=\frac{Y_i-\hat q_\tau(X_i)}{\hat s(X_i)},\quad i=1,\dots,n.
\end{eqnarray}
In the definition of the process on which test statistics are based we only consider those observations $(X_i,Y_i)$ such that $2h_n\leq X_i\leq 1-2h_n$ in order to avoid boundary problems of the estimators. The reason is that we first use $\hat q_\tau$ to build the residuals $\hat e_i$. For this, we need $h_n\leq X_i\leq 1-h_n$. The estimator $\hat s$ based on the pairs $(X_i,|\hat e_i|)$ is then used in the definition of the residuals $\hat \eps_i$. The estimation of $s$ requires us to again stay away from boundary points and thus we use the restriction  $2h_n\leq X_i\leq 1-2h_n$. 

For $y\in\mathbb{R}$, $t\in[2h_n,1-2h_n]$ we define the
joint empirical distribution function of pairs of covariates and residuals as
\begin{eqnarray}\label{FXen}
\hat F_{X,\eps,n}(t,y) &:=& \sum_{i=1}^nI\{\hat\eps_{i}\leq y\}I\{2h_n< X_i\leq t\}\frac{1}{\sum_{i=1}^nI\{2h_n< X_i\leq 1-2h_n\}}\\
&=& \frac{1}{n}\sum_{i=1}^nI\{\hat\eps_{i}\leq y\}I\{2h_n< X_i\leq t\}\frac{1}{\hat F_{X,n}(1-2h_n)-\hat F_{X,n}(2h_n)},\nb
\end{eqnarray}
where $\hat F_{X,n}$ denotes the usual empirical distribution function of the covariates $X_1,\dots,X_n$.
The empirical independence process compares the joint empirical distribution with the product of the corresponding marginal distributions. We thus define 
\begin{eqnarray}\label{Sn-neu}
S_n(t,y) &=& \sqrt{n}\Big(\hat F_{X,\eps,n}(t,y)-\hat F_{X,\eps,n}(1-2h_n,y)\hat F_{X,\eps,n}(t,\infty)\Big)
\end{eqnarray}
for $y\in\mathbb{R}$, $t\in[2h_n,1-2h_n]$, and $S_n(t,y)=0$ for $y\in\mathbb{R}$, $t\in[0,2h_n)\cup (1-2h_n,1]$. 
In the following theorem we state a weak convergence result for the independence process.

\begin{theo}\label{theo1} Under the location-scale model (\ref{mod-het}) and assumptions \ref{as:k1}-\ref{as:Gs}, \ref{as:fx}-\ref{as:fbound} and \ref{as:bw} we have the asymptotic expansion
\begin{eqnarray*}
S_n(t,y)&=& \frac{1}{\sqrt{n}}\sum_{i=1}^n\Big( I\{\eps_i\leq y\}-F_\eps(y)
-\phi(y)\Big(I\{\eps_i\leq 0\}-\tau\Big)-\psi(y)\Big(I\{|\eps_i|\leq 1\}-\frac{1}{2}\Big)\Big)\\
&&{}\qquad\quad\times\Big(I\{X_i\leq t\}-F_X(t)\Big)
 +o_P(1)
\end{eqnarray*}
uniformly with respect to $t\in[0,1]$ and $y\in\R$, where
\begin{eqnarray*}
\phi(y)&=&\frac{f_\eps(y)}{f_\eps(0)}\Big(1-y\frac{f_\eps(1)-f_\eps(-1)}{f_{|\eps|}(1)}\Big)
\, ,\quad \psi(y)\;=\;\frac{yf_\eps(y)}{f_{|\eps|}(1)} 
\end{eqnarray*}
and  $f_{|\eps|}(y)=(f_\eps(y)+f_\eps(-y))I_{[0,\infty)}(y)$ is the density of $|\eps_1|$.
The process $S_n$ converges weakly in $\ell^\infty([0,1]\times\mathbb{R})$ to a centered Gaussian process $S$ with covariance \rm
\begin{eqnarray*}
&& \Cov(S(s,y),S(t,z)) \;=\; (F_X(s\wedge t)-F_X(s)F_X(t))\\
&&{}\times\Big[ F_\eps(y\wedge z)-F_\eps(y)F_\eps(z) +\phi(y)\phi(z)(\tau-\tau^2)
+\frac{1}{4}\psi(y)\psi(z)\\
&&{}\quad-\phi(y)(F_\eps(z\wedge 0)-F_\eps(z)\tau)-\phi(z)(F_\eps(y\wedge 0)-F_\eps(y)\tau)\\
&&{}\quad-\psi(y) \Big( (F_\eps(z\wedge 1)-F_\eps(-1))I\{z>-1\}-\frac{1}{2}F_\eps(z)\Big)\\
&&{}\quad-\psi(z) \Big( (F_\eps(y\wedge 1)-F_\eps(-1))I\{y>-1\}-\frac{1}{2}F_\eps(y)\Big)\\
&&{}\quad
+(\phi(y)\psi(z)+\phi(z)\psi(y))\Big(F_\eps(0)-F_\eps(-1)-\frac{1}{2}\tau\Big)
\Big].
\end{eqnarray*}\it
\end{theo}
 
The proof is given in Appendix \ref{app-A}. 

\begin{rem}\rm\label{location-model}
The result can easily be adapted for location models $Y_i=q_\tau(X_i)+\eps_i$ with $\eps_i$ and $X_i$ independent. To this end we just set $\hat s\equiv 1$ in the definition of the estimators. The asymptotic covariance in Theorem \ref{theo1} then simplifies because the function $\phi$ reduces to $\phi(y)=f_\eps(y)/f_\eps(0)$ and $\psi(y)\equiv 0$. 
$\blacksquare$
\end{rem}

In the remainder of this section we discuss how the asymptotic result can be applied to test for validity of the location-scale model, i.\,e.\ testing the null hypothesis of independence of error $\eps_i$ and covariate $X_i$ in model (\ref{mod-het}). 

\begin{rem}\rm \label{rem-alt}
Assume that the location-scale model is not valid, i.\,e.\ $X_i$ and $\eps_i$ are dependent, but the other assumptions of Theorem \ref{theo1} are valid, where \ref{as:eps} is replaced by
\begin{enumerate}[label=(\textbf{A4'})]
\item The conditional error distribution function $F_\eps(\cdot|x)=P(\eps_i\leq \cdot|X_i=x)$ fulfills $F_\eps(0|x)=\tau$ and $F_\eps(1|x)-F_\eps(-1|x)=\frac 12$ for all $x$. It is strictly increasing and differentiable with density $f_\eps(\cdot|x)$ such that $\sup_{x,y}|yf_\eps(y|x)|<\infty$. 
\end{enumerate}
Then one can show that $S_n(t,y)/n^{1/2}$ converges in probability to $P(\eps_i\leq y,X_i\leq t)-F_\eps(y)F_X(t)$, uniformly with respect to $y$ and $t$. 
$\blacksquare$
\end{rem}

\begin{rem}\rm \label{tau-alpha}
If the location-scale model is valid for some $\tau$-th quantile regression function it is valid for every $\alpha$-th quantile regression function, $\alpha\in(0,1)$. This easily follows from $q_\alpha(x)=F_\eps^{-1}(\alpha)s(x)+q_\tau(x)$ which is a consequence from the representation of the conditional distribution function $F_Y(y|x)=F_\eps((y-q_\tau(x))/s(x))$ (compare Remark \ref{reskal}). A similar statement is even true for general location and scale measures, see e.\,g.\ Van Keilegom (1998), Prop.\ 5.1. Thus for testing the validity of the location-scale model one can restrict oneself to the median case $\tau=0.5$. 
$\blacksquare$
\end{rem}

\begin{rem}\rm 
Einmahl and Van Keilegom (2008a) consider a process similar to $S_n$ for general location and scale models. They define $q(x)=\int_0^1 F^{-1}(s|x)J(s)\,ds$ and $s^2(x)=\int_0^1 (F^{-1}(s|x))^2J(s)\,ds-q^2(x)$ with score function $J$, which rules out the quantile case $q(x)=F^{-1}(\tau|x)$. Einmahl and Van Keilegom (2008a) show that estimation of the errors has no influence in their context, i.\,e.\ they obtain a scaled completely tucked Brownian sheet as limit process and thus asymptotically distribution-free tests. This is clearly not the case in Theorem \ref{theo1}. 
$\blacksquare$
\end{rem}

To test for the validity of a location-scale model we reject the null hypothesis of independence of $X_i$ and $\eps_i$ for large values of, e.\,g., the Kolmogorov-Smirnov statistic 
$$
K_n=\sup_{t\in[0,1],y\in\R}|S_n(t,y)|
$$
or the CramÚr-von Mises statistic 
$$ 
C_n = \int_{\mathbb{R}}\int_{[0,1]} S_n^2(t,y) \,\hat F_{X,n}(dt)\,\hat F_{\eps,n}(dy),
$$
where $\hat F_{\eps,n}(\cdot)=\hat F_{X,\eps,n}(1-2h_n,\cdot)$. From Theorem \ref{theo1} we obtain the following asymptotic distributions.

\begin{cor}\label{cor-KS-CvM} Under the assumptions of Theorem \ref{theo1} we have
\begin{eqnarray*}
K_n&\stackrel{d}{\longrightarrow}& \sup_{t\in[0,1],y\in\R}|S(t,y)|= \sup_{x\in[0,1],y\in\R}|S(F_X^{-1}(x),y)|\\
C_n&\stackrel{d}{\longrightarrow}& \int_{\mathbb{R}}\int_{[0,1]} S^2(t,y)F_X(dt)F_\eps(dy)=\int_{\mathbb{R}}\int_{[0,1]} S^2(F_X^{-1}(x),y)\,dx\,F_\eps(dy).
\end{eqnarray*}
\end{cor}

The proof is given in Appendix A. 
The asymptotic distributions of the test statistics are  independent from the covariate distribution $F_X$, but depend in a complicated manner on the error distribution $F_\eps$. To overcome this problem we 
suggest a bootstrap version of the test.
To this end let $\mathcal{Y}_n = \{(X_1,Y_1),\dots,$ $(X_n,Y_n)\}$ denote the original sample. 
We generate bootstrap errors  as $\eps_i^*=\tilde{\eps}_i^*+\alpha_nZ_i$ ($i=1,\dots,n$),
where $\alpha_n$ denotes a positive smoothing parameter, $Z_1,\dots,Z_n$ are independent, standard normally distributed random variables (independent of $\Y_n$) and $\tilde{\eps}_1^*,\dots,\tilde{\eps}_n^*$ are randomly drawn with replacement from the set of residuals $\{\hat{\eps}_j\mid j\in\{1,\dots,n\}, X_j\in (2h_n,1-2h_n]\}$.
Conditional on the original sample $\mathcal{Y}_n$
the random variables $\eps_1^*,\dots,\eps_n^*$ are i.i.d.\ with  distribution function 
\beq \label{eq:feps}
\tilde{F}_\eps(y) = \frac{\frac{1}{n}\sum_{i=1}^n \Phi\Big(\frac{y-\hat{\eps}_i}{\alpha_n}\Big)I\{2h_n< X_i\leq 1-2h_n \}}{\hat F_{X,n}(1-2h_n)-\hat F_{X,n}(2h_n)},
\eeq
where $\Phi$ denotes the standard normal distribution function. 
Note that the bootstrap error's $\tau$-quantile is not exactly zero, but vanishes asymptotically. 
We use a smooth distribution to generate new bootstrap errors because smoothness of the error distribution is a crucial assumption for the theory necessary to derive Theorem \ref{theo1}; see also Neumeyer (2009a).

Now we build new bootstrap observations,
\begin{eqnarray*}
Y_i^*&=& \hat{q}_\tau(X_i)+\hat s(X_i)\eps_i^*,\quad i=1,\dots,n.
\end{eqnarray*}
Let $\hat q_\tau^*$ and $\hat s^*$ denote the quantile regression and scale function estimator defined analogously to $\hat q_\tau$ and $\hat s$, but based on the bootstrap sample $(X_1,Y_1^*),\dots,(X_n,Y_n^*)$. 
Analogously to (\ref{Sn-neu}) the bootstrap version of the independence process is defined as
\begin{eqnarray*}
S_n^*(t,y) &=& \sqrt{n}\Big(\hat F^*_{X,\eps,n}(t,y)-\hat F^*_{X,\eps,n}(1-4h_n,y)\hat F^*_{X,\eps,n}(t,\infty)\Big)
\end{eqnarray*}
for $t\in[4h_n,1-4h_n]$, $y\in\R$, and $S_n^*(t,y)=0$ for $t\in[0,4h_n)\cup (1-4h_n,1]$, $y\in\R$. Here, similar to (\ref{FXen}),
\begin{eqnarray*}
\hat F^*_{X,\eps,n}(t,y)&=& \frac{1}{n}\sum_{i=1}^nI\{\hat\eps_{i}^*\leq y\}I\{4h_n< X_i\leq t\}\frac{1}{\hat F_{X,n}(1-4h_n)-\hat F_{X,n}(4h_n)},\nb
\end{eqnarray*}
with $\hat\eps_{i}^*=(Y_i^*-\hat q_{\tau}^*(X_i))/\hat s^*(X_i)$, $i=1,\dots,n$. 

To obtain the conditional weak convergence we need the following additional assumptions.

\begin{enumerate}[label=(\textbf{B\arabic{*}})]
\item \label{as:b0} \label{as:b1} We have for some $\delta>0$
\[
\frac{nh_n^2\alpha_n^2}{\log h_n^{-1}\log n}\to\infty, \quad \frac{n\alpha_nh_n}{\log n}\to\infty,\quad
\frac{h_n}{\log n}=O(\alpha_n^{8\delta/3}), \quad n\alpha_n^4 = o(1)
\]
and there exists a $\lambda>0$ such that 
\[
\frac{nh_n^{1+\frac{1}{\lambda}}\alpha_n^{2+\frac{2}{\lambda}}}{\log h_n^{-1}(\log n)^{1/\lambda}}\to\infty.
\]
\item \label{as:b2} Let $E[|\eps_1|^{\max (\upsilon,2\lambda)}]<\infty$ for some  $\upsilon > 1+2 / \delta$ and with $\delta$ and $\lambda$ from assumption \ref{as:b1}. 
\end{enumerate}

Here, \ref{as:b2} can be relaxed to $E[|\eps_1|^{2\lambda}]<\infty$ if the process is only considered for $y\in[-c,c]$ for some $c>0$ instead of for $y\in\R$.

\begin{theo}\label{theo1-boot} Under the location-scale model (\ref{mod-het}) and assumptions \ref{as:k1}-\ref{as:Gs}, \ref{as:fx}-\ref{as:fbound}, \ref{as:bw} and \ref{as:b0}-\ref{as:b2}
conditionally on $\Y_n$, the process 
$S_n^*$  converges weakly in $\ell^\infty([0,1]\times\R)$ to the Gaussian process $S$ defined in Theorem \ref{theo1}, in probability. 
\end{theo}

A rigorous proof is given in Appendix \ref{app-boot}. 

\begin{rem}\rm \label{rem-bootKS}
Recall that the Kolmogorov-Smirnov test statistic is given by $K_n=\sup_{t,y}|S_n(t,y)|$ and define its bootstrap version as $K_n^*=\sup_{t,y}|S_n^*(t,y)|$. Let the critical value $k^*_{n,1-\alpha}$ be obtained from 
$$P ( K_n^* \geq k^*_{n,1-\alpha}\mid \Y_n)=1-\alpha,$$
and reject the location-scale model if $K_n\geq k^*_{n,1-\alpha}$. Then from Theorems \ref{theo1} and \ref{theo1-boot} it follows that the test has asymptotic level $\alpha$. Moreover if the location-scale model is not valid by Remark \ref{rem-alt} we have $K_n\to\infty$ in probability, whereas with the same methods as in the proof of Theorem \ref{theo1-boot} it can be shown that $k^*_{n,1-\alpha}$ converges to a constant. Thus the power of the test converges to one. A similar reasoning applies for the CramÚr-von Mises test. 
The finite sample performance of the bootstrap versions of both tests is studied in Section 5. 
$\blacksquare$
\end{rem}

\begin{rem}\rm  Recently, Sun (2006) and Feng, He and Hu (2011) proposed to use wild bootstrap in the setting of quantile regression. To follow the approach of the last-named authors, one would define $\eps_i^*=v_i\hat\eps_i$ such that $P^*(v_i\hat\eps_i\leq 0|X_i)=\tau$, e.\,g.\
$$v_i=\pm 1\mbox{ with probability } \left\{\begin{array}{cl} {1-\tau\atop\tau} & \mbox{ if } \hat\eps_i\geq 0\\
{\tau\atop 1-\tau} & \mbox{ if } \hat\eps_i< 0.\end{array}\right.$$
However, then when calculating the conditional asymptotic covariance (following the proof in Appendix \ref{app-boot}), instead of $\tilde F_\eps(y)$ the following term appears
$$\frac{1}{n}\sum_{i=1}^n P(v_i\hat\eps_i\leq y\mid\Y_n) \stackrel{n\to\infty}{\longrightarrow} (1-\tau)(F_\eps(y)-F_\eps(-y))+\tau.
$$
One obtains  $F_\eps(y)$ (needed to obtain the same covariance as in Theorem \ref{theo1}) only for $y=0$ or for median regression ($\tau=0.5$) with symmetric error distributions, but not in general. Hence, wild bootstrap cannot be applied in the general context of procedures using empirical processes in quantile regression. 
$\blacksquare$
\end{rem}

\begin{rem}\rm Under assumption of the location-scale model model (\ref{mod-het}) the result of Theorem \ref{theo1} can be applied to test for more specific model assumptions (e.\,g.\ testing goodness-of fit of a parametric model for the quantile regression function). The general approach is to build residuals $\hat\eps_{i,0}$ that only under $H_0$ consistently estimate the errors (e.\,g.\ using a parametric estimator for the conditional quantile function). Recall the definition of 
$\hat F_{X,\eps,n}$ in (\ref{FXen}) and define analgously $\hat F_{X,\eps_0,n}$ by using the residuals 
$\hat\eps_{i,0}$. Then, analogously to (\ref{Sn-neu}), define
\begin{eqnarray*}
S_{n,0}(t,y) &=& 
\sqrt{n}\Big(\hat F_{X,\eps_0,n}(t,y)-\hat F_{X,\eps,n}(1-2h_n,y)\hat F_{X,\eps,n}(t,\infty)\Big)
\end{eqnarray*}
for $y\in\mathbb{R}$, $t\in[2h_n,1-2h_n]$, and $S_{n,0}(t,y)=0$ for $y\in\mathbb{R}$, $t\in[0,2h_n)\cup (1-2h_n,1]$. 
With this process the discrepancy from the null hypothesis can be measured. 
This approach is considered in detail for the problem of testing monotonicity of conditional quantile functions in the next section. 

A related approach, which however does not assume the location-scale model, is suggested to test for significance of covariables in quantile regression models by Volgushev et al.\ (2013). 
$\blacksquare$
\end{rem}

\section{Testing for monotonicity of conditional quantile curves}\label{sec-mon}
\def\theequation{4.\arabic{equation}}
\setcounter{equation}{0}

\begin{figure}[t]
\epsfig{file=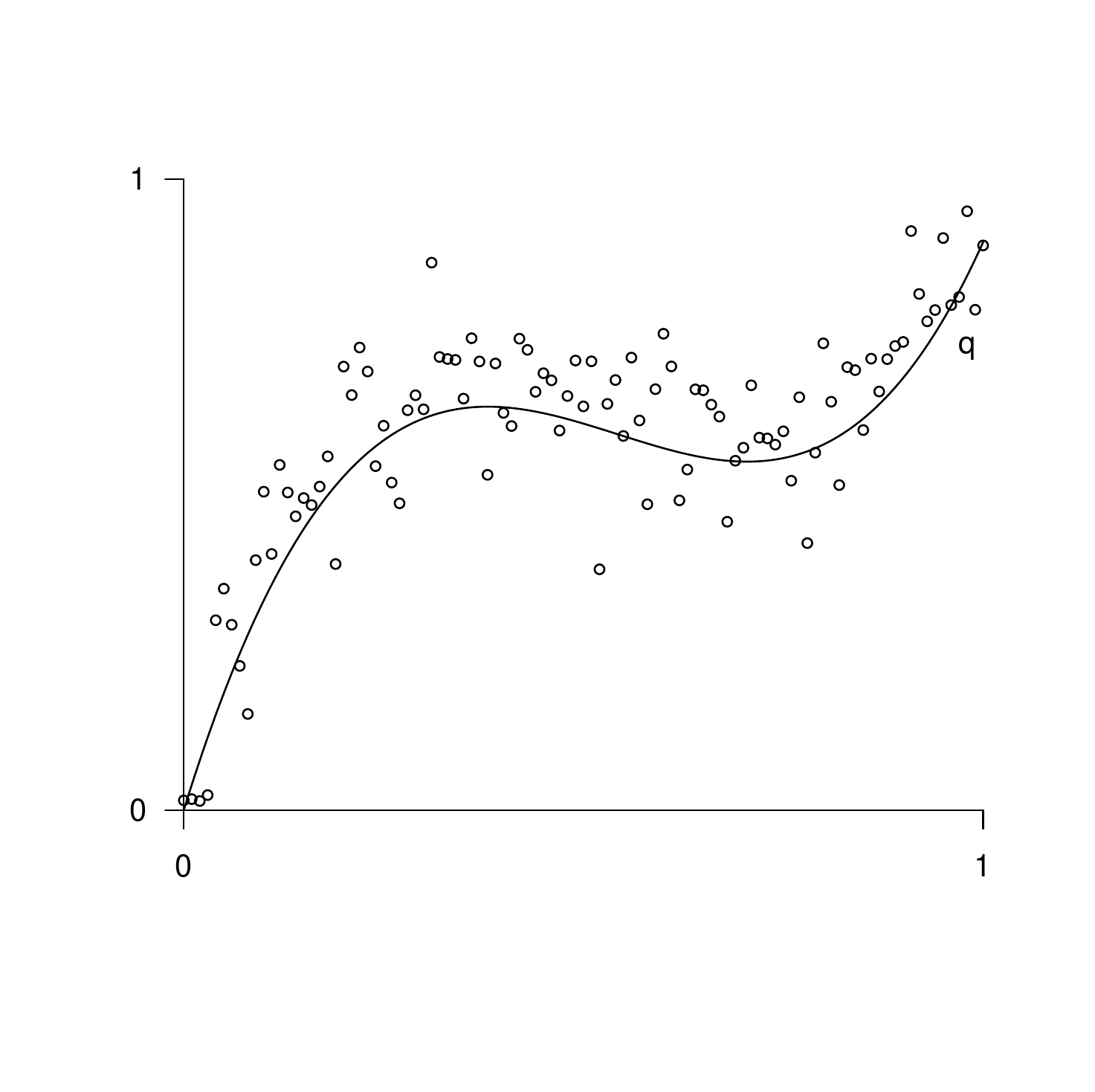,width=8.5cm} \epsfig{file=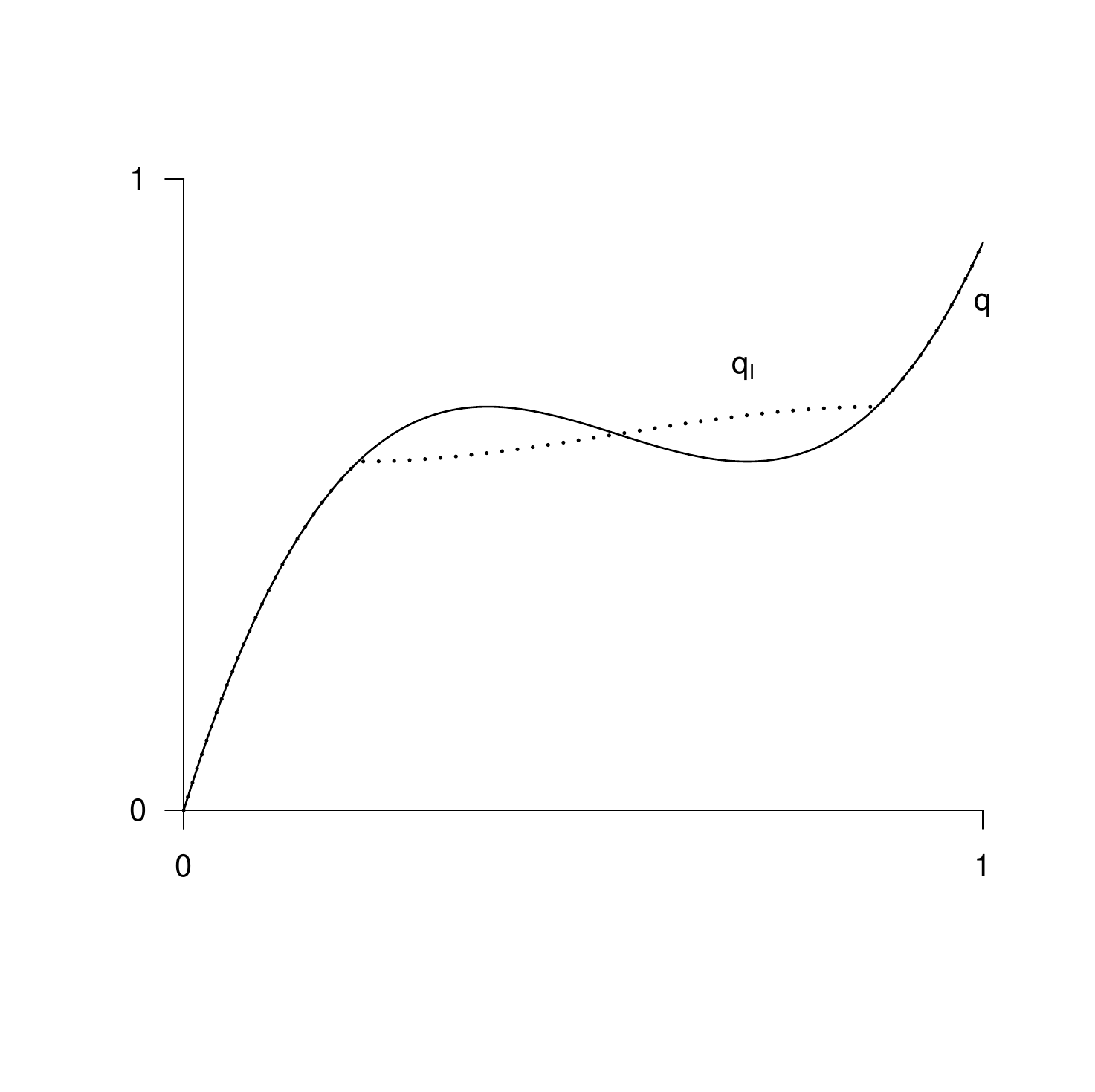,width=8.5cm}\vspace{-1.5cm}
\caption{\label{fig1}\sl Left part: True nonincreasing function $q_\tau$ for $\tau=0.25$ with scatter-plot of a typical sample. Right part: $q_\tau$ (solid line) and increasing rearrangement $q_{\tau,I}$ (dotted line).}
\end{figure}

In this section, we consider a test for the hypothesis 
$$H_0: q_\tau(x)\mbox{ is increasing in $x$}.$$
To this end we define an increasing estimator $\hat q_{\tau,I}$, which consistently estimates $q_\tau$ if the hypothesis $H_0$ is valid, and consistently estimates some increasing function $q_{\tau, I}\neq q_\tau$ under the alternative that $q_\tau$ is not increasing. 
For any function $h:[0,1]\to\R$ define the increasing rearrangement on $[a,b]\subset [0,1]$ as the function $\Gamma(h):[a,b]\to\R$ with
\begin{eqnarray*}
\Gamma(g)(x) &=& \inf\Big\{z\in\R\;\Big| a + \int_{a}^{b} I\{g(t)\leq z\}\, dt\geq x\Big\}.
\end{eqnarray*}
Note that if $g$ is increasing, then $\Gamma(g)=g|_{[a,b]}$. See Anevski and FougÞres (2007) and Neumeyer (2007) who consider increasing rearrangements of curve estimators in order to obtain monotone versions of unconstrained estimators.
We denote by $\Gamma_n$ the operator $\Gamma$ with $[a,b]=[h_n,1-h_n]$.
We define the increasing estimator as $\hat q_{\tau,I}=\Gamma_n(\hat q_\tau)$, where $\hat q_\tau$  denotes the unconstrained estimator of $q_\tau$ that was defined in Section \ref{sec-est}. The quantity $\hat q_{\tau,I}$ estimates the increasing rearrangement $q_{\tau,I}=\Gamma(q_\tau)$ of $q_\tau$ (with $[a,b]=[0,1]$). Only under the hypothesis $H_0$ of an increasing regression function we have $q_\tau=q_{\tau,I}$. In Figure \ref{fig1} (right part) a non-increasing function $q_\tau$ and its increasing rearrangement $q_{\tau,I}$ are displayed. 

Now we build (pseudo-) residuals 
\begin{eqnarray}\label{pse}
\hat\eps_{i,I}=\frac{Y_i-\hat q_{\tau,I}(X_i)}{\hat s(X_i)},
\end{eqnarray}
which estimate 
pseudo-errors $\eps_{i,I}=(Y_i-q_{\tau,I}(X_i))/s(X_i)$ that coincide with the true errors $\eps_i=(Y_i-q_\tau(X_i))/s(X_i)$ ($i=1,\dots,n$) in general only under $H_0$. Note that we use $\hat s$ from (\ref{def-szeta}) for the standardization and not an estimator built from the constrained residuals. 
 Let further 
$\hat\eps_{i}$ denote the unconstrained residuals as defined in (\ref{pse-neu}). 
The idea for the test statistic we suggest is the following. Compared to the true errors $\eps_1,\dots,\eps_n$, which are assumed to be i.i.d., the pseudo-errors $\eps_{1,I},\dots,\eps_{n,I}$ behave differently. 
\begin{figure}[t]
\epsfig{file=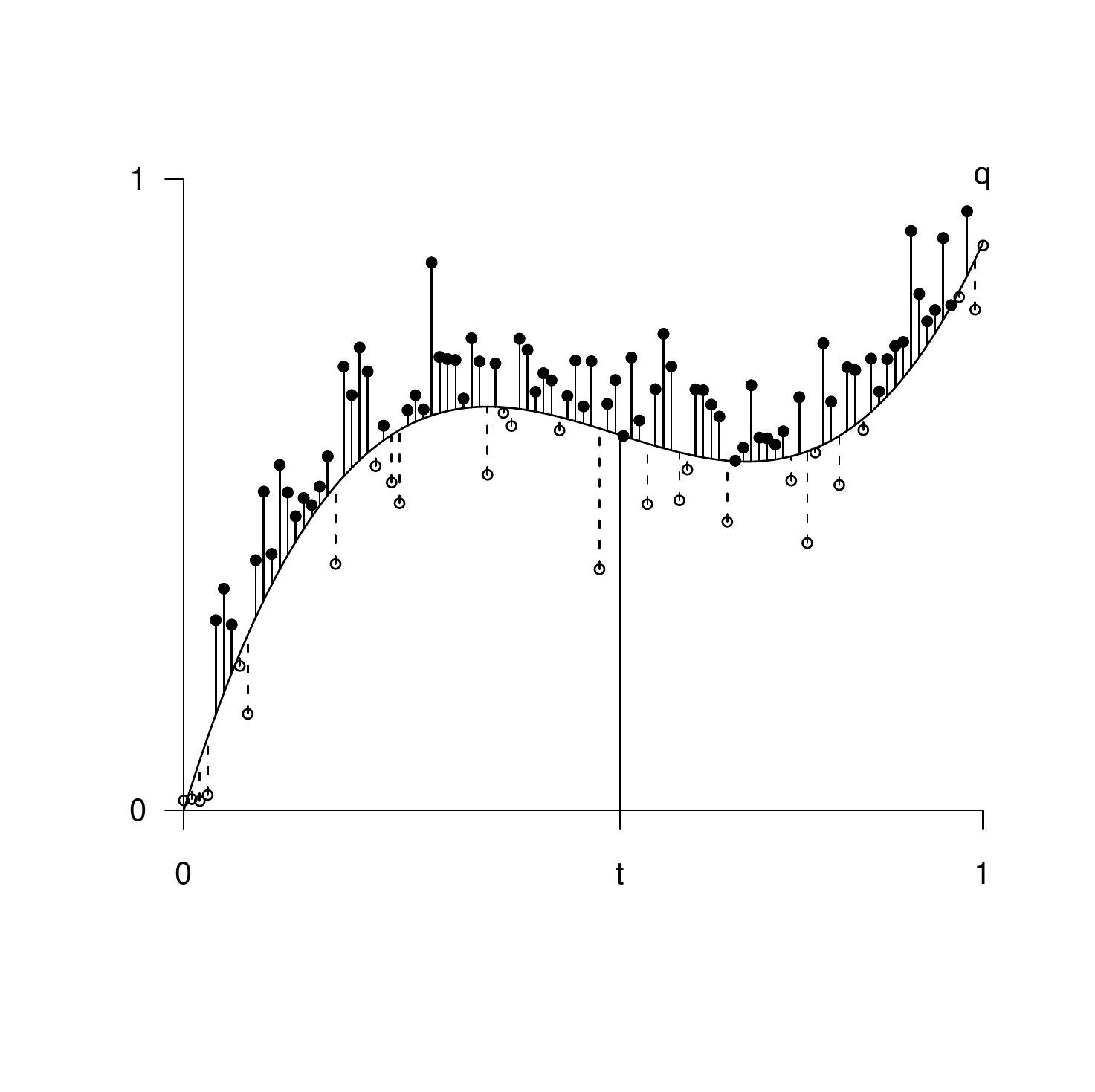,width=8.5cm} \epsfig{file=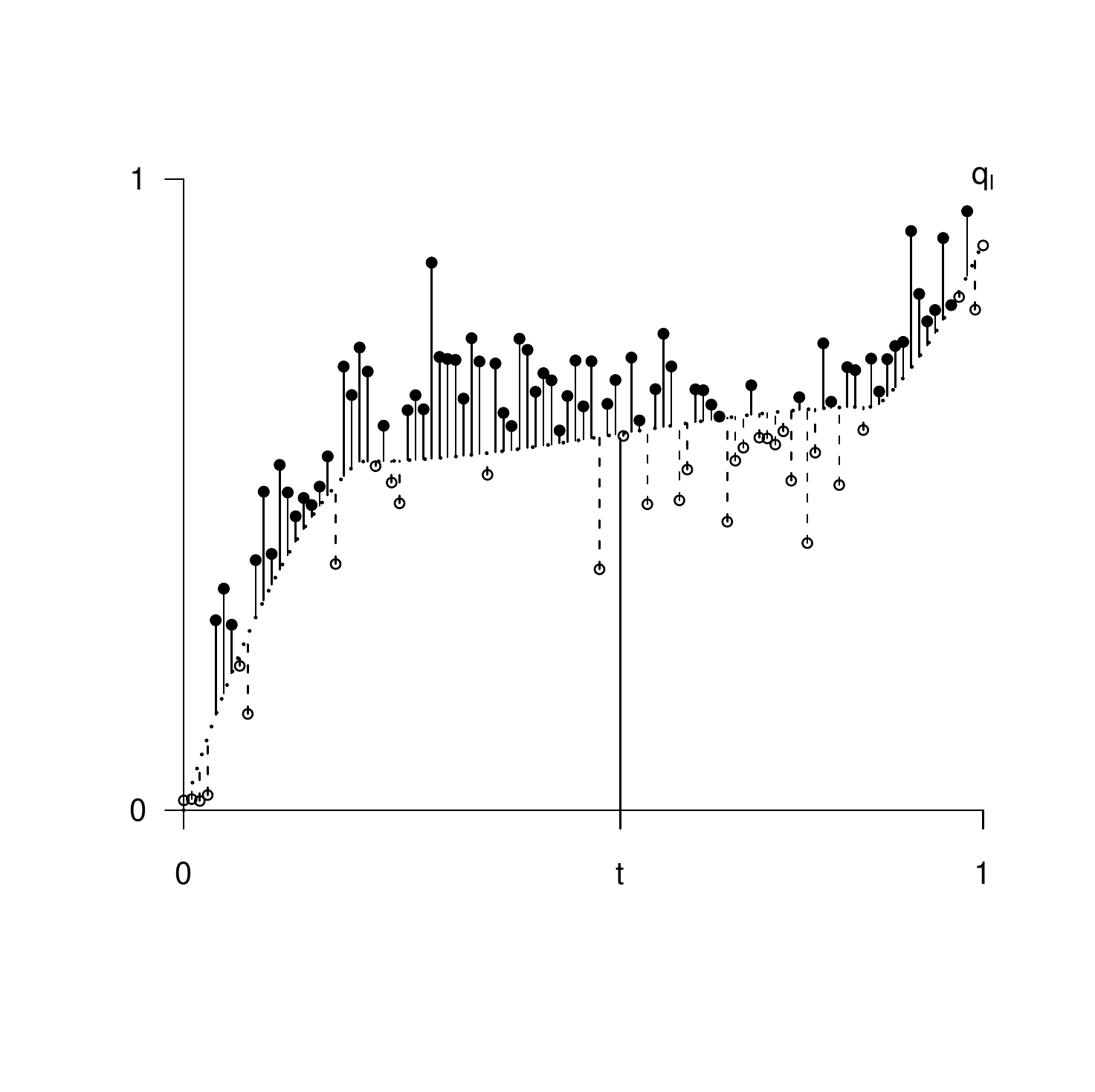,width=8.5cm}\vspace{-1.5cm}
\caption{\label{fig2} \sl Left part: True nonincreasing function $q_\tau$ for $\tau=0.25$ and errors for the sample shown in Figure \ref{fig1}. Right part: Increasing rearrangement  $q_{\tau,I}$  and pseudo-errors.
(Positive errors are marked by solid points and solid lines, negative errors marked by circles and dashed lines.)}
\end{figure}
If the true function $q_\tau$ is not increasing (e.g.~like in Figure \ref{fig1}) and we calculate the pseudo-errors from $q_{\tau,I}$, they are no longer identically distributed. This effect is demonstrated in Figure \ref{fig2} for a $\tau=0.25$-quantile curve. Consider for instance the interval $[t,1]$, where there are about $25\%$ negative errors (left part) and in comparison too many negative pseudo-errors (right part).
To detect such discrepancies from the null hypothesis, we estimate the pseudo-error distribution up to every $t\in [0,1]$ (i.\,e.\ for the covariate values $X_i\leq t$) and compare with what is expected under $H_0$. To this end recall the definition of $\hat F_{X,\eps,n}$ in (\ref{FXen}) and define $\hat F_{X,\eps_I,n}$ analogously, but using the constrained residuals $\hat\eps_{i,I}$, $i=1,\dots,n$. Analogously to (\ref{Sn-neu}) define the process
\begin{eqnarray}\label{Sn}
S_{n,I}(t,y) &=& \sqrt{n}\Big(\hat F_{X,\eps_I,n}(t,y)-\hat F_{X,\eps,n}(1-2h_n,y)\hat F_{X,\eps,n}(t,\infty)\Big)
\end{eqnarray}
for $y\in\mathbb{R}$, $t\in[2h_n,1-2h_n]$, and $S_{n,I}(t,y)=0$ for $y\in\mathbb{R}$, $t\in[0,2h_n)\cup (1-2h_n,1]$. 
For each fixed $t\in[0,1]$, $y\in\R$, for $h_n\to 0$ the statistic $n^{-1/2}S_{n,I}(t,y)$ consistently estimates the expectation
\begin{eqnarray*}
&&E[I\{\eps_{i,I}<y\}I\{X_i\leq t\}]-F_\eps ( y) F_{X}(t)\\
&=& E\Big[I\Big\{\eps_{i}<y+\frac{(q_{\tau,I}-q_{\tau})(X_i)}{s(X_i)}\Big\}I\{X_i\leq t\}\Big]-F_\eps(y) F_{X}(t) \label{equ1}. 
\end{eqnarray*}
Define a  Kolmogorov-Smirnov type statistic as $K_n=\sup_{y\in\mathbb{R},t\in[0,1]}|S_{n,I}(t,y)|$. Then $n^{-1/2}K_n$ estimates
$$K=\sup_{t\in[0,1],y\in\mathbb{R}}\left|\int_0^t\Big(F_\eps\Big(y+\frac{(q_{\tau,I}-q_{\tau})(x)}{s(x)}\Big)-F_\eps(y) \Big) f_{X}(x)\,dx\right|.$$
Note that under $H_0: q_{\tau,I}=q_{\tau}$ we have $K=0$. On the other hand, if $K=0$ then also 
\begin{eqnarray*}
\sup_{t\in[0,1]}\Big|\int_0^t\Big(F_\eps\Big(\frac{(q_{\tau,I}-q_{\tau})(x)}{s(x)}\Big)-F_\eps(0) \Big) f_{X}(x)\,dx\Big| &=& 0
\end{eqnarray*}
and from this it follows that $q_{\tau,I}=q_\tau$ is valid $F_X$-a.\,s.\ by the strict monotonicity of $F_\eps$. Thus under the alternative we have $K>0$ and $K_n$ converges to infinity. Define $c$ as the $(1-\alpha)$-quantile of the distribution of $\sup_{t\in[0,1],y\in\mathbb{R}}|S(t,y)|$ with $S$ from Theorem \ref{theo1}. Then the test that rejects $H_0$ for $K_n>c$ is consistent by the above argumentation and has asymptotic level $\alpha$ by the next theorem and an application of the continuous mapping theorem. 

\begin{theo}\label{theo1-mon} Under model (\ref{mod-het}) and assumptions \ref{as:k1}-\ref{as:Gs}, \ref{as:fx}-\ref{as:fbound} and \ref{as:bw}, under the null hypothesis $H_0$ and the assumption $\inf_{x\in[0,1]} q_\tau'(x)>0$ the process $S_{n,I}$ converges weakly in $\ell^\infty([0,1]\times \R)$ to the Gaussian process $S$ defined in Theorem \ref{theo1}. 
\end{theo}
 
The proof is given in Appendix \ref{app-A}. 

\begin{rem}\rm
Note that we use non-smooth monotone rearrangement estimators $\hat q_{\tau,I}$. Dette et al.\ (2006) and Birke and Dette (2008) consider smooth versions of the increasing rearrangements in the context of monotone mean regression. Corresponding increasing quantile curve estimators could be defined as 
\begin{eqnarray*}
\hat q_{\tau,I} (x) &=& \inf\Big\{z\in\mathbb{R}\;\Big|\;  \frac{1}{b_n} \int_0^1 \int^z_{-\infty} k \Bigl( \frac{{\hat q_\tau} (v) - u}{b_n}\Bigr) \,dudv \geq x\Big\}.
\end{eqnarray*}
 Under suitable assumptions on the kernel $k$ and bandwidths $b_n$ it can be shown that the same weak convergence as in Theorem \ref{theo1-mon} holds for $S_{n,I}$ based on this estimator. 
 $\blacksquare$
\end{rem}

For the application of the test for monotonicity we suggest a bootstrap version of the test analogously to the one considered in Section \ref{sec-asy}, but applying the increasing estimator to build new observations, i.\,e.\ $Y_i^*=\hat q_{\tau,I}(X_i)+\hat s(X_i)\eps_i^*$, $i=1,\dots,n$. We have the following theoretical result.

\begin{theo}\label{theo1-mon-boot} Under the assumptions of Theorem \ref{theo1-mon} and \ref{as:b1}--\ref{as:b2} the process $S_{n,I}^*$, conditionally on $\Y_n$, converges weakly in $\ell^\infty([0,1]\times \R)$ to the Gaussian process $S$ defined in Theorem \ref{theo1}, in probability. 
\end{theo}

The proof is given in Appendix \ref{app-boot}.
A consistent asymptotic level-$\alpha$ test is constructed as in Remark \ref{rem-bootKS}.

\begin{rem}\rm 
In the context of testing for monotonicity of mean regression curves Birke and Neumeyer (2013) based their tests on the observation that too many of the pseudo-errors are positive (see solid lines in Figure \ref{fig2}) on some subintervals of $[0,1]$ and too many are negative (see dashed lines) on other subintervals. Transferring this idea to the quantile regression model, one would consider a stochastic process
\begin{eqnarray*}
\tilde S_n(t,0) &=& \frac{1}{\sqrt{n}}\sum_{i=1}^n\Big(I\{\hat\eps_{i,I}\leq0\}I\{2h_n<X_i\leq t\}-\hat F_{X,\eps,n}(1-2h_n,0) I\{2h_n<X_i\leq t\}\Big)
\end{eqnarray*}
or alternatively (because $\hat F_{X,\eps,n}(1-2h_n,0)$ estimates the known $F_\eps(0)=\tau$)
\begin{eqnarray*}
R_n(t) &=& \frac{1}{\sqrt{n}}\sum_{i=1}^n\Big(I\{\hat\eps_{i,I}\leq0\}I\{X_i\leq t\}-\tau I\{X_i\leq t\}\Big)
\end{eqnarray*}
where $t\in [0,1]$.
For every $t\in[2h_n,1-2h_n]$ the processes count how many pseudo-residuals are positive up to covariates $\leq t$. This term is then centered with respect to the estimated expectation under $H_0$ and scaled with $n^{-1/2}$. 
However, as can be seen from Theorem \ref{theo1-mon} the limit is degenerate for $y=0$, and hence we have under $H_0$ that
\begin{eqnarray}\label{neg}
\sup_{t}|\tilde S_n(t,0)|=o_P(1).
\end{eqnarray} 
Also, $\sup_{t\in[0,1]}|R_n(t)|=o_P(1)$ can be shown analogously. 
Hence, no critical values can be obtained for the Kolmogorov-Smirnov test statistics, and those test statistics are not suitable for our testing purpose. To explain the negligibility (\ref{neg}) heuristically, consider the case $t=1$ (now ignoring the truncation of covariates for simplicity of explanation). Then, under $H_0$,  $n^{-1}\sum_{i=1}^nI\{\hat\eps_{i,I}\leq0\}$ estimates $F_\eps(0)=\tau$. But the information that $\eps_i$ has $\tau$-quantile zero was already applied to estimate the $\tau$-quantile function $q_\tau$.  Hence, one obtains $n^{-1}\sum_{i=1}^nI\{\hat\eps_{i,I}\leq0\}-\tau=o_P(n^{-1/2})$. This observation is in accordance to the fact that $n^{-1}\sum_{i=1}^n\hat\eps_i=o_P(n^{-1/2})$, when residuals are built from a mean regression model with centered errors [see M³ller et al.\ (2004) and Kiwitt et al.\ (2008)]. \\
Finally, consider the process 
\begin{eqnarray*}
\tilde S_n(1-2h_n,y) &=& \frac{1}{\sqrt{n}}\sum_{i=1}^n\Big(I\{\hat\eps_{i,I}\leq y\}I\{2h_n<X_i\leq 1-2h_n\}\\
&&{}\qquad -\hat F_{X,\eps,n}(1-2h_n,y) I\{2h_n<X_i\leq 1-2h_n\}\Big)
\end{eqnarray*}
i.\,e.\ the difference between the estimated distribution functions of pseudo-residuals $\hat\eps_{i,I}$ and unconstrained residuals $\hat\eps_i$ ($i=1,\dots,n$), respectively, 
scaled with $n^{1/2}$. An analogous process has been considered by Van Keilegom et al.\ (2008)
 for testing for parametric classes of mean regression functions. However, as can be seen from Theorem \ref{theo1-mon}, in our case of testing for monotonicity the limit again is degenerate, i.\,e.\ $\Var(S(1,y))=0$ for all $y$, and hence $\sup_{y\in\R}|\tilde S_n(1,y)|=o_P(1)$. Similar observations can be made when typical distance based tests from lack-of-fit literature [for instance $L^2$-tests or residual process based procedures by Hõrdle and Mammen (1993) and Stute (1997), respectively] are considered in the problem of testing monotonicity of regression function, see Birke and Neumeyer (2013). The reason is that under $H_0$ the unconstrained and constrained estimators, $\hat q_\tau$ and $\hat q_{\tau,I}$, typically are first order asymptotically equivalent. This for estimation purposes  very desirable property limits the possibilities to apply the estimator $\hat q_{\tau,I}$ for hypotheses testing. 
 $\blacksquare$
\end{rem}


\section{Simulation results}
\def\theequation{5.\arabic{equation}}
\setcounter{equation}{0}

In this section we show some simulation results for the bootstrap based tests introduced in this paper. If available we compare the results to already existing methods. Throughout the whole section we choose the bandwidths according to condition \ref{as:bw} as $d_n = 2(\hat \sigma^2/n)^{1/7}$, $h_n = (\hat \sigma^2/n)^{1/7}$, $b_n = \hat \sigma^2 (1/n)^{2/7}$ and $\hat \sigma^2$ is the difference estimator proposed in Rice (1984) [see Yu and Jones (1997) for a related approach]. The degree of the local polynomial estimators of location and scale [see equation (\ref{def-Fdach})] was chosen to be 3, the Kernel $K$ is the Gauss Kernel while $\kappa$ was chosen to be the Epanechnikov Kernel. The function $\Omega$ was defined through $\Omega(t) = \int_{-\infty}^t \omega(x)dx$ where $\omega(x) := (15/32)(3 - 10x^2 + 7x^4)I\{|x|\leq 1\}$, which is a kernel of order $4$ [see Gasser et al.\ (1985)]. For the choice of the distribution functions $G$ and $G_s$, we follow the procedure described in Dette and Volgushev (2008) who suggested a normal distribution such that the $5\%$ and $95\%$ quantiles coincide with  the corresponding empirical quantities of the sample $Y_1,...,Y_n$. Finally, the parameter $\alpha_n$ for generating the bootstrap residuals was chosen as $\alpha_n = 0.1n^{-1/4}\sqrt{2}\,\mbox{median}(|\hat\eps_1|,...,|\hat\eps_n|)$. All the results under $H_0$ are based on $1000$ simulation runs and $200$ bootstrap replications while the results under alternatives are simulated with $500$ simulation runs and $200$ bootstrap replications. 

\subsection{Testing for location and location-scale models} \label{sim:locscale}

The problem of testing the validity of location and location-scale models has previously been considered by Einmahl and Van Keilegom (2008a) and Neumeyer (2009b), and we therefore compare the properties of our test statistic with theirs. In testing the validity of location models (see Remark \ref{location-model}), we considered the following data generation processes
\bean
\mbox{(model 1)} && Y|X=x \sim (x - 0.5x^2) + \frac{(1+ax)^{1/2}}{10}\mathcal{N}(0,1), \quad X \sim U[0,1],
\\
\mbox{(model 2a)} && Y|X=x \sim (x - 0.5x^2) + \frac{1}{10}\Big( 1-\frac{1}{2c}\Big)^{1/2}t_{c}, \quad X \sim U[0,1],
\\
\mbox{(model 2b)} && Y|X=x \sim (x - 0.5x^2) + \frac{1}{10}\Big( 1-(cx)^{1/4}\Big)^{1/2}t_{2/(cx)^{1/4}}, \quad X \sim U[0,1],
\\
\mbox{(model 3)} && Y|X=x,U=u \sim (x - 0.5x^2) + \Big(U - 0.5 - \frac{b}{6}(2x-1)\Big), \quad (X,U) \sim C(b).
\eean
Note that model 1 with parameter $a=0$, model 2a with arbitrary parameter $c$, and model 3 with parameter $b=0$ correspond to a location model, while models 1, 2b and 3 with parameters $a,b,c\neq 0$ describe models that are not of this type.
Here $t_c$ denotes a $t-$distribution with $c$ degrees of freedom ($c$ not necessarily integer) and models 1 and 2b have also been considered by Einmahl and Van Keilegom (2008a).
Model 3 is from Neumeyer (2009b) and $(X,U)\sim C(b)$ are generated as follows. Let $X,V,W$ be independent $U[0,1]$-distributed random variables and define $U=\min(V,W/(b(1-2X)))$ if $X\leq \frac 12$, and $U=\max(V, 1+ W/(b(1-2X)))$ else. 
Note that this data generation produces observations from the Farlie-Gumbel-Morgenstern copula if the parameter $b$ is between $-1$ and $1$.  

Simulation results under the null are summarized in Table \ref{sizeloc}.
As we can see, both the Kolmogorov-Smirnov (KS) and the CramÚr-von Mises (CvM) bootstrap versions of the test hold the level quite well in all models considered and both for $n=100$ and $n=200$ observations.

Next, we take a look at the power properties of the tests in models 1, 2a and 3. The rejection probabilities are reported in Table \ref{sizeloc1}, Table \ref{sizeloc2} and Table \ref{sizeloc3}, respectively. For the sake of comparison, we have also included the results reported in Neumeyer (2009b) (noted N in the tables) and Einmahl and Van Keilegom (2008a) (noted EVK in the tables), where available. Note that Neumeyer (2009b) considers several bandwidth parameters, while Einmahl and Van Keilegom (2008a) consider various types of test statistics (KS, CvM and Anderson-Darling) and two types of tests (difference and estimated residuals). We have included the \textit{best} values of all the possible tests in Neumeyer (2009b) and Einmahl and Van Keilegom (2008a). Note that this does not correspond to a practical data-driven test since typically the best test is unknown. 

An inspection of Table \ref{sizeloc1} and Table \ref{sizeloc2} reveals that the tests of Neumeyer (2009b) and Einmahl and Van Keilegom (2008a) perform better for normal errors (Table \ref{sizeloc1}), while our test seems to perform better for $t$ errors (Table \ref{sizeloc2}). This corresponds to intuition since for normal errors the mean provides an optimal estimator of location, while for heavier tailed distributions the median has an advantage. Additionally, we see that in almost all cases the CvM test outperforms the KS test. In model 3, the test of Neumeyer (2009b) performs better than the tests proposed here, with significantly higher power for $b=1,2$ and $n=200$. The CvM version again has somewhat higher power than the KS version of the test. Overall, we can conclude that the newly proposed testing procedures show a competitive performance and can be particularly recommended for error distributions with heavier tails. The CvM test seems to always be preferable.\\

\begin{center}
\textbf{Please insert Tables \ref{sizeloc}, \ref{sizeloc1}, \ref{sizeloc2} and \ref{sizeloc3} here}
\end{center}

To evaluate the test for location-scale models, we considered the following settings 
\bean
\mbox{(model $1_h$)} && Y|X=x \sim (x - 0.5x^2) + \frac{2+x}{10}\mathcal{N}(0,1), \quad X \sim U[0,1],
\\
\mbox{(model $2a_h$)} && Y|X=x \sim (x - 0.5x^2) + \frac{2+x}{10}\Big( 1-\frac{1}{2c}\Big)^{1/2}t_{c}, \quad X \sim U[0,1],
\\
\mbox{(model $2b_h$)} && Y|X=x \sim (x - 0.5x^2) + \frac{2+x}{10}\Big( 1-(cx)^{1/4}\Big)^{1/2}t_{2/(cx)^{1/4}}, \quad X \sim U[0,1],
\\
\mbox{(model $3_h$)} && Y|(X,U)=(x,u) \sim (x - 0.5x^2) + \frac{2+x}{10}\Big(U - 0.5 - b(2x-1)\Big)
\\ && \quad\quad\quad\quad\quad(X,U) \sim C(b),
\eean
Models $1_h$ and $2b_h$ have also been considered in Einmahl and Van Keilegom (2008a), while model $3_h$ is from Neumeyer (2009b). Simulation results corresponding to different null models are collected in Table \ref{sizelocscal}. We observe that in all three models both the KS and the CvM test hold their level quite well for all sample sizes, with both tests being slightly conservative for $n=50$ and in model $1_h$.

The power against alternatives in model $2b_h$ and $3_h$ is investigated in Table \ref{sizelocscal2} and Table \ref{sizelocscal3}, respectively. From Table \ref{sizelocscal3}, we see that the CvM version of the proposed test has higher (sometimes significantly so) power than the test of Neumeyer (2009b). One surprising fact is that the power of the test of Neumeyer (2009b) decreases for large values of $b$, while the power of our test continues to increase. This might be explained by the fact that for larger values of $b$ (in particular for $b=5$), the variance of the residuals is extremely small, which probably leads to an instability of variance estimation. 

Inspecting Table \ref{sizelocscal2}, we see that the situation differs dramatically from the results in the homoscedastic model $2b$. In this particular setting, the tests proposed in this paper have no power for $n=50,n=100$, even for the most extreme setting $b=1$ (only this setting is shown here since for smaller values of $b$ the test also does not have any power). The test of Einmahl and Van Keilegom (2008a) has less power than in the homoscedastic case, but is still able to detect that this model corresponds to the alternative. An intuitive explanation of those differences is that Einmahl and Van Keilegom (2008a) scale their residuals to have the same variances while our residuals are scaled to have the same median absolute deviation (note that the mean and median of a t-distribution coincide provided that the mean exists). Under various alternative distributions, this leads to different power curves for the location-scale test. This difference is particularly extreme in the case of $t$-distributions. To illustrate this fact, recall in models which are not of location-scale structure, $n^{-1/2}S_n(t,y)$ converges in probability to $P(\eps_i^{ad}\leq y, X_i\leq t) - F_{\eps^{ad}}(y)F_X(t)$, see Remark \ref{rem-alt}. Here, the residuals $\eps^{ad}_i$ are defined as $(Y_i-F_Y^{-1}(\tau|X_i))/s(X_i)$ with $s(x)$ denoting the conditional median absolute deviation of $Y_i-F_Y^{-1}(\tau|X_i)$ given $X_i = x$. A similar result holds for the residuals in EVK which take the form $\eps_i^\sigma := (Y_i - m(X_i))/\sigma(X_i)$ where $\sigma^2$ denotes the conditional variance. One thus might expect that computing the quantities $K^{ad} := \sup_{t,y}|P(\eps_i^{ad}\leq y, X_i\leq t) - F_{\eps^{ad}}(y)F_X(t)|$ and $K^{\sigma} := \sup_{t,y}|P(\eps_i^\sigma\leq y,X_i\leq t) - F_{\eps^\sigma}(y)F_X(t)|$ will give some insights into the power properties of the KS test for residuals that are scaled in different ways. Indeed, numerical computations show that $K^{\sigma}/K^{ad} \approx 4.5$ which explains the large difference in power (note that the power for EVK reported in Table \ref{sizelocscal2} is in fact the power of their Anderson-Darling test, the power of the KS test in EVK is lower). For a corresponding version of the CvM distance the ratio is roughly ten. We suspect that using a different scaling for the residuals would improve the power of the test in this particular model. However, since the optimal scaling depends on the underlying distribution of the residuals which is typically unknown, it seems difficult to implement an optimal scaling in practice. We leave this interesting question to future research.\\
Note that we do not present simulation results for the models with $\chi^2$-distributed errors considered by Einmahl and Van Keilegom (2008a) and Neumeyer (2009b) for power simulations. The reason is that for error distributions that are $\chi^2_b$ with $b < 2$, tests based on residuals do not hold their level, and the power characteristics described in the aforementioned papers are a consequence of this fact. The intuitive reason for this fact is that weak convergence of the residual process requires the errors to have a uniformly bounded density, which is not the case for chi-square distributions with degrees of freedom less than two. This phenomenon is not related to non-parametric estimation of the location function and can already be observed in a simple linear model.

\begin{center}
\textbf{Please insert Tables \ref{sizelocscal}, \ref{sizelocscal2} and \ref{sizelocscal3} here}
\end{center}

\subsection{Testing for monotonicity of quantile curves in a location-scale setting} \label{sec:monsim}

\begin{center}
\textbf{Please insert Figure \ref{fig3} here}
\end{center}  

Next, we considered the test for monotonicity of quantile curves that is introduced in Section \ref{sec-mon}. Here, we simulated the following two models that are both of location-scale type
\bean
(\mbox{model 4}) &\quad& Y|X=x \sim 1 + x - \beta e^{-50(x-0.5)^2} + 0.2\mathcal{N}(0,1), \quad X\sim U[0,1]
\\ 
(\mbox{model 5}) &\quad& Y|X=x \sim \frac{x}{2} + 2(0.1 - (x-0.5)^2)\mathcal{N}(0,1), \quad X\sim U[0,1].
\eean 
The results for models 4 and 5 are reported in Table \ref{sizepower4} and Table \ref{sizepower5}, respectively. In model 4, all quantile curves are parallel and so all quantile curves have a similar monotonicity behavior. In particular, the parameter value $\beta = 0$ corresponds to strictly increasing quantile curves, for $\beta = 0.15$ the curves have a flat spot, and for $\beta > 0.15$ the curves have a small decreasing bump that gets larger for larger values of $\beta$. The median curves for different values of $\beta$ are depicted in Figure \ref{fig3}, and the $25\%$ quantile curves are parallel to the median curves with exactly the same shape. We performed the tests for two different quantile  curves ($\tau=0.25$ and $\tau=0.5$) and see that in both cases the test has a slowly increasing power for increasing values of $\beta$ and sample size. The case $\beta = 0.45$ is already recognized as alternative for $n=50$, while for $\beta = 0.25$ the test only starts to show some power for $n=200$. Note also that for very large sample sizes, even the flat function corresponding to $\beta = 0.15$ should be recognized as alternative since all the results under $H_0$ require that the quantile curves are strictly increasing. However, with a sample of size $n=200$ this effect is not visible in the simulations. 
\\
In model 5, the median is a strictly increasing function while the outer quantile curves are not increasing. In Table \ref{sizepower5}, we report the simulation results for three different quantile values ($\tau=0.25, \tau=0.5$ and $\tau=0.75$) and two sample sizes $n=50, 100, 200$. For $n=50$, the observed rejection probabilities are slightly above the nominal critical values (for $\tau=0.5$), and the cases $\tau=0.25$ and $\tau=0.75$ are recognized as alternatives. For $n=100, 200$, the test holds its level for $\tau=0.5$ and also shows a slow increase in power at the other quantiles. The increase is not really significant when going from $n=50$ to $n=100$ for $\tau=.25$ and not present for $\tau=.75$. For $n=200$, the test clearly has more power compared to $n =50$. Overall, we can conclude that the proposed test shows a satisfactory behavior.\\

\begin{center}
\textbf{Please insert Tables \ref{sizepower4} and \ref{sizepower5} here}
\end{center}

\section{Conclusion}

The paper at hand considered location-scale models in the context of nonparametric quantile regression. For the first time a test for model validity was investigated.
It is based on the empirical independence process of covariates and residuals built from nonparametric estimators for the location and scale functions. The process converges weakly to a Gaussian process. A bootstrap version of the test was investigated in theory and by means of a simulation study. 
The theoretical results open a new toolbox to test for various model hypotheses in location-scale quantile models. As example we considered in detail the testing for monotonicity of a conditional quantile function in theory as well as in simulations. Similarly other structural assumptions on the  location or the scale function can be tested. All weak convergence results are proved in the appendix and supplementary material in a detailed manner. A small simulation study demonstrated that the proposed method works well.\\

\newpage

\begin{appendix}

\section{Proof of weak convergence results}\label{app-A}
\def\theequation{A.\arabic{equation}}
\setcounter{equation}{0}

Before beginning with the proof, we give a brief overview of the results. The proofs of the main results (Theorem \ref{theo1}, Corollary \ref{cor-KS-CvM} and Theorem \ref{theo1-mon}) and the bootstrap versions (Theorems \ref{theo1-boot} and \ref{theo1-mon-boot}) are contained in Appendixes \ref{app-A} and \ref{app-boot}, respectively. Technical details needed in the proofs of those results can be found in the supplementary material in Appendix \ref{app-tech}. Finally, Appendix \ref{app-quant} in the supplement contains basic results on linearized versions and differentiability of the quantile estimator $\hat q_\tau$, scale estimator $\hat s$ and the corresponding bootstrap versions, while Appendix \ref{sec:tec} contains additional technical details.\\
\\
{\bf Proof of Theorem \ref{theo1}.}
For the numerator $\bar F_{X,\eps,n}(t,y)=\hat F_{X,\eps,n}(t,y)(\hat F_{X,n}(1-2h_n)-\hat F_{X,n}(2h_n))$ of the joint empirical distribution function defined in (\ref{FXen}) we have
\begin{eqnarray*}
\bar F_{X,\eps,n}(t,y) &=& \frac{1}{n}\sum_{i=1}^n I\Big\{\eps_i\leq y\frac{\hat s(X_i)}{s(X_i)}+\frac{\hat q_{\tau}(X_i)-q_\tau (X_i)}{s(X_i)}\Big\}I\{2h_n< X_i\leq t\}.
\end{eqnarray*} 
Note that in Lemma \ref{lem:proclin} in the supplement it is shown that without changing the asymptotic distribution of the process the residuals $\hat\eps_i$ can be replaced by their versions obtained from linearized estimators $\hat q_{\tau,L}$, $\hat s_L$ instead of $\hat q_\tau$, $\hat s$ (see Appendix \ref{app-quant} for the definitions). Thus we have
\begin{eqnarray*}
\bar F_{X,\eps,n}(t,y) &=& \frac{1}{n}\sum_{i=1}^n I\Big\{\eps_i\leq y\frac{\hat s_L(X_i)}{s(X_i)}+\frac{\hat q_{\tau,L}(X_i)-q_\tau (X_i)}{s(X_i)}\Big\}I\{2h_n< X_i\leq t\}+o_P(\frac{1}{\sqrt{n}}).
\end{eqnarray*} From this we obtain the expansion 
\begin{eqnarray}\nb
\bar F_{X,\eps,n}(t,y)&=& \frac{1}{n}\sum_{i=1}^n I\{\eps_i\leq y\}I\{2h_n< X_i\leq t\}\\
\label{avk}&&{}+\int_{2h_n}^{1-2h_n} \Big(F_\eps\Big(y\frac{\hat s_L(x)}{s(x)}+\frac{\hat q_{\tau,L}(x)-q_\tau (x)}{ s(x)}
\Big)-F_\eps(y)\Big)I\{x\leq t\}f_X(x)\,dx\\
\nb
&&{} +o_P(\frac{1}{\sqrt{n}})
\end{eqnarray}
uniformly with respect to $t\in[2h_n,1-2h_n]$ and $y\in\R$ by the following argumentation. Consider the empirical process
\begin{eqnarray*}
G_n(\varphi)&=&\frac{1}{\sqrt{n}}\sum_{i=1}^n \Big( \varphi(X_i,\eps_i)-E[\varphi(X_i,\eps_i)] \Big),\quad \varphi\in\mathcal{F},
\end{eqnarray*}
indexed by the following class of functions,
\begin{eqnarray*}
\mathcal{F} &=& \Big\{(X,\eps)\mapsto I\{\eps\leq yd_2(X)+d_1(X)\}I\{h< X\}I\{X\leq t\}-I\{{\eps}\leq y\}I\{h< X\}I\{X\leq t\}\; \\ && \Big|\; y\in\mathbb{R},h,t\in [0,1], d_1\in C_1^{1+\delta}([0,1]), d_2\in \tilde C_2^{1+\delta}([0,1])\Big\},
\end{eqnarray*}
for some arbitrary $\delta\in (0,1)$, 
 where the function class $C_c^{1+\delta}([0,1])$ is defined as the set of differentiable functions $g: [0,1]\to\R$ with derivatives $g'$ such that
\begin{eqnarray*}
\max\Big\{\sup_{x\in [0,1]}|g(x)|,\sup_{x\in [0,1]}|g'(x)|\Big\}+\sup_{x,z\in [0,1]}\frac{|g'(x)-g'(z)|}{|x-z|^\delta} 
&\leq& c
\end{eqnarray*}
[see van der Vaart and Wellner (1996, p.\ 154)].
We further by slight abuse of notation define the subset $\tilde{C}_{2}^{1+\delta}([0,1])$ of $C_1^{1+\delta}([0,1])$ by the additional constraint $\inf _{x\in [0,1]}g(x)\geq 1/2$.
Now $\mathcal{F}$ is a product of the uniformly bounded Donsker classes $\{(X,\eps)\mapsto I\{h< X\}I\{X\leq t\}|h,t\in [0,1]\}$ and $\{(X,\eps)\mapsto I\{\eps\leq yd_2(X)+d_1(X)\}-I\{{\eps}\leq y\}| y\in\mathbb{R}, d_1\in C_1^{1+\delta}([0,1]), d_2\in \tilde C_2^{1+\delta}([0,1])\}$ [the Donsker property for the second class is shown in Lemma 1 by Akritas and Van Keilegom (2001)] and is therefore Donsker as well (Ex.\ 2.10.8, van der Vaart and Wellner (1996), p.\ 192). The remaining part of the proof for equality (\ref{avk}) follows exactly the lines of the end of the proof of Lemma 1, Akritas and Van Keilegom (2001), p.\ 567, using the inequality
\begin{eqnarray*}
&&\Var\Big(I\{\eps_1\leq yd_2(X_1)+d_1(X_1)\}I\{h< X_1\}I\{X_1\leq s\}-I\{{\eps_1}\leq y\}I\{h< X_1\}I\{X_1\leq s\}\Big) \\
&\leq&
E\Big[\Big(I\{\eps_1\leq yd_2(X_1)+d_1(X_1)\}-I\{{\eps_1}\leq y\}\Big)^2\Big].
\end{eqnarray*} 
Here one also needs $\hat s_L/s\in \tilde C_2^{1+\delta}([0,1])$, $(\hat q_{\tau,L}-q_\tau)/s\in C_1^{1+\delta}([0,1])$ with probability converging to one, which follows from uniform consistency results in Lemma \ref{lem:unifratesquant} in the supplement.
 For $\varphi=\varphi_{h,t,y,d_1,d_2}$ we obtain
$$\sup_{y\in\mathbb{R},\atop t\in[2h_n,1-2h_n]}\Big|G_n\Big(\varphi_{2h_n,t,y,\frac{\hat q_{\tau,L}-q_\tau}{s},\frac{\hat s_L}{s}}\Big)\Big|=o_P(1)$$
and thus (\ref{avk}).

Further, by a Taylor expansion we obtain from (\ref{avk}) together with assumption \ref{as:eps} that 
\begin{eqnarray*}
\bar F_{X,\eps,n}(t,y) 
&=& \frac{1}{n}\sum_{i=1}^nI\{\eps_i\leq y\}I\{2h_n< X_i\leq t\} + yf_\eps(y)\int_{2h_n}^{1-2h_n}\frac{\hat s_L(x)-s(x)}{s(x)}I\{x\leq t\}f_X(x)\,dx\\
&&{}
+f_\eps(y)\int_{2h_n}^{1-2h_n}\frac{\hat q_{\tau,L}(x)-q_\tau(x)}{s(x)}I\{x\leq t\}f_X(x)\,dx
 +o_P(\frac{1}{\sqrt{n}})\\
\end{eqnarray*}
uniformly with respect to $t\in[2h_n,1-2h_n]$ and $y\in\R$.
In  Lemma \ref{lem:B2boot} in the supplementary material expansions of the integrals in this decomposition are derived and  it follows that
\begin{eqnarray}
&&\bar F_{X,\eps,n}(t,y) \label{netto}\\
\nonumber&=& \frac{1}{n}\sum_{i=1}^nI\{\eps_i\leq y\}I\{2h_n< X_i\leq t\}
-\phi(y) \frac{1}{n}\sum_{i=1}^n(I\{\eps_i\leq 0\}-\tau)I\{2h_n< X_i\leq t\}\\
&&{}-\psi(y)\frac{1}{n}\sum_{i=1}^n\Big(I\{|\eps_i|\leq 1\}-\frac{1}{2}\Big)I\{2h_n< X_i\leq t\}
 +o_P(\frac{1}{\sqrt{n}}) ,\nonumber
\end{eqnarray}
where $\phi$ and $\psi$ are defined in the assertion of the theorem. 
Thus noting that $\hat F_{X,n}(1-2h_n)-\hat F_{X,n}(2h_n)=F_{X}(1-2h_n)-F_{X}(2h_n)+o_P(1)=1+o_P(1)$, from the definition (\ref{Sn-neu}) we obtain by Slutsky's lemma that
\begin{eqnarray*}
 S_n(t,y)&=& 
\frac{1}{\sqrt{n}}\sum_{i=1}^n\Big( I\{\eps_i\leq y\}-F_\eps(y)-\phi(y) (I\{\eps_i\leq 0\}-\tau)-\psi(y)\Big(I\{|\eps_i|\leq 1\}-\frac{1}{2}\Big)\Big)\\
&&\qquad {}\times\Big(I\{2h_n< X_i\leq t\}-I\{2h_n< X_i\leq 1-2h_n\}\frac{\hat F_{X,n}(t)-\hat F_{X,n}(2h_n)}{\hat F_{X,n}(1-2h_n)-\hat F_{X,n}(2h_n)}\Big)\\
&&{}  +o_P(1) .
\end{eqnarray*}
uniformly with respect to $t\in[2h_n,1-2h_n]$ and $y\in\R$. Note that the dominating part of this process vanishes in the boundary points $t=2h_n$ and $t=1-2h_n$. Further, from  $\hat F_{X,n}(t)=F_X(t)+O_P(n^{-1/2})$ uniformly in $t\in[0,1]$ and $F_X(2h_n)\to 0$, $F_X(1-2h_n)\to 1$ we have
\begin{eqnarray*}
 S_n(t,y)&=& 
S_{n,1}(t,y) +o_P(1) ,
\end{eqnarray*}
uniformly with respect to $t\in[0,1]$, $y\in\R$, where $S_{n,1}(t,y)=0$ for $t\in[0,2h_n)\cup (1-2h_n,1]$ and 
\begin{eqnarray*}
 S_{n,1}(t,y)&=& 
\frac{1}{\sqrt{n}}\sum_{i=1}^ng(\eps_i,y)\Big(I\{2h_n< X_i\leq t\}-I\{2h_n< X_i\leq 1-2h_n\}F_X(t)\Big)
\end{eqnarray*}
for $t\in[2h_n,1-2h_n]$ and $y\in\R$,
where $g(\eps_i,y)=I\{\eps_i\leq y\}-F_\eps(y)-\phi(y) (I\{\eps_i\leq 0\}-\tau)-\psi(y)(I\{|\eps_i|\leq 1\}-\frac{1}{2})$ is centered and independent of $X_i$. 
The first assertion of the theorem now follows if we show that for 
\begin{eqnarray*}
 S_{n,2}(t,y)&=& 
\frac{1}{\sqrt{n}}\sum_{i=1}^ng(\eps_i,y)\Big(I\{X_i\leq t\}-F_X(t)\Big)\, ,\quad t\in[0,1],y\in\R,
\end{eqnarray*}
we have 
$\sup_{t\in[0,1],y\in\R}|S_{n,1}(t,y)-S_{n,2}(t,y)|=o_P(1),$
which is equivalent to 
\begin{eqnarray}\label{bo1}
&&\sup_{t\in[2h_n,1-2h_n],y\in\R}|S_{n,1}(t,y)-S_{n,2}(t,y)|\;=\;o_P(1)
\end{eqnarray}

together with 
\begin{eqnarray}
&&\sup_{t\in[0,2h_n)\cup (1-2h_n,1],y\in\R}|S_{n,2}(t,y)|\;=\;o_P(1)\label{bo2}.
\end{eqnarray}
We will only show (\ref{bo1}); (\ref{bo2}) follows by similar arguments. 
Note that $S_{n,1}(t,y)-S_{n,2}(t,y)=G_n(h_n,t,y)$  for $t\in[2h_n,1-2h_n]$, $y\in\R$,
where the process
$$G_n(h,t,y)=\frac{-1}{\sqrt{n}}\sum_{i=1}^n g(\eps_i,y)(I\{X_i\leq t\}-F_X(t))I\{X_i\in [0,2h)\cup (1-2h,1]\}$$
indexed in $h\in[0,\frac{1}{4}]$, $t\in[0,1]$, $y\in\R$, converges weakly to a centered Gaussian process $G$  with asymptotic variance 
\begin{eqnarray*}
\Var(G(h,t,y))&=&E[g^2(\eps_1,y)]\Big((F_X(t\wedge 2h)+F_X(t)-F_X(t\wedge (1-2h)))(1-2F_X(t))\\
&&{}+F_X^2(t)(F_X(2h)+1-F_X(1-2h))\Big).
\end{eqnarray*}
For $h=h_n\to 0$ this asymptotic variance vanishes uniformly with respect to $y$ and $t$. From asymptotic equicontinuity of $G_n$ (confer van der Vaart and Wellner, 1996, p.\  89/90), using the asymptotic variance as semi-metric, with $G_n(0,t,y)\equiv 0$ it follows that $\sup_{t,y}|G_n(h_n,t,y)|=o_P(1)$ and thus (\ref{bo1}).

Hence, we have shown the first assertion of the theorem, i.\,e.\ $S_n=S_{n,2}+o_P(1)$ uniformly.
Weak convergence of $S_{n,2}$ (and thus of $S_n$) to a centered Gaussian process with the asserted covariance structure follows by standard arguments. 
\hfill $\Box$\\
\\

{\bf Proof of Corollary \ref{cor-KS-CvM}.}
The asymptotic distribution of $K_n$ directly follows from Theorem \ref{theo1} and the continuous mapping theorem. From those theorems also follows that
$$\tilde C_n = \int_{\mathbb{R}}\int_{[0,1]} S_n^2(t,y) \,F_{X}(dt)\,F_{\eps}(dy)
$$
converges in distribution to the desired limit. It therefore remains to show that $C_n-\tilde C_n=o_P(1)$. To this end denote 
\begin{eqnarray*}
\tilde C_n^{(1)} = \int_{\mathbb{R}}\int_{[0,1]} S_n^2(F_X^{-1}(\hat F_{X,n}(t)),F_\eps^{-1}(\hat F_{\eps,n}(y))) \,\hat F_{X,n}(dt)\,\hat F_{\eps,n}(dy)
\end{eqnarray*}
and let $\varrho_n$ be some sequence specified later with $\varrho_n\to\infty$ for $n\to\infty$. Then
\begin{eqnarray*}
|C_n-\tilde C_n^{(1)}|&\leq& \Big|\int_{[-\varrho_n,\varrho_n]}\int_{[0,1]} \Big( S_n^2(t,y)-S_n^2(F_X^{-1}(\hat F_{X,n}(t)),F_\eps^{-1}(\hat F_{\eps,n}(y)))\Big)\,\hat F_{X,n}(dt)\,\hat F_{\eps,n}(dy)\Big|\\
&&{}+2\sup_{t,y}|S_n^2(t,y)|\int_{\R\setminus[-\varrho_n,\varrho_n]}\hat F_{\eps,n}(dy).
\end{eqnarray*}
The second term on the right hand side is $O_P(1)(1-\hat F_{\eps,n}(\varrho_n)+\hat F_{\eps,n}(-\varrho_n))=o_P(1)$ due to the results from Theorem \ref{theo1} and because $\varrho_n\to\infty$ and $\hat F_{\eps,n}$ converges to $F_\eps$ uniformly in probability (this follows from the proof of Theorem \ref{theo1}). The first term on the right hand side can further be bounded by
 \begin{eqnarray*}
2\sup_{t,y}|S_n(t,y)| \sup_{t\in[0,1]\atop y\in [-\varrho_n,\varrho_n]}\Big|S_n(t,y)-S_n(F_X^{-1}(\hat F_{X,n}(t)),F_\eps^{-1}(\hat F_{\eps,n}(y)))\Big|.
\end{eqnarray*}
From Theorem \ref{theo1} it follows that the process $S_n$ is asymptotically stochastic equicontinuous such that we obtain the desired rate $o_P(1)$ from
$$\sup_{t\in[0,1]|}|t-F_X^{-1}(\hat F_{X,n}(t))|\leq\sup_{\xi\in[0,1]}\frac{1}{f_X(\xi)}\sup_{t\in[0,1]}|\hat F_{X,n}(t)-F_X(t)| =o_P(1)$$
by assumption \ref{as:fx}
and
$$\sup_{y\in [-\varrho_n,\varrho_n]}|y-F_\eps^{-1}(\hat F_{\eps,n}(y)))|\leq \sup_{y\in[-\varrho_n,\varrho_n]}\sup_{\zeta\mbox{ \tiny{between} }\atop F_\eps(y)\mbox{\tiny{ and} }\hat F_{\eps,n}(y)}\frac{1}{f_\eps(F_\eps^{-1}(\zeta))} \sup_{y\in\R}|\hat F_{\eps,n}(y)-F_\eps(y)|=o_P(1).
$$
The latter rate follows because $\sup_{y\in\R}|\hat F_{\eps,n}(y)-F_\eps(y)|=O_P(n^{-1/2})$ (which can be deduced by $\hat F_{\eps,n}(\cdot)=\overline{F}_{X,\eps,n}(1-2h_n,\cdot)/(\hat F_{X,n}(1-2h_n)-\hat F_{X,n}(2h_n))$ and (\ref{netto}) in the proof of Theorem \ref{theo1}) if we choose a sequence $\varrho_n$ such that
$n^{1/2}\inf_{y\in[-2\varrho_n,2\varrho_n]}f_\eps(y)\to\infty$ for $n\to\infty$. This is possible by assumption \ref{as:eps}.
\\
We have shown $C_n-\tilde C_n^{(1)}=o_P(1)$ and it remains to show that $\tilde C_n-\tilde C_n^{(1)}=o_P(1)$. To this end, note that almost surely
\begin{eqnarray*}
\tilde C_n-\tilde C_n^{(1)}&=& \int_{[0,1]}\int_{[0,1]} S_n^2(F_X^{-1}(s),F_\eps^{-1}(z))\,ds\,dz-\frac{1}{n^2}\sum_{i=1}^n\sum_{j=1}^n S_n^2(F_X^{-1}(\textstyle{\frac in}),F_\eps^{-1}(\textstyle{\frac jn}))\\
&=& \sum_{i=1}^n\sum_{j=1}^n\int_{[\frac{i-1}{n},\frac in)}\int_{[\frac{j-1}{n},\frac jn)}\Big( S_n^2(F_X^{-1}(s),F_\eps^{-1}(z))-S_n^2(F_X^{-1}(\textstyle{\frac in}),F_\eps^{-1}(\textstyle{\frac jn}))\Big)\,ds\,dz.
\end{eqnarray*}
We decompose the second sum into $\sum_{j=1}^{j_n}\dots+\sum_{j=J_n+1}^n\dots+\sum_{j=j_n+1}^{J_n}\dots$ for sequences of integers with $1\leq j_n<J_n\leq n$ and $j_n/n\to 0$, $J_n/n\to 1$ for $n\to\infty$. We obtain
\begin{eqnarray*}
|\tilde C_n-\tilde C_n^{(1)}|&\leq&
2 \frac{j_n+n-J_n}{n}\sup_{t,y}|S_n^2(t,y)| \\
&&{}
+2\sup_{t,y}|S_n(t,y)|\sup_{|s-u|\leq\frac 1n\atop s,u\in[0,1]}\sup_{|z-v|\leq \frac 1n\atop z,v\in \left[\frac{j_n}{n},\frac{J_n}{n}\right]}|S_n(F_X^{-1}(s),F_\eps^{-1}(z))-S_n(F_X^{-1}(u),F_\eps^{-1}(v))|.
\end{eqnarray*}
By asymptotic stochastic equicontinuity of $S_n$ this converges to zero in probability if 
$$\sup_{|s-u|\leq\frac 1n\atop s,u\in[0,1]}|F_X^{-1}(s)-F_X^{-1}(u)|\to 0$$
which follows from assumption \ref{as:fx} and the mean value theorem, and 
$$\sup_{|z-v|\leq \frac 1n\atop z,v\in \left[\frac{j_n}{n},\frac{J_n}{n}\right]}|F_\eps^{-1}(z)-F_\eps^{-1}(v)|\to 0$$
which can be guaranteed by assumption  \ref{as:eps} and the mean value theorem if $j_n/n$ and $J_n/n$ converge slowly enough. 
\hfill $\Box$\\
\\

{\bf Proof of Theorem \ref{theo1-mon}.}
The assertion follows from Theorem \ref{theo1} if we show that uniformly with respect to $t\in[0,1]$ and $y\in\R$, $S_n(t,y)=S_{n,I}(t,y)+o_P(1)$. To this end, observe that as in the proof of Theorem \ref{theo1} we can replace the estimators $\hat q_\tau$ and $\hat s$ by their linearized versions $\hat q_{\tau,L}$ and $\hat s_L$ in the definition of $S_n$ without changing the asymptotic properties. Denote the corresponding version of the process by $S_{n,L}$. Similarly, in the definition of $S_{n,I}$ the estimators $\hat q_{\tau,I}$ and $\hat s$ can be replaced by $\hat q_{\tau,L,I}=\Gamma_n(\hat q_{\tau,L})$ and $\hat s_L$, where $\hat q_{\tau,L,I}$ denotes the increasing rearrangement of the linearized estimator $\hat q_{\tau,L}$. More precisely, denoting this version of the process by $S_{n,L,I}$, we will show that
\beq \label{eq:h20}
\sup_{t\in[0,1] ,y\in\R} |S_{n,L,I}(t,y)- S_{n,I}(t,y)| = o_P(1).
\eeq
To see this, let $c=\inf_{x\in[0,1]}q_\tau'(x)$ and note that by our assumptions $c>0$ and by Lemma \ref{lem:unifratesquant} in the supplement we have for the set $\Omega_n := \{\sup_{x\in[h_n,1-h_n]}|\hat q_{\tau,L}'(x)-q_\tau'(x)|>\frac{c}{2}\}$ that $P(\Omega_n)\to 0$ for $n\to\infty$. Observe that by a straightforward modification of the proof of Theorem 3.1 (a) in Neumeyer (2007), we have on the set $\Omega_n$
\[
\sup_{x \in [h_n,1-h_n]}|\Gamma_n(\hat q_{\tau,L})(x) - \Gamma_n(\hat q_{\tau})(x)| \leq C \sup_{x \in [h_n,1-h_n]}|\hat q_{\tau,L}(x) - \hat q_{\tau}(x)|
\] 
for a universal constant $C$ which is independent of $n$. Thus Lemma \ref{lem:quantlin} in the supplement together with $P(\Omega_n) \to 1$ implies that 
\[
\sup_{x \in [h_n,1-h_n]}|\Gamma_n(\hat q_{\tau,L})(x) - \Gamma_n(\hat q_{\tau})(x)| = o_P(n^{-1/2}).
\]
Additionally, observe that the estimator $\hat q_{\tau,L}$ is strictly increasing provided that the event $\Omega_n$ holds, which implies that $P(\hat q_{\tau,L} \equiv \Gamma_n(\hat q_{\tau,L})) \geq P(\Omega_n) \to 1$. Now similar arguments as those used in the proof of Lemma \ref{lem:proclin} in the supplement show that, defining $F_{X,\eps_{L,I},n}$ in the same manner as $\hat F_{X,\eps_I,n}$ but with $\hat \eps_{i,L,I} := (Y_i - \Gamma_n(\hat q_{\tau,L})(x))/\hat s(X_i)$ instead of $\eps_{i,I}$, we have 
\[
\hat F_{X,\eps_I,n}(t,y) = \hat F_{X,\eps_{L,I},n}(t,y) + o_P(n^{-1/2})
\]
uniformly on $x \in [2h_n,1-2h_n], y \in \R$. Combining this with arguments which are similar to those in the proof of Theorem \ref{theo1}, this shows the validity of (\ref{eq:h20}). Next, note that on $\Omega_n$ the estimator $\hat q_{\tau,L}$ is strictly increasing. For every $\epsilon>0$ it follows that
\begin{eqnarray}\nb
&&P\Big(\sup_{t\in[2h_n,1-2h_n],y\in\R} |S_{n,L,I}(t,y)- S_{n,L}(t,y)|>\epsilon\Big) \\
\nb&=&P\Big(\sup_{t\in[2h_n,1-2h_n],y\in\R} |S_{n,L,I}(t,y)- S_{n,L}(t,y)|>\epsilon\Big)+o(1) \\
\nb
&\leq& 
P\Big(\sup_{t\in[2h_n,1-2h_n],y\in\R} |S_{n,L,I}(t,y)- S_{n,L}(t,y)|>\epsilon\; ,\;\sup_{x\in[h_n,1-h_n]}|\hat q_{\tau,L}'(x)-q_\tau'(x)|\leq\frac{c}{2} \Big) +o(1)
\\ \nb
&\stackrel{(*)}{\leq}& P\Big(\sup_{t\in[2h_n,1-2h_n],y\in\R} |S_{n,L,I}(t,y)- S_{n,L}(t,y)|>\epsilon\; ,\; \inf_{x\in[h_n,1-h_n]}\hat q_{\tau,L}'(x)>0\Big) +o(1)
\\
\nb
&=& o(1).
\end{eqnarray}
Here the last equality is due to the following argumentation. If $\inf_{x\in[h_n,1-h_n]}\hat q_{\tau,L}'(x)>0$, then $\hat q_{\tau,L}$ is strictly increasing, and for any increasing function the increasing rearrangement equals the original function function and we have $\hat q_{\tau,L,I}=\hat q_{\tau,L}$ (see Section \ref{sec-mon}).
But then, $S_{n,L}(t,y)=S_{n,L,I}(t,y)$ for all $t \in[2h_n,1-2h_n],y\in\R$ and the probability in ($*$) is zero. Finally, similar arguments as those in the proof of Theorem \ref{theo1} show that, uniformly with respect to $t\in [0,2h_n) \cup (1-2h_n,1],y\in\R$, we have $S_{n,L,I}(t,y) = S_{n,L}(t,y) + o_P(1)$. This completes the proof.
\hfill $\Box$

\section{Validity of bootstrap}\label{app-boot}
\def\theequation{B.\arabic{equation}}
\setcounter{equation}{0}

{\bf Preliminaries.}

Let $\tilde f_\eps$ denote the density corresponding to $\tilde F_\eps$.
Then under assumptions \ref{as:b1} analogous to Lemma 2 in Neumeyer (2009a) it can be shown that 
\begin{eqnarray}\label{le2n09}
\sup_{y\in\R}|\tilde f_\eps(y)-f_\eps(y)|=o_P((\frac{h_n}{\log n})^{1/2}),\quad
\sup_{y\in\R}|y\tilde f_\eps(y)-yf_\eps(y)|=o_P(1)\\
\sup_{y,z\in\R}\frac{|\tilde{f}_\eps(y)-f(y)-\tilde{f}_\eps(z)+f(z)|}{|y-z|^{\delta/2}}=o_P(1),\quad \sup_{y\in\R}|\tilde{F}_\eps(y)-F(y)|=o_P(1)\nonumber
\end{eqnarray}
(with $\delta$ from assumption \ref{as:b1}).
Further note that under assumption \ref{as:b2},
 Proposition 4 in Neumeyer (2009a) is valid (with $\upsilon$ from assumption \ref{as:b2})
and it follows that (for some constants $d$ and $L$) we have $\tilde F_\eps\in\mathcal{D}$  with probability converging to one. Here the function class is defined as
\begin{eqnarray}\label{cal-D}
\qquad\quad \mathcal{D}= \Big\{ F:\mathbb{R}\to [0,1]\;\Big|\; F \mbox{ increasing and continuously differentiable with derivative} &&\\
\nonumber \mbox{  } f  \mbox{ such that } 
\sup_{x\in\mathbb{R}}|f(x)|+\sup_{x,x'}\frac{|f(x)-f(x')|}{|x-x'|^{\delta/2}}\leq L, &&\\
 |1-F(x)|\leq \frac{d}{x^\upsilon}\forall x>0 \mbox{ and } |F(x)|\leq \frac{d}{|x|^\upsilon}\forall x<0
\Big\}.&&
\nonumber
\end{eqnarray}
From Lemma 4 in Neumeyer (2009a) and the conditions on $\delta$ and $\upsilon$ in assumption \ref{as:b2}  it follows that
\begin{equation}\label{cal-D-cov}
\log N(\epsilon,\mathcal{D},||\cdot||_\infty)=O(\epsilon^{-a}) \mbox{ for some }a<1.
\end{equation} 

\medskip

{\bf Proof of Theorem \ref{theo1-boot}.}

In the supplementary material in Lemma \ref{lem:proclin} it is shown that in the process $\hat F_{X,\eps,n}^*$ the residuals $\hat \eps_i^*$ can be replaced by  linearized versions $\hat\eps_{i,L}^*$ (see Appendix \ref{app-quant} in the supplement for the definitions). 
Using this, the preliminaries above as well as Lemma \ref{lem:unifratesquant} in the supplement (instead of Lemma 3 in Neumeyer (2009a)) we obtain analogously to the proofs of Lemma 1(i) and Theorem 2 in the reference that
\begin{eqnarray*}
&& \hat F_{X,\eps,n}^*(t,y)\\
 &=&
\frac{1}{n}\sum_{i=1}^n I\{\hat\eps_{i,L}^*\leq y\}I\{4h_n< X_i\leq t\}+o_P(\frac{1}{\sqrt{n}})\\
&=&\frac{1}{n}\sum_{i=1}^n I\{\eps_i^*\leq y\}I\{4h_n< X_i\leq t\}\\
&&{}+\int \Big(\tilde F_\eps\Big(y\frac{\hat s_L^*(x)}{\hat s_L(x)}+\frac{\hat q_{\tau,L}^*(x)-\hat q_{\tau,L}(x)}{\hat s_L(x)}\Big)-\tilde F_\eps(y)\Big)I\{4h_n<x\leq t\} f_X(x)\,dx
\\
&&{}+o_P(\frac{1}{\sqrt{n}})
\end{eqnarray*}
uniformly with respect to $t\in(4h_n,1-4h_n]$, $y\in\R$.
One can further apply a Taylor expansion for $\tilde F_\eps$. Lemma \ref{lem:B2boot} in the supplement gives expansions for the remaining integrals and we obtain
\begin{eqnarray*}
\hat F_{X,\eps,n}^*(t,y)
&=& \frac{1}{n}\sum_{i=1}^n I\{4h_n< X_i\leq t\}\Big( I\{\eps_i^*\leq y\}-\tilde\psi_n(y)\Big(I\{|\eps_i^*|\leq 1\}-\frac{1}{2}\Big)\\
&&{}\qquad\quad -\tilde\phi_n(y)\Big(I\{\eps_i^*\leq 0\}-\tau\Big)\Big)\\
&&{}+o_P(\frac{1}{\sqrt{n}})
\end{eqnarray*}
uniformly with respect to $t\in(4h_n,1-4h_n]$, $y\in\R$, 
where
$$\tilde \psi_n(y)=\frac{y\tilde f_\eps(y)}{f_{|\eps|}(1)}\, ,\quad \tilde\phi_n(y)=\frac{\tilde f_\eps(y)}{f_\eps(0)}\Big(1-y\frac{f_\eps(1)-f_\eps(-1)}{f_{|\eps|}(1)}\Big).$$
By the definition of the process $S_n^*$ one now directly has
\begin{eqnarray*}
&&S_n^*(t,y)\\
&=& \frac{1}{\sqrt{n}}\sum_{i=1}^n \Big( I\{\eps_i^*\leq y\}-\tilde\psi_n(y)\Big(I\{|\eps_i^*|\leq 1\}-\frac{1}{2}\Big) -\tilde\phi_n(y)\Big(I\{\eps_i^*\leq 0\}-\tau\Big)\Big)\\
&&\times \Big(I\{4h_n< X_i\leq t\}-I\{4h_n< X_i\leq 1-4h_n\}\frac{\hat F_{X,n}(t)-\hat F_{X,n}(4h_n)}{\hat F_{X,n}(1-4h_n)-\hat F_{X,n}(4h_n)}\Big)\\
&&{}+o_P(1)\\
&=& \frac{1}{\sqrt{n}}\sum_{i=1}^n g_n(\eps_i^*,y)\Big(I\{4h_n< X_i\leq t\}-I\{4h_n< X_i\leq 1-4h_n\}\frac{\hat F_{X,n}(t)-\hat F_{X,n}(4h_n)}{\hat F_{X,n}(1-4h_n)-\hat F_{X,n}(4h_n)}\Big)\\
&&{}+o_P(1)
\end{eqnarray*}
uniformly with respect to $t\in(4h_n,1-4h_n]$, $y\in\R$, with
\begin{eqnarray*}
&&g_n(\eps_i^*,y)\\
&=&I\{\eps_i^*\leq y\}-\tilde F_\eps(y)-\tilde\phi_n(y)\Big(I\{\eps_i^*\leq 0\}-\tilde F_\eps(0)\Big)-\tilde\psi_n(y)\Big(I\{|\eps_i^*|\leq 1\}-\tilde F_\eps(1)+\tilde F_\eps(-1)\Big).
\end{eqnarray*}
Note that $E[g_n(\eps_i^*,y)\mid\Y_n]=0$ and the dominating part of the process $S_n^*$ vanishes in the boundary points $t=4h_n$ and $t=1-4h_n$, for all $y\in\R$. Similarly to the corresponding arguments in the proof of Theorem \ref{theo1} (but with more technical effort) it can be shown that this process is equivalent in terms of conditional weak convergence in $\ell^\infty([0,1]\times\R)$ in probability to the process
\begin{eqnarray*}
S_{n,2}^*(t,y)
&=& \frac{1}{\sqrt{n}}\sum_{i=1}^n g_n(\eps_i^*,y)\Big(I\{X_i\leq t\}-\hat F_{X,n}(t)\Big),\quad t\in[0,1], y\in\R.
\end{eqnarray*}
Details are omitted for the sake of brevity. 

To finish the proof we have to show that (conditional on $\Y=((X_1,Y_1),(X_2,Y_2),\dots)$) the process $S_{n,2}^*$ converges weakly to $S$ in probability ($n\to\infty$). To this end we may show that for each subsequence $(n_k)_k$ there exists a further subsequence $(n_{k_\ell})_\ell$  such that (conditional on $\Y$) $S_{n_{k_\ell},2}^*$ converges weakly to $S$ almost surely ($\ell\to\infty$), cf.\ Sweeting (1989), p.\ 463. To this end we choose a subsequence $(n_{k_\ell})_\ell$ such that along this subsequence the convergences in (\ref{le2n09}) hold almost surely ($\ell\to\infty$). To simplify notation for the remainder of the proof we simply assume that the sequences in (\ref{le2n09})  converge almost surely ($n\to\infty$) and show that then $S_{n,2}^*$ converges weakly to $S$ almost surely ($n\to\infty$). 

It is easy to see that the conditional covariances $\Cov(S_{n,2}^*(s,y),S_{n,2}^*(t,z)\mid \Y)$ converge almost surely to $\Cov(S(s,y),S(t,z))$ as defined in Theorem \ref{theo1}. Thus it remains to show  conditional tightness and conditional fidi convergence of $S_{n,2}^*$. To obtain the latter we use CramÚr-Wold's device. Let $k\in\N$, $(y_1,t_1),\dots,(y_k,t_k)\in\R\times[0,1]$, $a_1,\dots,a_k\in\R$ and $Z_n=\sum_{j=1}^k a_jS_{n,2}^*(t_j,y_j)=n^{-1/2}\sum_{i=1}^n z_{n,i}$. Note that for some constant $c$, $|g_n(\eps_i^*,y)(I\{X_i\leq t\}-\hat F_{X,n}(t))|\leq 1+c(1+y)\tilde f_\eps(y)$, which converges almost surely to $1+c(1+y)f_\eps(y)$ due to (\ref{le2n09}) and thus is almost surely bounded. From this the validity of the conditional Lindeberg condition easily follows, i.\,e.\
\begin{eqnarray*}
L_n(\delta) &=& \frac{1}{n}\sum_{i=1}^n E[z_{n,i}^2
I\{|z_{n,i}|>n^{1/2}\delta\}\mid\mathcal{Y}]
\rightarrow 0 \mbox{ almost surely, for all } \delta>0.
\end{eqnarray*}

Finally, to prove conditional tightness we use the decomposition 
$S_{n,2}^*(t,y) = \sum_{k=0}^3 U_n^{(k)}(t,y)$, 
where
\begin{eqnarray*}
U_n^{(0)}(t,y)&=& \frac{1}{\sqrt{n}}\sum_{i=1}^n\Big(I\{\eps_i^*\leq y\}-\tilde F_\eps(y)\Big)I\{X_i\leq t\}\\
U_n^{(1)}(t,y)&=& -\tilde\phi_n(y)V_{n,1}(t)\\
&&\mbox{with }V_{n,1}(t)\;=\; \frac{1}{\sqrt{n}}\sum_{i=1}^n\Big(I\{\eps_i^*\leq 0\}-\tilde F_\eps(0)\Big)\Big(I\{X_i\leq t\}-\hat F_{X,n}(t)\Big)\\
U_n^{(2)}(t,y)&=&-\tilde \psi_n(y)V_{n,2}(t)\\
&&\mbox{with }V_{n,2}(t)\;=\; \frac{1}{\sqrt{n}}\sum_{i=1}^n\Big(I\{|\eps_i^*|\leq 1\}-\tilde F_\eps(1)+\tilde F_\eps(-1)\Big)\Big(I\{X_i\leq t\}-\hat F_{X,n}(t)\Big)\\
U_n^{(3)}(t,y)&=&-\hat F_{X,n}(t)W_n(y)\mbox{ with }W_n(y)\;=\; \frac{1}{\sqrt{n}}\sum_{i=1}^n\Big(I\{\eps_i^*\leq y\}-\tilde F_\eps(y)\Big).
\end{eqnarray*}
Note that conditional weak convergence of $V_{n,1}$ and $V_{n,2}$ to centered Gaussian processes, almost surely, can be shown analogously to the proof of bootstrap validity in Birke and Neumeyer (2013).
Further conditional weak convergence of $W_n$ is  completely analogous to Theorem 4 by Neumeyer (2009a).
From uniform almost sure convergence of $\phi_n$, $\psi_n$ and $\hat F_{X,n}$ to bounded functions, conditional tightness of $U_n^{(k)}$ follows for $k=1,2,3$. 

It remains to consider $U_n^{(0)}$. Applying Corollary 1 from Shorack and Wellner (1986), p.\ 622, (set $a=n^{-1}$, $b=\delta=\frac 12$, $\lambda=\sqrt{n}$) and the Borel-Cantelli lemma one obtains the existence of $c\in (0,\infty)$ such that with probability one
\begin{equation}\label{***neu}
 \quad |\hat F_{X,n}(t)-\hat F_{X,n}(s)|  \leq c|s-t|^{1/2} \quad \forall s,t \mbox{ with } n^{-1}\Delta_1^{-1}\leq |s-t| \leq \frac 12\Delta_2
\end{equation}
 for all but finitely many $n$, where 
 $\Delta_1=\inf_x f_X(x)>0$, $\Delta_2=\sup_x f_X(x)<\infty$.

We proceed by applying Theorem 2.11.9 by van der Vaart and Wellner (1996). Define $\mathcal{F} := [0,1]\times\R$ and for $f = (t,y)$ let
\[
Z_{ni}(f) := \frac{1}{\sqrt{n}}\Big(I\{\eps_i^*\leq y\}-\tilde F_\eps(y)\Big)I\{X_i\leq t\}.
\]

Let $\eta > 0$ and let $N_{[]}(\eta,\mathcal{F},L_2^n)$ denote the minimal number of sets $N_\eta$ in a partition of $\mathcal{F}$  in subsets $\mathcal{F}^n_{\eta j}$, $j=1,\dots,N_\eta$, such that for every $\mathcal{F}^n_{\eta j}$ 
\begin{equation}\label{*neu}
\sum_{i=1}^n E\Big[\sup_{f,g\in \mathcal{F}^n_{\eta j}} |Z_{ni}(f)-Z_{ni}(g)|^2\;\Big| \mathcal{Y}\Big]\leq \eta^2.
\end{equation}
Here the subsets are allowed to depend on $n$. 
Note also that we consider the conditional probability measure $P(\cdot\mid \mathcal{Y})$, so the sequence $(X_1,Y_1),(X_2,Y_2),\dots$ is given and the subsets are allowed to depend on it. 
We distinguish two cases.

\smallskip

{\bf 1.} Let $n\geq \Delta_1^{-1}\eta^{-4}$. 

Partition $[0,1]$ into $L=O(\eta^{-4})$ intervals $[t_{\ell-1},t_\ell]$, $\ell=1,\dots,L$ of length $\eta^{4}\leq t_\ell-t_{\ell-1}\leq 2\eta^{4}$ ($\forall \ell$).
Partition $\mathbb{R}$ into $K=O(\eta^{-2})$ intervals  $[y_{k-1},y_k]$, $k=1,\dots,K$, with $\tilde F_\eps(y_k)-\tilde F_\eps(y_{k-1})\leq \eta^2$ (using quantiles of the smooth distribution function $\tilde F_\eps$). 
The $N_\eta=LK$ intervals $[t_{\ell-1},t_\ell]\times[y_{k-1},y_k]$ define the subsets $\mathcal{F}^n_{\eta j}$, $j=1,\dots,N_\eta$. 

Now fix one subset and let  $f,g\in\mathcal{F}^n_{\eta j}=[t_{\ell-1},t_\ell]\times[y_{k-1},y_k]$. Then for monotonicity reasons $Z_{ni}(f)$ as well as $Z_{ni}(g)$ are elements of the bracket $[Z_{ni}^{k,\ell,l},Z_{ni}^{k,\ell,u}]$,
where
\begin{eqnarray*} 
Z_{ni}^{k,\ell,l} &=& \frac{1}{\sqrt{n}}\Big(I\{\eps_i^*\leq y_{k-1}\}I\{X_i\leq t_{\ell -1}\} -\tilde F_\eps(y_k)I\{X_i\leq t_{\ell}\}\Big)\\
Z_{ni}^{k,\ell,u} &=& \frac{1}{\sqrt{n}}\Big(I\{\eps_i^*\leq y_{k}\}I\{X_i\leq t_{\ell }\} -\tilde F_\eps(y_{k-1})I\{X_i\leq t_{\ell -1}\}\Big).
\end{eqnarray*}
Thus the left hand side of (\ref{*neu}) can be bounded by 
\begin{eqnarray} 
&&\sum_{i=1}^n E\Big[(Z_{ni}^{k,\ell,u}-Z_{ni}^{k,\ell,l})^2\;\Big| \mathcal{Y}\Big] \nonumber
\\
&\leq & \frac{2}{n}\sum_{i=1}^n(I\{X_i\leq t_{\ell }\}-I\{X_i\leq t_{\ell -1}\})^2 \nonumber
\\
&&{}+ \frac{2}{n}\sum_{i=1}^n E\Big[\Big(I\{\eps_i^*\leq y_{k}\}-\tilde F_\eps(y_{k-1})-I\{\eps_i^*\leq y_{k-1}\}+\tilde F_\eps(y_{k})\Big)^2\;\Big| \mathcal{Y}\Big] \nonumber
\\
&\leq & \frac{2}{n}\sum_{i=1}^n(I\{X_i\leq t_{\ell }\}-I\{X_i\leq t_{\ell -1}\})^2 \nonumber
\\
&&{}+ \frac{4}{n}\sum_{i=1}^n E\Big[I\{\eps_i^*\leq y_{k}\}-I\{\eps_i^*\leq y_{k-1}\} + \tilde F_\eps(y_{k}) - \tilde F_\eps(y_{k-1})\;\Big| \mathcal{Y}\Big] \nonumber
\\
&\leq&
2 (\hat F_{X,n}(t_{\ell})-\hat F_{X,n}(t_{\ell-1}))+8(\tilde F_\eps(y_k)-\tilde F_\eps(y_{k-1})) \nonumber
\\
&\leq & 2 (\hat F_{X,n}(t_{\ell})-\hat F_{X,n}(t_{\ell-1}))+8\eta^2\;\leq\; C\eta^2, \label{eq:boundz3}
\end{eqnarray}
where we have used (\ref{***neu}) and $t_\ell-t_{\ell-1}\geq\eta^{4}\geq n^{-1}\Delta_1^{-1}$, and the constant $C$ does not depend on $n$ and $\eta$.

\smallskip

{\bf 2.} Let $n < \Delta_1^{-1}\eta^{-4}$. 

As before we partition  $\mathbb{R}$ into $K=O(\eta^{-4})$ intervals  $[y_{k-1},y_k]$, $k=1,\dots,K$, with $\tilde F_\eps(y_k)-\tilde F_\eps(y_{k-1})\leq \eta^2$.
We partition $[0,1]$ into $n+2=O(\eta^{-4})$ intervals $I_\ell=[t_{\ell-1},t_\ell)$, $\ell=1,\dots,n+1$, and $I_{n+2}=\{1\}$, where $t_0=0$, $t_\ell=X_{(\ell)}$ for $\ell=1,\dots,n$ and $t_{n+1}=1$. 
Here $X_{(1)},\dots,X_{(n)}$ denote the order statistics of $X_1,\dots,X_n$. 
Now we proceed as in case 1 but replacing $Z_{ni}^{k,\ell,u},Z_{ni}^{k,\ell,l}$ with 
\begin{eqnarray*} 
\tilde{Z}_{ni}^{k,\ell,l} &=& \frac{1}{\sqrt{n}}\Big(I\{\eps_i^*\leq y_{k-1}\}I\{X_i\leq t_{\ell -1}\} -\tilde F_\eps(y_k)I\{X_i < t_{\ell}\}\Big)\\
\tilde{Z}_{ni}^{k,\ell,u} &=& \frac{1}{\sqrt{n}}\Big(I\{\eps_i^*\leq y_{k}\}I\{X_i < t_{\ell }\} -\tilde F_\eps(y_{k-1})I\{X_i\leq t_{\ell -1}\}\Big).
\end{eqnarray*}
By definition, $\tilde{Z}_{ni}^{k,\ell,l} \leq Z_{ni}(f) \leq \tilde{Z}_{ni}^{k,\ell,u}$ for $f = (t,y) \in [t_{\ell-1},t_\ell)\times[y_{k-1},y_k]$. Noting that $\hat F_{X,n}(t_\ell-)-\hat F_{X,n}(t_{\ell-1})=0$ for all $\ell=1,\dots,n+1$, we obtain by similar arguments as used to derive (\ref{eq:boundz3})
\begin{eqnarray*} 
\sum_{i=1}^n E\Big[(\tilde{Z}_{ni}^{k,\ell,u}-\tilde{Z}_{ni}^{k,\ell,l})^2\;\Big| \mathcal{Y}\Big]
&\leq & 2 (\hat F_{X,n}(t_{\ell}-)-\hat F_{X,n}(t_{\ell-1}))+ 8\eta^2 = 8\eta^2.
\end{eqnarray*}

\smallskip

The partitionings in both cases depend on $n$, but the bracketing number $N_{[]}(\eta,\mathcal{F},L_2^n)$ can be bounded by $O(\eta^{-8})$, independent of $n$, such that the condition
$$\int_0^{\delta_n}\sqrt{\log N_{[]}(\eta,\mathcal{F},L_2^n)}\,d\eta\longrightarrow 0\mbox{ for every }\delta_n\searrow 0$$
is fulfilled (this corresponds to the third condition in Theorem 2.11.9 by van der Vaart and Wellner (1996)). 
Further, because $|Z_{ni}(f)|\leq n^{-1/2}$ $\forall f$ we have
$$\sum_{i=1}^n E\Big[\sup_{f\in\mathcal{F}}|Z_{ni}(f)|I\{\sup_{f\in\mathcal{F}}|Z_{ni}(f)|>\eta\}\;\Big|\mathcal{Y}\Big]\longrightarrow 0\mbox{ for every }\eta>0$$
(this corresponds to the first condition in Theorem 2.11.9 by van der Vaart and Wellner (1996)).
Moreover, $(\mathcal{F},\rho)$ is a totally bounded semimetric space with $\rho((s,y),(t,z))=|t-s|+|F_\eps(z)-F_\eps(y)|$. Now for $\delta_n\searrow 0$ we obtain similarly to the calculation in case 1 above (for some constant $c$),
\begin{eqnarray*} 
&&\sup_{\rho(f,g)<\delta_n}\sum_{i=1}^n E\Big[(Z_{ni}(f)-Z_{ni}(g))^2 \;\Big|\mathcal{Y}\Big]\\
&\leq & c \Big(\sup_{|t-s|\leq\delta_n}|\hat F_{X,n}(t)-\hat F_{X,n}(s)|+\sup_{z,y:\atop |F_\eps(z)-F_\eps(y)|\leq\delta_n}|\tilde F_\eps(z)-\tilde F_\eps(y)|\Big)\\
&=&o(1) \mbox{ almost surely}
\end{eqnarray*}
by uniform convergence of $\hat F_{X,n}$ to $F_X$ and $\tilde F_\eps$ to $F_\eps$
(this corresponds to the second condition in Theorem 2.11.9 by van der Vaart and Wellner (1996)) and uniform continuity of $F_X$.
From Theorem 2.11.9 one obtains
$$
\lim_{\delta\searrow 0}\lim_{n\to\infty}P\Bigg(  \sup_{\rho((s,y),(t,z))<\delta}|\tilde U_n^{(0)}(s,y)-\tilde U_n^{(0)}(t,z)|>\eta
\;\Bigg|\; \mathcal{Y}\Bigg) =0 \mbox{ for all $\eta>0$} 
$$
for almost all sequences $\Y$.
This completes the proof. \hfill $\Box$\\
\\

{\bf Proof of Theorem \ref{theo1-mon-boot}.}

 Theorem \ref{theo1-mon-boot} follows from Theorem \ref{theo1-boot} in the same manner as Theorem \ref{theo1-mon} follows from Theorem \ref{theo1}.
\hfill $\Box$


\vspace{1cm}

{\LARGE \textbf{References}}

\begin{description}
\item[J.\ Abrevaya] (2005). Isotonic quantile regression: asymptotics and bootstrap. Sankhy\=a 67, 187--199.
\item[M.\ Akritas and I.\ Van Keilegom] (2001). \emph{Nonparametric estimation of the residual distribution.} Scand. J. Statist. 28, 549--567.
\item[D.\ Anevski and A.-L. FougÞres] (2007). \emph{Limit properties of the monotone rearrangement for density and regression function estimation.} arXiv:0710.4617v1
\item[Y.\ Baraud, S.\ Huet and B.\ Laurent] (2003). \emph{Adaptive tests of qualitative hypotheses}.  ESAIM, Probab. Statist. 7, 147--159.
\item[M.\ Birke and H.\ Dette] (2007). \emph{Testing strict monotonicity in nonparametric regression.} Math. Meth. Statist. 16, 110--123.
\item[M.\ Birke and H.\ Dette] (2008). \emph{A note on estimating a smooth monotone regression by combining kernel and density estimates.} J. Nonparam. Statist. 20, 679--690.
\item[M.\ Birke and N.\ Neumeyer] (2013). \emph{Testing monotonicity of regression functions - an empirical process approach.} Scand. J. Statist. 40, 438--454.
\item[J.R.\ Blum, J. Kiefer and M. Rosenblatt] (1961). \emph{Distribution free tests of independence based on the sample distribution functions.} Ann. Math. Stat. 32, 485--498.
\item[A.\,W.\ Bowman, M.\,C.\ Jones and I.\ Gijbels] (1998). \emph{Testing monotonicity of regression.}  J. Comput. Graph. Stat. 7, 489--500.
\item[R.J.\ Casady and J.D.\ Cryer] (1976). \emph{Monotone percentile regression.} Ann. Stat. 4, 532--541.
\item[S. Chen, G.B. Dahl and S. Khan] (2005). \emph{Nonparametric identification and estimation of a censored location-scale regression model.}
J. Amer. Statist. Assoc. 100, 212--221.
\item[J.D.\ Cryer, T.\ Robertson, F.T.\ Wright and R.J.\ Casady] (1972). \emph{Monotone median regression.} Ann. Math. Stat. 43, 1459--1469.
\item[H.\ Dette, N.\ Neumeyer and K.\,F.\ Pilz] (2006). \emph{A simple nonparametric estimator of a strictly monotone regression function.}  Bernoulli  12, 469--490.
\item[H. Dette and S. Volgushev] (2008). \emph{Non--crossing nonparametric estimates of
quantile curves.} J. Roy. Stat. Soc. B 70, 609--627.
\item[J.\ Dom\'inguez-Menchero, G.\  Gonz\'alez-Rodr\'iguez and M.\,J.\ L\'opez-Palomo] (2005). \emph{An L2 point of view in testing monotone regression}  J. Nonparam. Statist. 17, 135--153.
\item[L.\ Duembgen] (2002). \emph{Application of local rank tests to nonparametric regression.}  J. Nonparametr. Stat.  14, 511--537.
\item[C.\ Durot] (2003). \emph{A Kolmogorov-type test for monotonicity of regression.} Statist. Probab. Lett. 63, 425--433.
\item[S.\ Efromovich] (1999). Nonparametric curve estimation. Methods, theory, and applications. Springer, New York.
\item[J.H.J.\ Einmahl and I.\ Van Keilegom] (2008a). \emph{Specification tests in nonparametric regression.} J. Econometr. 143, 88--102.
\item[J.H.J.\ Einmahl  and I.\ Van Keilegom] (2008b). \emph{Tests for independence in nonparametric regression.} Statist. Sinica 18, 601--616. 
\item[X.\ Feng, X.\ He and J.\ Hu] (2011). \emph{Wild bootstrap for quantile regression.} Biometrika 98, 995--999.
\item[B. Fitzenberger, K. Kohn and A. Lembcke] (2008). \emph{Union Density and Varieties of Coverage: The Anatomy of Union  Wage Effects in Germany.} IZA Working Paper No. 3356. Available at SSRN: http://ssrn.com/abstract=1135932 
\item[T.\ Gasser, H.\ M\"uller and V. Mammitzsch] (1985). \emph{Kernels for nonparametric curve estimation.} J. Roy. Stat. Soc. B, 238--252.
\item[I.\ Gijbels, P.\ Hall, M.\,C.\ Jones and I.\ Koch] (2000). \emph{Tests for monotonicity of a regression mean with guaranteed level.} Biometrika 87, 663--673.
\item[S.\ Ghosal, A.\ Sen and A.\,W.\ van der Vaart] (2000). \emph{Testing monotonicity of regression.} Ann. Statist. 28, 1054--1082.
\item[P.\ Hall and N.\,E.\ Heckman] (2000).  \emph{Testing for monotonicity of a regression mean by calibrating for linear functions.} Ann. Statist. 28, 20--39.
\item[W.\ Hõrdle and E.\ Mammen] (1993). \emph{Comparing nonparametric versus parametric\linebreak regression fits.} Ann.
Statist. 21, 1926--1947.
\item[X. He] (1997). \emph{Quantile Curves without Crossing.} Am. Stat. 51, 186--192.
\item[W. Hoeffding] (1948). \emph{A nonparametric test of independence.} Ann. Math. Statist. 19, 546--557. 
\item[S.\ Kiwitt, E.-R.\ Nagel and N.\ Neumeyer] (2008). \emph{Empirical Likelihood Estimators for the Error Distribution in Nonparametric Regression Models.} Math. Meth. Statist. 17, 241--260.
\item[R. Koenker] (2005). \emph{Quantile Regresssion.} Cambridge University Press, Cambridge.
\item[R. Koenker and G. Bassett] (1978). \emph{Regression quantiles.} Econometrica 46, 33--50.
\item[U.\,U.\ M³ller, A.\ Schick and W.\ Wefelmeyer] (2004). \emph{Estimating linear functionals of the error distribution in nonparametric regression.}  J. Statist. Plann. Inf.  119, 75--93.
\item[N. Neumeyer] (2007). \emph{A note on uniform consistency of monotone function estimators.} Statist. Probab. Lett. 77, 693--703. 
\item[N.\ Neumeyer] (2009a). \emph{Smooth residual bootstrap for empirical processes of nonparametric regression residuals.}  Scand. J. Statist. 36, 204--228.
\item[N.\ Neumeyer] (2009b). \emph{Testing  independence in  nonparametric regression.} J. Multiv. Anal. 100, 1551--1566.
\item[J. Rice](1984). \emph{Bandwidth choice for nonparametric regression.} Ann. Statist. 12, 1215--1230.
\item[T.\ Robertson and F.T.\ Wright] (1973). \emph{Multiple isotonic median regression.} Ann. Statist. 1, 422--432.
\item[J.\ Shim, C.\ Hwang and H.K.\ Seok] (2009). \emph{Non-crossing quantile regression via doubly penalized kernel machine.} Comput. Statist. 24, 83--94.
\item[G.\,R.\ Shorack and J. A.\ Wellner] (1986). \emph{Empirical Processes with Applications to Statistics.} Wiley, New York.
\item[W.\ Stute] (1997). \emph{Nonparametric model checks for regression.} Ann. Statist. 25, 613--641.
\item[Y.\ Sun] (2006). \emph{A consistent nonparametric equality test of conditional quantile functions.} Econometric Th. 22, 614--632.
\item[T.\,J.\ Sweeting] (1989). \emph{On Conditional Weak Convergence.} J. Theoret. Probab. 2, 461--474.
\item[I. Takeuchi, Q.V. Le, T.D. Sears, and A.J. Smola](2006). \emph{Nonparametric
quantile regression.} Journal of Machine Learning Research, 7:1231--1264.
\item[A.\,W.\ van der Vaart and J. A. Wellner] (1996). \emph{Weak convergence and empirical processes.} Springer, New York.
\item[I.\ Van Keilegom] (1998). \emph{Nonparametric estimation of the conditional distribution in regression with censored data.} PhD thesis, University of Hasselt, Belgium. available at http.//www.ibiostat.be/publications/.
\item[I.\ Van Keilegom, W.\ Gonzßlez--Manteiga and C.\  Sßnchez Sellero] (2008). \emph{Goodness-of-fit tests in parametric regression based on the estimation of the error distribution.} TEST 17, 401--415.
\item[S.\ Volgushev, M. Birke, H. Dette and N. Neumeyer] (2013). \emph{Significance testing in quantile regression.} Electron. J. Stat. 7, 105--145. 
\item[J.C. Wang and M.C. Meyer] (2011). \emph{Testing the monotonicity or convexity of a function using regression splines.} Can. J. Stat. 39, 89--107.
\item[K. Yu, and M.C. Jones](1997). \emph{A comparison of local constant and local linear
regression quantile estimators.} Comput. Stat. Data Anal. 25, 159--166.
\item[K. Yu and M.C. Jones] (1998). \emph{Local linear quantile regression.} J. Amer. Statist. Assoc. 93, 228--237.
\item[K. Yu, Z. Lu, and J. Stander] (2003). \emph{Quantile regression: applications and current research areas.} J. Roy. Stat. Soc. D (The Statistician) 52, 331--350.
\end{description}

\bigskip

\hrulefill

\bigskip

\noindent \sc Melanie Birke, \rm  Universit\"{a}t Bayreuth, Fakult\"{a}t f\"{u}r Mathematik, Physik und Informatik , 95440 Bayreuth, Germany, e-mail:  Melanie.Birke@uni-bayreuth.de

\smallskip

\noindent \sc Natalie Neumeyer, \rm Universitõt Hamburg, Fachbereich  Mathematik, Bundesstra▀e 55, 20146 Hamburg, Germany, e-mail: neumeyer@math.uni-hamburg.de

\smallskip

\noindent\sc Stanislav Volgushev, \rm University of Toronto, Department of Statistical Sciences 100 St.
George Street, Toronto, Ontario, Canada, e-mail: volgushe@utstat.toronto.edu

\newpage

\begin{figure}[t]
\begin{center}
\epsfig{file=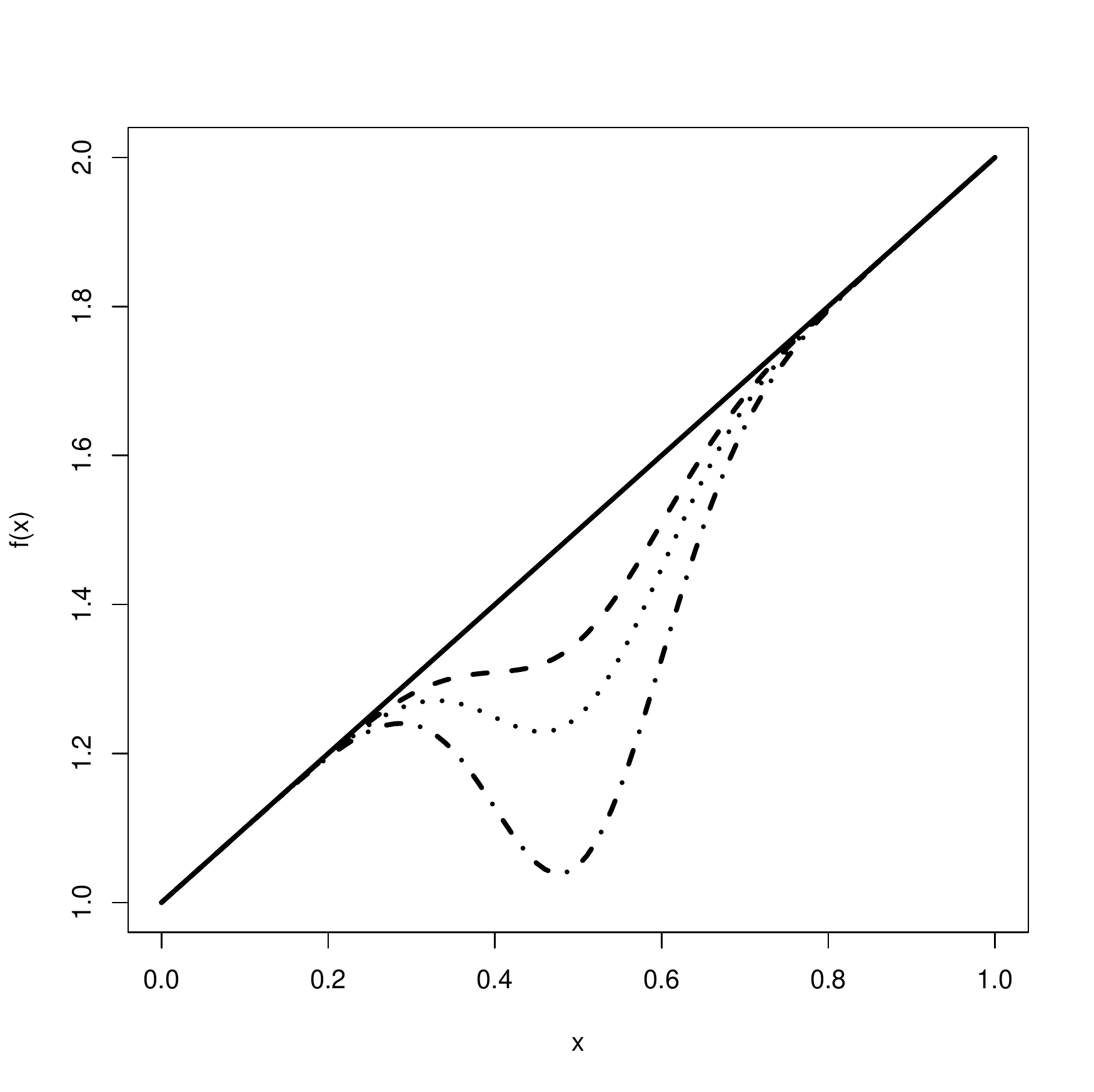,width=8.5cm}
\end{center}
\caption{\label{fig3} \sl The function $x \mapsto 1 + x - \beta e^{-50(x-0.5)^2}$ for values $\beta = 0$ (solid line), $\beta = 0.15$ (dashed line), $\beta = 0.25$ (dotted line), $\beta = 0.45$ (dash-dotted line), respectively. This is the median function in model 4, see Section \ref{sec:monsim}.}
\end{figure}

\begin{table}[h!]
\begin{center}
\begin{tabular}{|c|c|ccc|}
\hline
& &  model 1 &  model 2a & model 3  
\\
& &  $a=0$ &  $c = 2$ & $b=0$
\\
\hline 
KS & $n = 100$ & 0.034 & 0.039 & 0.045
\\
CvM & $n = 100$ & 0.029 & 0.044 & 0.053
\\
\hline 
KS & $n = 200$ & 0.034 & 0.039 & 0.050
\\
CvM & $n = 200$ & 0.046 & 0.049 & 0.062
\\
\hline
\end{tabular}
\end{center}
\caption  {\label{sizeloc}
\it Rejection probabilities for testing the validity of a location model under various $H_0$ scenarios, the nominal level is $\alpha = 5\%$.}
\end{table}

\begin{table}[h!]
\begin{center}
\begin{tabular}{|c|c|ccccc|}
\hline 
& a& 0 &1 & 2.5 & 5 & 10  
\\
\hline 
\hline 
KS & $n = 100$ & 0.032 & 0.078 & 0.16 & 0.23 & 0.444
\\
CvM & $n = 100$ & 0.038 & 0.128 & 0.364 & 0.568 & 0.746 
\\ 
\hline 
N & $n = 100$ & 0.054 & 0.190& 0.506& 0.734& 0.884
\\
\hline 
EVK & $n = 100$ & 0.072 &0.132& 0.316& 0.524& 0.668  
\\
\hline 
\hline 
KS & $n = 200$ & 0.034 & 0.144& 0.292& 0.586& 0.784  
\\
CvM & $n = 200$ & 0.046 & 0.296& 0.632&  0.9& 0.976 
\\ 
\hline 
N & $n = 200$ & 0.044 & 0.390& 0.860& 0.976& 0.972
\\
\hline 
EVK & $n = 200$ & 0.066 &0.376& 0.788& 0.960& 1.00
\\
\hline
\end{tabular}
\end{center}
\caption  {\label{sizeloc1}
\it  Rejection probabilities for testing the validity of a location model under the alternative in model 1 for different values of the parameter $a$, nominal level is $\alpha = 5\%$}
\end{table}

\begin{table}[h!]
\begin{center}
\begin{tabular}{|c|c|ccccc|}
\hline 
& c& .2 & .4 & .6 & .8 & 1 
\\
\hline 
\hline
KS & $n = 100$ & 0.044 & 0.074 & 0.120 & 0.194 & 0.390
\\
CvM & $n = 100$ & 0.082 & 0.124 & 0.218 & 0.414 & 0.768
\\ 
\hline 
N & $n = 100$ & 0.096& 0.120& 0.224& 0.420& 0.676
\\
\hline 
EVK & $n = 100$ &0.116& 0.160& 0.224& 0.360 & 0.612
\\
\hline 
\hline 
KS & $n = 200$ & 0.08 & 0.136 & 0.222 & 0.4 & 0.762
\\
CvM & $n = 200$ & 0.118 & 0.29 & 0.49 & 0.792 & 0.996 
\\ 
\hline 
N & $n = 200$ & 0.156 & 0.216 & 0.412 & 0.688 & 0.904
\\
\hline 
EVK & $n = 200$ & 0.124 & 0.216 & 0.344 & 0.584 & 0.944
\\
\hline
\end{tabular}
\end{center}
\caption  {\label{sizeloc2}
\it  Rejection probabilities for testing the validity of a location model under the alternative in model 2b for different values of the parameter $c$, nominal level is $\alpha = 5\%$}
\end{table}

\newpage

\begin{table}[h!]
\begin{center}
\begin{tabular}{|c|c|ccccc|}
\hline
& b& 0 & 1 & 2 & 3 & 5 
\\
\hline
\hline  
KS & $n=100$ & 0.045 & 0.094 & 0.154 & 0.306 & 0.712
\\
CvM & $n=100$ &  0.053 & 0.128 & 0.240 & 0.576 & 0.968
\\ 
\hline 
N & $n=100$ & 0.024& 0.172& 0.284& 0.452& 0.662
\\
\hline 
\hline 
KS & $n=200$ &  0.050  & 0.134 & 0.31 & 0.518 & 0.906 
\\
CvM & $n=200$ &  0.062 & 0.254 & 0.538 & 0.92 &  1
\\ 
\hline 
N & $n=200$ & 0.034& 0.620& 0.926& 0.998& 1.000
\\
\hline
\end{tabular}
\end{center}
\caption{\label{sizeloc3}
\it  Rejection probabilities for testing the validity of a location model under the alternative in model 3 for different values of the parameter $b$, nominal level is $\alpha = 5\%$}
\end{table}

\begin{table}[h!]
\begin{center}
\begin{tabular}{|c|c|ccc|}
\hline
& &  model $1_h$ &  model $2a_h$ & model $3_h$  
\\
& &    &  $c = 2$ & $b=0$
\\
\hline
KS & $n = 50$ & 0.025 & 0.026 & 0.023
\\
CvM & $n = 50$ & 0.022 & 0.026 & 0.034
\\
\hline 
KS & $n = 100$ & 0.031 &  0.037 & 0.037
\\
CvM & $n = 100$ & 0.029 & 0.031 & 0.041
\\
\hline 
KS & $n = 200$ & 0.024 & 0.044 & 0.057
\\
CvM & $n = 200$ & 0.028 & 0.044 & 0.062
\\
\hline
\end{tabular}
\end{center}
\caption  {\label{sizelocscal}
\it  Rejection probabilities for the testing the validity of a location-scale model under various $H_0$ scenarios, the nominal level is $\alpha = 5\%$. }
\end{table}

\begin{table}[h!]
\begin{center}
\begin{tabular}{|c|c|c|}
\hline
& c & 1 
\\
\hline 
\hline
KS & $n=50$ &  0.032 
\\
CvM & $n=50$ & 0.034 
\\ 
\hline 
EVK & $n=50$ &  0.262
\\
\hline
\hline 
KS & $n=100$ & 0.046
\\
CvM & $n=100$ &  0.04
\\
\hline 
EVK & $n=100$ & 0.478
\\
\hline
\end{tabular}
\end{center}
\caption{\label{sizelocscal2}
\it  Rejection probabilities for testing the validity of a location-scale model under the alternative in model $2b_h$, the nominal level is $\alpha = 5\%$}
\end{table}

\begin{table}[h!]
\begin{center}
\begin{tabular}{|c|c|ccccc|}
\hline
& b& 0 & 1 & 2 & 3 & 5 
\\
\hline
\hline 
KS & $n=100$ &  0.037 & 0.212 & 0.344 & 0.546 & 0.878
\\
CvM & $n=100$ & 0.041 & 0.368 & 0.658 & 0.922 & 0.992
\\ 
\hline 
N & $n=100$ & 0.036& 0.278 & 0.388 & 0.190 & 0.156   
\\
\hline 
\hline 
KS & $n=200$ & 0.057  & 0.452 & 0.646 & 0.8 & 0.972 
\\
CvM & $n=200$ & 0.062  & 0.802 & 0.966 &   1 & 1
\\ 
\hline 
N & $n=200$ & 0.035 & 0.630 & 0.774 & 0.402 & 0.268 
\\
\hline
\end{tabular}
\end{center}
\caption{\label{sizelocscal3}
\it  Rejection probabilities for testing the validity of a location-scale model under the alternative in model $3_h$ for different values of the parameter $b$, the nominal level is $\alpha = 5\%$.}
\end{table}

\begin{table}[h!]
\begin{center}
\begin{tabular}{|c|c|c|c|c|c|c|}
\hline
&\multicolumn{3}{|c|}{$\tau = 0.25$}&
\multicolumn{3}{|c|}{$\tau = 0.5$}
\\
\hline
& \multicolumn{1}{|c|}{$ n = 50 $}&\multicolumn{1}{|c|}{$ n = 100 $}&\multicolumn{1}{|c|}{$ n = 200 $}&
\multicolumn{1}{|c|}{$ n = 50 $}&\multicolumn{1}{|c|}{$ n = 100 $}&\multicolumn{1}{|c|}{$ n = 200 $}
\\
\hline
$\beta=0$
&0.020  & 0.020 & 0.026 
&0.025  & 0.023 & 0.026
\\
\hline
$\beta=0.15$
&0.024  & 0.027 & 0.050 
&0.027  & 0.047 & 0.060
\\
\hline
$\beta=0.25$
&0.028  & 0.057 & 0.126 
&0.037  & 0.053 & 0.154
\\
\hline
$\beta=0.45$
&0.140  & 0.202 & 0.410 
&0.084  & 0.154 & 0.344
\\
\hline
\end{tabular}
\end{center}
\caption  {\label{sizepower4}
\it Rejection probabilities for the test for monotonicity of quantile curves in model 4. The nominal level is $\alpha = 5\%$.}
\end{table}

\begin{table}[h!]
\begin{center}
\begin{tabular}{|c|c|c|c|}
\hline
& \multicolumn{1}{|c|}{$ n = 50 $}&\multicolumn{1}{|c|}{$ n = 100 $}&\multicolumn{1}{|c|}{$ n = 200 $}\\
\hline
$\tau=0.25$ & 0.23 & 0.262 & 0.376
\\
\hline
$\tau=0.5$ & 0.073 & 0.061 & 0.043
\\
\hline
$\tau=0.75$ & 0.181 & 0.180 & 0.296
\\
\hline
\end{tabular}
\end{center}
\caption  {\label{sizepower5}
\it Rejection probabilities for the test for monotonicity of quantile curves in model 5. Different rows correspond to the $0.25,0.5$ and $0.75$ quantile curves, respectively. The nominal level is $\alpha = 5\%$.}
\end{table}


\clearpage

\section{Supplement to ``The independence process in conditional quantile location-scale models and an application to testing for monotonicity'' by Melanie Birke, Natalie Neumeyer and Stanislav Volgushev --- Technical results} 
\pagenumbering{roman}
\setcounter{page}{1}
\def\theequation{C.\arabic{equation}}
\setcounter{equation}{0}

We begin by recalling some notation from the main body of the paper that will be used throughout the proofs. 

One fact that we will use throughout is that the bootstrap residuals $\eps_i^*$ can be represented as $\eps_i^* = \tilde F_\eps^{-1}(U_i)$ where $U_1,...,U_n$ denote a sample of i.i.d. $\mathcal{U}[0,1]$ random variables that are independent of the original sample and 
\[
\tilde{F}_\eps(y) = \frac{\frac{1}{n}\sum_{i=1}^n \Phi\Big(\frac{y-\hat{\eps}_i}{\alpha_n}\Big)I\{2h_n< X_i\leq 1-2h_n \}}{\hat F_{X,n}(1-2h_n)-\hat F_{X,n}(2h_n)}
\]
denotes the distribution function of $\eps_1^*$ conditional on the sample, see (\ref{eq:feps}).
Additionally, we will use the abbreviation 
\[
r_n:= \Big(\frac{\log n}{n h_n}\Big)^{1/2}.
\]

Next, we introduce some additional notation that will be used throughout. First, introduce the functional
\[
Q_{G,\kappa,\tau,b_n}(F) := G^{-1}\Big(\frac{1}{b_n}\int_0^1 \int_{-\infty}^\tau \kappa\Big(\frac{F(G^{-1}(u)) - v}{b_n} \Big) dvdu \Big)
\]
which is defined for arbitrary functions $F$ that are uniformly bounded. Some properties of this functional are collected in Lemma \ref{lem:genlin}. Additionally, define the quantities
\bean
\hat F_Y^*(y|x) &:=& \sum_{i=1}^n W_i(x) \Omega \Big(\frac{y-Y_i^*}{d_n}\Big), \quad \hat q_{\tau}^*(x) :=  Q_{G,\kappa,\tau,b_n}(\hat F_{Y}^*(\cdot|x)),
\\
\hat F_{|e|}^*(y|x) &:=& \sum_{i=1}^n W_i(x) \Omega \Big(\frac{y-|Y_i^*-\hat q_{\tau}^*(X_i)|}{d_n}\Big), \quad \hat s^*(x) :=  Q_{G,\kappa,1/2,b_n}(\hat F_{|e|}^*(\cdot|x)).
\eean
where the weights $W_i$ are the same as in equation (\ref{def-LL}). Observe that the estimators $\hat q_\tau, \hat s$ which we introduced in the main body of the paper admit the representations
\[
\hat q_{\tau}(x) =  Q_{G,\kappa,\tau,b_n}(\hat F_{Y}(\cdot|x)), \quad \hat s(x) :=  Q_{G,\kappa,1/2,b_n}(\hat F_{|e|}(\cdot|x)).
\]
 
In appendix \ref{app-quant}, we will introduce linearized versions of the estimators $\hat q_\tau, \hat q_{\tau}^*, \hat s, \hat s^*$, those will be denoted by $\hat q_{\tau,L}, \hat q_{\tau,L}^*, \hat s_L, \hat s^*_L$. Key results there are Lemma \ref{lem:unifratesquant} and Lemma \ref{lem:quantlin} which state that the linearized versions are uniformly close to the original estimators and that the linearized versions have certain smoothness properties, respectively. The rest of the Appendix is organized as follows. Section \ref{app-quant} contains results about the estimators $\hat q_\tau, \hat q_{\tau}^*, \hat s, \hat s^*$ and their linearizations. The proofs of those results require additional technical Lemmas, that we collect and prove in Section \ref{sec:tec}. Finally, some key results which are used in the main body of the paper and whose proofs rely on findings in Sections \ref{app-quant} and \ref{sec:tec} can be found in Section \ref{app-tech}. 


\subsection{Properties of $\hat q_\tau$ and $\hat s$}\label{app-quant}
We start this section by introducing some notation and giving an overview of the derived results. Our first key result is an asymptotic representation of the form 
\bean
\hat F_Y(y|x) &=& \hat F_{Y,L,S}(y|x) + o_P(1/\sqrt n), \qquad \hat F_{|e|}(y|x) = F_{|e|,L,S}(y|x) + o_P(1/\sqrt n),
\\ 
\hat F_Y^*(y|x) &=& \hat F_{Y,L,S}^*(y|x) + o_P(1/\sqrt n), \qquad \hat F_{|e|}^*(y|x) = F_{|e|,L,S}^*(y|x) + o_P(1/\sqrt n),
\eean
holding uniformly over $x,y$ where the expressions on the right-hand side of the above equations are defined as 
\bean
\hat F_{Y,L,S}(y|x) &:=& F_Y(y|x) + u_1^t\mathcal{M}(K)^{-1} \Big( T_{n,0,L,S}(x,y),\dots, T_{n,p,L,S}(x,y) \Big)^t,
\\
\hat F_{|e|,L,S}(y|x) &:=& F_{|e|}(y|x) + u_1^t\mathcal{M}(K)^{-1} \Big( T_{|e|,n,0,L,S}(x,y),\dots, T_{|e|,n,p,L,S}(x,y) \Big)^t,
\\
\hat F_{Y,L,S}^*(y|x) &:=& F_Y(y|x) + u_1^t\mathcal{M}(K)^{-1} \Big( T_{n,0,L,S}^*(x,y),\dots, T_{n,p,L,S}^*(x,y) \Big)^t,
\\
\hat F_{|e|,L,S}^*(y|x) &:=& F_{|e|}(y|x) + u_1^t\mathcal{M}(K)^{-1} \Big( T_{|e|,n,0,L,S}^*(x,y),\dots, T_{|e|,n,p,L,S}^*(x,y) \Big)^t,
\eean
$u_1^t := (1,0,...,0)$ denotes the first unit vector in $\R^{p+1}$, $\mathcal{M}(K)$ denotes a $(p+1)\times(p+1)$ matrix with entries 
\[
\mathcal{M}(K)_{ij} = \mu_{i+j-2}(K) := \int u^{i+j-2}K(u)du,
\]
and
\bean
T_{n,k,L,S}(x,y) &:=& \frac{1}{nh}\sum_{i=1}^n \frac{1}{f_X(X_i)}K_{h,k}(x-X_i) \Big( \Omega\Big(\frac{y-Y_i}{d_n} \Big) - F_Y(y|X_i)\Big),
\\
T_{|e|,n,k,L,S}(x,y) &:=& \frac{1}{nh}\sum_{i=1}^n \frac{1}{f_X(X_i)}K_{h,k}(x-X_i) \Big(\Omega\Big(\frac{y-|Y_i-\hat q_{\tau,L}(X_i)|}{d_n} \Big) - F_{|e|}(y|X_i)\Big),
\\
T_{n,k,L,S}^*(x,y) &:=& \frac{1}{nh}\sum_{i=1}^n \frac{1}{f_X(X_i)}K_{h,k}(x-X_i) \Big(\Omega\Big(\frac{y-Y_i^*}{d_n} \Big) - F_Y(y|X_i)\Big),
\\
T_{|e|,n,k,L,S}^*(x,y) &:=& \frac{1}{nh}\sum_{i=1}^n \frac{1}{f_X(X_i)}K_{h,k}(x-X_i) \Big(\Omega\Big(\frac{y-|Y_i^*-\hat q_{\tau,L}^*(X_i)|}{d_n} \Big) - F_{|e|}(y|X_i)\Big).
\eean
This, and further properties as differentiability and convergence rates of $\hat F_{Y,L,S}(y|x),\hat F_{|e|,L,S}, \hat F_{Y,L,S}^*$, $\hat F_{|e|,L,S}^*$ is the subject of Lemma \ref{lem:propfdach}.\\
The results in Lemma \ref{lem:genlin} and properties of the estimators $\hat F_Y, \hat F_{|e|}, \hat F_Y^*, \hat F_{|e|}^*$ yield representations of the form 
\bean
\hat q_\tau(x) &=& \hat q_{\tau,L}(x) + o_P(n^{-1/2}), \qquad \hat s(x) = \hat s_L(x) + o_P(n^{-1/2}),\\
\hat q_\tau^*(x) &=& \hat q_{\tau,L}^*(x) + o_P(n^{-1/2}), \qquad \hat s^*(x) = \hat s_L^*(x) + o_P(n^{-1/2})
\eean
uniformly in $x$ [see Lemma \ref{lem:quantlin}] where
\bean
\hat q_{\tau,L}(x) &:=& q_\tau(x) - \frac{1}{f_e(0|x)}\int_{-1}^1 \Big(\hat F_{Y,L,S}(q_{\tau+vb_n}(x)|x) -  F_Y(q_{\tau+vb_n}(x)|x)  \Big) \kappa(v) dv\\
&=& q_\tau(x) - \frac{u_1^t\mathcal{M}(K)^{-1}}{f_e(0|x)} \int_{-1}^1 \kappa(v)\Big( T_{n,0,L,S}(x,q_{\tau+vb_n}(x)),\dots, T_{n,p,L,S}(x,q_{\tau+vb_n}(x)) \Big)^t
dv
\\
\hat s_L(x) &:=& s(x) - \frac{1}{f_{|\eps|}(1|x)}\int_{-1}^1 \Big(\hat F_{|e|,L,S}(s_{1/2+vb_n}(x)|x) -  F_{|e|}(s_{1/2+vb_n}(x)|x)  \Big) \kappa(v) dv\\
&=& s(x) - \frac{u_1^t\mathcal{M}(K)^{-1}}{f_{|\eps|}(1)} \int_{-1}^1 \kappa(v)\Big( T_{|e|,n,0,L,S}(x,s_{1/2+vb_n}(x)),\dots, T_{|e|,n,p,L,S}(x,s_{1/2+vb_n}(x)) \Big)^t
dv
\\
\hat q_{\tau,L}^*(x) &:=& q_\tau(x) - \frac{1}{f_e(0|x)}\int_{-1}^1 \Big(\hat F_{Y,L,S}^*(q_{t+vb_n}(x)|x) -  F_Y(q_{t+vb_n}(x)|x)  \Big) \kappa(v) dv\\
&=& q_\tau(x) - \frac{u_1^t\mathcal{M}(K)^{-1}}{f_e(0|x)} \int_{-1}^1 \kappa(v)\Big( T_{n,0,L,S}^*(x,q_{\tau+vb_n}(x)),\dots, T_{n,p,L,S}^*(x,q_{\tau+vb_n}(x)) \Big)^t
dv
\\
\hat s_L^*(x) &:=& s(x) - \frac{1}{f_{|\eps|}(1)}\int_{-1}^1 \Big(\hat F_{|e|,L,S}^*(s_{1/2+vb_n}(x)|x) -  F_{|e|}(s_{1/2+vb_n}(x)|x)  \Big) \kappa(v) dv\\
&=& s(x) - \frac{u_1^t\mathcal{M}(K)^{-1}}{f_{|\eps|}(1)} \int_{-1}^1 \kappa(v)\Big( T_{|e|,n,0,L,S}^*(x,s_{1/2+vb_n}(x)),\dots, T_{|e|,n,p,L,S}^*(x,s_{1/2+vb_n}(x)) \Big)^t
dv
\eean
where $s_\alpha(x) := F_{|e|}^{-1}(\alpha|x)$. Differentiability properties and convergence rates of derivatives of these estimators can obviously be derived from the corresponding properties of the underlying distribution function estimators, see Lemma \ref{lem:unifratesquant}.

\begin{lemma}
\label{lem:unifratesquant}  
Let \ref{as:k1}-\ref{as:Gs}, \ref{as:fx}-\ref{as:fbound}, \ref{as:bw} hold. Then for any $k \leq 2$
\bean
\sup_{x\in [h_n,1-h_n]}|\hat q_{\tau,L}^{(k)}(x)-q_\tau^{(k)}(x)| \;=\; O_P\Big(\frac{\log h_n^{-1}}{nh_n(h_n\wedge d_n)^{2k}}\Big)^{1/2} &=& o_P(1),
\\
\sup_{x\in [2h_n,1-2h_n]}|\hat s_L^{(k)}(x)-s^{(k)}(x)| \;=\; O_P\Big(\frac{\log h_n^{-1}}{nh_n(h_n\wedge d_n)^{2k}}\Big)^{1/2} &=& o_P(1),
\eean
and under \ref{as:b0}-\ref{as:b2} it follows that
\bean
\sup_{x\in [3h_n,1-3h_n]}|(\hat q_{\tau,L}^*)^{(k)}(x)-q_\tau^{(k)}(x)| \;=\; O_P\Big(\frac{\log h_n^{-1}}{nh_n(h_n\wedge d_n)^{2k}}\Big)^{1/2} &=& o_P(1),
\\
\sup_{x\in [4h_n,1-4h_n]}| (\hat s_L^*)^{(k)}(x)-s^{(k)}(x)| \;=\; O_P\Big(\frac{\log h_n^{-1}}{nh_n(h_n\wedge d_n)^{2k}}\Big)^{1/2} = o_P(1).
\eean
\end{lemma}
{\bf Proof of Lemma \ref{lem:unifratesquant}} Since all claims share the same structure, we will only establish that
\[
\sup_{x\in [h_n,1-h_n]}|\hat q_{\tau,L}^{(k)}(x)-q_\tau^{(k)}(x)| \;=\;
 O_P\Big(\frac{\log h_n^{-1}}{nh_n(h_n\wedge d_n)^{2k}}\Big)^{1/2} = o_P(1).
\]
Observe that by definition of $\hat q_{\tau,L}$ we have
\bean
\hat q_{\tau,L}^{(k)}(x)-q_\tau^{(k)}(x)
= -\frac{\partial^k}{\partial x^k} \Big( \frac{1}{f_e(0|x)}\int_{-1}^1 \Big(\hat F_{Y,L,S}(q_{\tau+vb_n}(x)|x) -  F_Y(q_{\tau+vb_n}(x)|x)  \Big) \kappa(v) dv\Big).
\eean
Observing that $f_e(0|x) = f_\eps(0)/s(x)$, it suffices to show that
\[
\sup_{x \in [h_n,1-h_n] \atop v \in [-1,1]}\sup_{m \leq k} \Big|\frac{\partial^m}{\partial x^m}\Big(\hat F_{Y,L,S}(q_{\tau+vb_n}(x)|x) -  F_Y(q_{\tau+vb_n}(x)|x)  \Big)\Big| \;=\; O_P\Big(\frac{\log h_n^{-1}}{nh_n(h_n\wedge d_n)^{2k}}\Big)^{1/2}. 
\]
Now by Remark \ref{rem:qdiff} in the main body of the paper, the function $x \mapsto q_{\tau+vb_n}(x)$ is $2$ times continuously differentiable and its derivatives are bounded uniformly over $x\in (0,1), v \in [-1,1]$. Thus the above assertion follows from (i) of Lemma \ref{lem:propfdach} combined with the chain rule for derivatives. 
\hfill $\Box$

\begin{lemma} \label{lem:quantlin}
Let \ref{as:k1}-\ref{as:Gs}, \ref{as:fx}-\ref{as:fbound}, \ref{as:bw} hold. Then
\bean
&(i)& \sup_{x\in [h_n,1-h_n]}|\hat q_\tau(x) - \hat q_{\tau,L}(x)| = o_P(1/\sqrt{n}),
\\
&(ii)& \sup_{x\in [2h_n,1-2h_n]}|\hat s(x) - \hat s_L(x)| = o_P(1/\sqrt{n}),
\\
\eean
and if additionally \ref{as:b0}-\ref{as:b2} hold, we also have
\bean
&(iii)& \sup_{x\in [3h_n,1-3h_n]}|\hat q_\tau^*(x) - \hat q_{\tau,L}^*(x)| = o_P(1/\sqrt{n}),
\\
&(iv)& \sup_{x\in [4h_n,1-4h_n]}|\hat s^*(x) - \hat s_L^*(x)| = o_P(1/\sqrt{n}).
\eean
\end{lemma}
\textbf{Proof}
Since all assertions share a similar structure, we will only prove (iii). We begin by stating and intermediate result which we will establish in the end.
\begin{equation} \label{eq:lemc3_1}
\sup_{y \in\R}\sup_{x\in [3h_n,1-3h_n]} |\hat F_Y^*(y|x) - F_Y(y|x)| = o_P(1).
\end{equation}
Note that, in contrast to the statements in Lemma \ref{lem:propfdach} part (iii), the range for $y$ is $\R$ instead of a bounded set. Now let $\delta>0, c_0 >0$ be such that $\inf_{x \in [0,1]}\inf_{|y-q_\tau(x)|\leq 2\delta} f_Y(y|x) \geq c_0$ and define
\[
F_Y^*(y|x) := \hat F_Y^*(y|x)I\{|y-q_\tau(x)|\leq 2\delta/c_0\} + F_Y(y|x)I\{|y-q_\tau(x)| > 2\delta/c_0\}. 
\]
By the results in Lemma \ref{lem:propfdach} parts (iii), (iii)' we have
\begin{equation} \label{eq:lemc3_2}
\sup_{y \in\R}\sup_{x\in [3h_n,1-3h_n]} |F_Y^*(y|x) - F_Y(y|x)| = O_P\Big(\frac{\log n}{nh_n} \Big)^{1/2},
\end{equation}
and 
\begin{equation} \label{eq:lemc3_3}
\sup_{x\in [3h_n,1-3h_n]}\sup_{|y-q_\tau(x)|\leq 2\delta/c_0} |F_Y^*(y|x) - \hat F_{Y,L,S}^*(y|x)| = o_P(n^{-1/2}).
\end{equation}
Moreover, as we shall prove later, we have 
\begin{equation} \label{eq:lemc3_4}
P\Big( Q_{G,\kappa,\tau,b_n}(\hat F_Y^*(\cdot|x)) = Q_{G,\kappa,\tau,b_n}(F_Y^*(\cdot|x)) \ \forall x\in [3h_n,1-3h_n] \Big) \to 1.
\end{equation}
Now apply part (c) of Lemma \ref{lem:genlin} with $F = F_1 = F_Y(\cdot|x), F_2 = F_Y^*(\cdot|x)$. A careful inspection of the remainder terms in the statement of Lemma \ref{lem:genlin} part (c) shows that, uniformly in $x \in [3h_n,1-3h_n]$, 
\begin{eqnarray} \nonumber
&&Q_{G,\kappa,\tau,b_n}(F_Y^*(\cdot|x)) - Q_{G,\kappa,\tau,b_n}(F_Y(\cdot|x))
\\
&=& - \frac{1}{f_e(0|x)} \int_{-1}^{1} \kappa(v)\Big(F_Y^*(q_{\tau+vb_n}(x)|x) - F_Y(q_{\tau+vb_n}(x)|x)\Big)dv + o_P(n^{-1/2}). \label{eq:lemc3_5}
\end{eqnarray}
An application of Lemma \ref{lem:genlin}, part (a) with $F = F_Y(\cdot|x)$ shows that
\[
Q_{G,\kappa,\tau,b_n}(F_Y(\cdot|x)) = q_\tau(x) + O(b_n^2) = q_\tau(x) + o(n^{-1/2})
\]
uniformly in $x \in [0,1]$. Combining this with (\ref{eq:lemc3_3}), (\ref{eq:lemc3_4}) and (\ref{eq:lemc3_5}) and observing that $\hat q_\tau^*(x) = Q_{G,\kappa,\tau,b_n}(\hat F_Y^*(\cdot|x))$ we obtain, uniformly in $x \in [3h_n,1-3h_n]$, 
\begin{eqnarray*}
&&\hat q_\tau^*(x) - q_\tau(x)
\\
&=& - \frac{1}{f_e(0|x)} \int_{-1}^{1} \kappa(v)\Big(\hat F_{Y,L,S}^*(q_{\tau+vb_n}(x)|x) - F_Y(q_{\tau+vb_n}(x)|x)\Big)dv + o_P(n^{-1/2}).
\end{eqnarray*}
Note that, by the definition of $\hat q_{\tau,L}^*(x)$, the leading term in this representation is equal to $\hat q_{\tau,L}^*(x) - q_\tau(x)$. This implies statement (iii), and thus it remains to prove (\ref{eq:lemc3_1}) and (\ref{eq:lemc3_4}).\\

\textbf{Proof of (\ref{eq:lemc3_1})} Define (with $W_i$ the same as defined in (\ref{def-LL}))
\[
\hat F_{Y,U}^*(y|x) := \sum_{i=1}^n W_i(x) I\{ Y_i^* \leq y\}.
\] 
Since 
\[
\hat F_{Y}^*(y|x) = (\hat F_{Y,U}^*(\cdot|x)*\frac{1}{d_n}\omega(\cdot/d_n))(y) 
\] 
and by the smoothness of $F_Y$, it suffices to prove that
\beq \label{eq:lemc3_6}
\sup_{y \in\R}\sup_{x\in [3h_n,1-3h_n]} |\hat F_{Y,U}^*(y|x) - F_Y(y|x)| = o_P(1).
\eeq
Now by the definition of $Y_i^*$ we have
\[
\hat F_{Y,U}^*(y|x) = \sum_{i=1}^n W_i(x) I\{ \hat q_\tau(X_i) + \hat s(X_i)\tilde F_\eps^{-1}(U_i) \leq y\} = \sum_{i=1}^n W_i(x) I\Big\{ U_i \leq  \tilde F_\eps \Big(\frac{y - \hat q_\tau(X_i)}{\hat s(X_i)}\Big)\Big\}.
\]
From (\ref{le2n09}) in the main body of the paper we obtain after a Taylor expansion
\[
\sup_{x \in [3h_n,1-3h_n]}\sup_{y \in \R}\Big| \tilde F_\eps \Big(\frac{y - \hat q_\tau(x)}{\hat s(x)}\Big) - \tilde F_\eps \Big(\frac{y - q_\tau(x)}{s(x)}\Big) \Big| = o_P(1).
\]
Since the conclusion of Lemma 2 in Neumeyer (2009a) remains valid in our setting [see the discussion in the beginning of Section \ref{app-boot}], it follows that $\sup_{z \in \R} |\tilde F_\eps(z) - F_\eps(z)| = o_P(1)$ and thus
\[
\sup_{x \in [3h_n,1-3h_n]}\sup_{y \in \R}\Big| \tilde F_\eps \Big(\frac{y - \hat q_\tau(x)}{\hat s(x)}\Big) - F_\eps \Big(\frac{y - q_\tau(x)}{s(x)}\Big) \Big| = o_P(1).
\] 
Thus there exists a deterministic sequence $\gamma_n \to 0$ such that $P(D_n) \to 1$ where we defined the event
\[
D_n := \Big\{ \sup_{x \in [3h_n,1-3h_n]}\sup_{y \in \R}\Big| \tilde F_\eps \Big(\frac{y - \hat q_\tau(x)}{\hat s(x)}\Big) - F_\eps \Big(\frac{y - q_\tau(x)}{s(x)}\Big) \Big| \leq \gamma_n\Big\}.
\]
Additionally, define the event 
\[
\tilde D_n := \{\sup_i \sup_{x \in [h_n,1-h_n]}|W_i(x)| \leq C (nh_n)^{-1}I\{|x-X_i|\leq h_n\}\}
\] 
and observe that $P(\tilde D_n) \to 1$ by the definition of $W_i(x)$ and Lemma \ref{lem:propweight}. Thus on $D_n \cap\tilde D_n$ we have
\begin{eqnarray} \nonumber
\quad &&\sup_{y \in\R}\sup_{x\in [3h_n,1-3h_n]} \Big|\hat F_{Y,U}^*(y|x) -\sum_{i=1}^n W_i(x) I\Big\{ U_i \leq  F_Y(y|X_i)\Big\} \Big|
\\
\quad &\leq& \frac{C}{nh_n}\sup_{y \in \R}\sup_{x\in [3h_n,1-3h_n]} \sum_{i=1}^n I\{|X_i-x|\leq h_n\}I\Big\{\Big|U_i-F_Y(y|X_i)\Big|\leq \gamma_n\Big\} = o_P(1) \label{eq:lemc3_7}
\end{eqnarray}  
where the last equality follows by a combination of parts 1, 4-6 of Lemma \ref{lem:entrnum} with Lemma \ref{lem:base}. Similarly, applying Lemma \ref{lem:propweight}, parts 1,2 4-6 of Lemma \ref{lem:entrnum} with Lemma \ref{lem:base} shows that
\beq \label{eq:lemc3_8}
\sum_{i=1}^n W_i(x) \Big( I\Big\{ U_i \leq  F_Y(y|X_i)\Big\} - F_Y(y|X_i)\Big) = o_P(1)
\eeq
uniformly in $x \in [3h_n,1-3h_n], y \in \R$. Finally, by similar arguments as used in the proof of (\ref{eq:help1}) one can show that
\beq \label{eq:lemc3_9}
\sum_{i=1}^n W_i(x) F_Y(y|X_i) = F_{Y}(y|x) + o_P(1)
\eeq
uniformly in $x \in [3h_n,1-3h_n], y \in \R$. Combining (\ref{eq:lemc3_7})-(\ref{eq:lemc3_9}) yields (\ref{eq:lemc3_6}) and completes the proof of (\ref{eq:lemc3_1}).\\

\textbf{Proof of (\ref{eq:lemc3_4})} Define the events 
\bean
D_{n1} &:=& \Big\{ \hat F_Y^*(y|x) = F_Y^*(y|x)\ \forall (x,y) \in \{(x,y): |F_Y^*(y|x) - \tau| \leq \delta, x \in [3h_n,1-3h_n]\} \Big\}
\\
D_{n2} &:=& \Big\{\sup_{x \in [3h_n,1-3h_n], y \in \R}|\hat F_Y^*(y|x) - F_Y^*(y|x)| \leq \delta/2\Big\}
\\
D_{n3} &:=& \Big\{\sup_{x \in [3h_n,1-3h_n], y \in \R}|\hat F_Y^*(y|x) - F_Y(y|x)| \leq \delta/2\Big\}.
\eean 
Observe that on $D_{n1} \cap D_{n2} \cap D_{n3}$ we have $F_Y^*(y|x) \leq \tau - \delta \Rightarrow \hat F_Y^*(y|x) \leq \tau - \delta/2$, $F_Y^*(y|x) \geq \tau + \delta \Rightarrow \hat F_Y^*(y|x) \geq \tau + \delta/2$ and $|F_Y^*(y|x) - \tau| \leq \delta \Rightarrow F_Y^*(y|x) = \hat F_Y^*(y|x)$. Thus on $D_{n1} \cap D_{n2} \cap D_{n3}$ we obtain $Q_{G,\kappa,\tau,b_n}(\hat F_Y^*(\cdot|x)) = Q_{G,\kappa,\tau,b_n}(F_Y^*(\cdot|x))$ provided that $b_n \leq \delta/2$. It remains to prove that $P(D_{n1} \cap D_{n2} \cap D_{n3}) \to 1$. The fact that $P(D_{n2} \cap D_{n3}) \to 1$ follows from (\ref{eq:lemc3_1}), (\ref{eq:lemc3_2}), so that it remains to prove $P(D_{n1}) \to 1$ which follows from
\bean
&&P\Big( \{(x,y): |F_Y^*(y|x) - \tau| \leq \delta, x \in [3h_n,1-3h_n]\} \subset
\\
&&\quad\quad\quad\quad\quad\quad\quad\quad\quad\quad\{(x,y): |y - q_\tau(x)| \leq 2\delta/c_0, x \in [3h_n,1-3h_n]\}\Big) \to 1.
\eean
This in turn is a consequence of the fact that on $D_{n3}$ (note that $|F_Y^*(y|x)-F_Y(y|x)| \leq |\hat F_Y^*(y|x)-F_Y(y|x)|$)
\[
|F_Y^*(y|x) - \tau| \leq \delta \Rightarrow |F_Y(y|x) - \tau| \leq 3\delta/2 \Rightarrow |y- q_\tau(x)| \leq 3\delta/(2c_0)
\]
by the definition of $\delta,c_0$. This completes the proof.

\hfill $\Box$

\newpage

\begin{lemma}\label{lem:propfdach}
Assume that conditions \ref{as:k1}-\ref{as:Gs}, \ref{as:fx}-\ref{as:fbound} and \ref{as:bw} hold. Denote by $\tilde T_{n,0,L,S}, \tilde T_{|e|,n,0,L,S}, \tilde T_{n,0,L,S}^*, \tilde T_{|e|,n,0,L,S}^*$ versions of $T_{n,0,L,S}, T_{|e|,n,0,L,S}, T_{n,0,L,S}^*, T_{|e|,n,0,L,S}^*$ where $1/f_X(X_i)$ is replaced by $1/f_X(x)$.\\
Then for any bounded $\Y_1\subset\R,\Y_2 \subset \R^+$ such that $\Y_2$ is bounded away from zero we have
\bean
(i)'\quad \hat F_Y(y|x) &=& \hat F_{Y,L,S}(y|x) + o_P(1/\sqrt{n}), \quad T_{n,0,L,S} = \tilde T_{n,0,L,S} + o_P(1/\sqrt{n}),
\eean
uniformly in $y\in\Y_1, x\in [h_n,1-h_n]$ and 
\bean
(ii)'\quad \hat F_{|e}|(y|x) &=& \hat F_{|e|,L,S}(y|x) + o_P(1/\sqrt{n}), \quad T_{|e|,n,0,L,S} = \tilde T_{|e|,n,0,L,S} + o_P(1/\sqrt{n}),
\eean
uniformly in $y\in\Y_2, x\in [2h_n,1-2h_n]$. If additionally \ref{as:b0}-\ref{as:b2} hold,
\bean
(iii)'\quad \hat F_Y^*(y|x) &=& \tilde T_{n,0,L,S}^* + o_P(1/\sqrt{n}), \quad T_{n,0,L,S}^* = \tilde T_{n,0,L,S}^* + o_P(1/\sqrt{n}), 
\eean
uniformly in $y\in\Y_1, x\in [3h_n,1-3h_n]$ and
\bean
(iv)'\quad \hat F_{|e|}^*(y|x) &=& \hat F_{|e|,L,S}^*(y|x) + o_P(1/\sqrt{n}), \quad T_{|e|,n,0,L,S}^* = \tilde T_{|e|,n,0,L,S}^* + o_P(1/\sqrt{n}).
\eean
uniformly in $y\in\Y_2, x\in [4h_n,1-4h_n]$.\\

Moreover, (i)-(iv) hold under the assumptions of $(i)'-(iv)'$, respectively.
\bean
&(i)& \forall k+l \leq 2 \quad \sup_{y \in \mathcal{Y}_1, x \in [h_n,1-h_n]} |\partial_x^k\partial_y^l\hat F_{Y,L,S}(y|x) - \partial_x^k\partial_y^lF_Y(y|x)| = O_P\Big(\frac{\log n}{nh_n^{2k+1}d_n^{2l}} \Big)^{1/2},
\\
&(ii)& \forall k+l \leq 2 \quad \sup_{y \in \mathcal{Y}_2, x \in [2h_n,1-2h_n]} |\partial_x^k\partial_y^l\hat F_{|e|,L,S}(y|x) - \partial_x^k\partial_y^l F_{|e|}(y|x)| =  O_P\Big(\frac{\log n}{nh_n^{2k+1}d_n^{2l}} \Big)^{1/2},
\\
&(iii)& \forall k+l \leq 2 \quad \sup_{y \in \mathcal{Y}_1, x \in [3h_n,1-3h_n]} |\partial_x^k\partial_y^l\hat F_{Y,L,S}^*(y|x) - \partial_x^k\partial_y^l F_Y(y|x)| =  O_P\Big(\frac{\log n}{nh_n^{2k+1}d_n^{2l}} \Big)^{1/2},
\\
&(iv)& \forall k+l \leq 2 \quad \sup_{y \in \mathcal{Y}_2, x \in [4h_n,1-4h_n]} |\partial_x^k\partial_y^l\hat F_{|e|,L,S}^*(y|x) - \partial_x^k\partial_y^l F_{|e|}(y|x)| =  O_P\Big(\frac{\log n}{nh_n^{2k+1}d_n^{2l}} \Big)^{1/2}.
\eean

\end{lemma}

\textbf{Proof of Lemma \ref{lem:propfdach}}\\
We will only provide the arguments for (iv) and (iv)' since all other assertions can be derived analogously. Since $\Y_2$ is bounded away from zero, and since $d_n \to 0$, the fact that $\omega = \Omega'$ is symmetric and has support $[-1,1]$ implies that for $n$ sufficiently large
\[
\Omega\Big(\frac{y-|z|}{d_n} \Big) = \Omega\Big(\frac{y-z}{d_n} \Big) - \Omega\Big(\frac{-y - z}{d_n} \Big) \quad \forall y \in \Y_2,\ z \in \R.
\]

Thus we find that for $n$ sufficiently large
\bean
\hat F_{|e|}^*(y|x) &=& \hat F_{e}^*(y|x) - \hat F_{e}^*(-y|x),
\\
\hat F_{|e|,L,S}^*(y|x) &=& \hat F_{e,L,S}^*(y|x) - \hat F_{e,L,S}^*(-y|x),
\\
T_{|e|,n,0,L,S}^*(x,y) &=& T_{e,n,0,L,S}^*(x,y) - T_{e,n,0,L,S}^*(x,-y),
\\
\tilde T_{|e|,n,0,L,S}^* &=& \tilde T_{e,n,0,L,S}^*(x,y) - \tilde T_{e,n,0,L,S}^*(x,-y),
\eean
where
\bean
\hat F_{e}^*(y|x) &:=& \sum_i W_i(x) \Omega\Big(\frac{y-(Y_i^*-\hat q_{\tau}^*(X_i))}{d_n} \Big),
\\
\hat F_{e,L,S}^*(y|x) &:=& F_{e}(y|x) + u_1^t\mathcal{M}(K)^{-1} \Big( T_{e,n,0,L,S}^*(x,y),\dots, T_{e,n,p,L,S}^*(x,y) \Big)^t,
\\
T_{e,n,0,L,S}^*(x,y) &:=& \frac{1}{nh_n}\sum_{i=1}^n \frac{1}{f_X(X_i)}K_{h_n,k}(x-X_i) \Big(\Omega\Big(\frac{y-(Y_i^*-\hat q_{\tau,L}^*(X_i))}{d_n} \Big) - F_e(y|X_i)\Big),
\\
\tilde T_{e,n,0,L,S}^* &:=& \frac{1}{nh_n}\sum_{i=1}^n \frac{1}{f_X(x)}K_{h_n,k}(x-X_i) \Big(\Omega\Big(\frac{y-(Y_i^*-\hat q_{\tau,L}^*(X_i))}{d_n} \Big) - F_e(y|X_i)\Big).
\eean

It thus suffices to establish, uniformly in $y \in \mathcal{Y}:= \Y_2 \cup (-\Y_2), x \in [4h_n,1-4h_n]$,
\bea \label{eq:iv1}
&& \quad\quad \hat F_{e}^*(y|x) = F_{e}(y|x) + u_1^t\mathcal{M}(K)^{-1}
\Big(  T_{e,n,0,L,S}^*(x,y),\dots, T_{e,n,p,L,S}^*(x,y) \Big)^t + o_P(n^{-1/2}),
\\ \label{eq:iv2}
&& \quad\quad T_{e,n,0,L,S}^* = \tilde T_{e,n,0,L,S}^* + o_P(n^{-1/2}),
\\ \label{eq:iv'}
&& \quad\quad\sup_{y \in \mathcal{Y}_2, x \in [4h_n,1-4h_n]} |\partial_x^k\partial_y^l\hat F_{e,L,S}^*(y|x) - \partial_x^k\partial_y^l F_{e}(y|x)| =  O_P\Big(\frac{\log n}{nh_n^{2k+1}d_n^{2l}} \Big)^{1/2}.
\eea

Define the quantities
\begin{eqnarray*}
T_{e,n,k,L}^*(x,y) &:=& \frac{1}{nh_n}\sum_{i=1}^n \frac{1}{f_X(X_i)}K_{h_n,k}(x-X_i) \Big(I\{Y_i^* \leq y+\hat q_{\tau,L}^*(X_i)\} - F_e(y|X_i)\Big),
\\
\tilde T_{e,n,k,L}^*(x,y) &:=& \frac{1}{nh_n}\sum_{i=1}^n \frac{1}{f_X(x)}K_{h_n,k}(x-X_i) \Big(I\{Y_i^* \leq y+\hat q_{\tau,L}^*(X_i)\} - F_e(y|X_i)\Big),
\end{eqnarray*}
and note that,uniformly in $y \in \Y, x \in [4h_n,1-4h_n]$, 
\bea \label{eq:help00a}
(T_{e,n,k,L}^*(x,\cdot)*\frac{1}{d_n}\omega(\cdot/d_n))(y) &=& T_{e,n,k,L,S}^*(x,y) + o(1/\sqrt n),
\\ \label{eq:help00b}
(\tilde T_{e,n,k,L}^*(x,\cdot)*\frac{1}{d_n}\omega(\cdot/d_n))(y) &=& \tilde T_{e,n,k,L,S}^*(x,y) + o(1/\sqrt n).
\eea

Also, let
\bean
\hat F_{e,U}^*(y|x) &:=& \sum_{i=1}^n W_i(x) I\{ Y_i^* - \hat q_{\tau}^*(X_i) \leq y\}
\\
&=& \frac{1}{nh_n} u_1^t (\XX^t\WW\XX)^{-1}\left(
\begin{array}{c}
\sum_i K_{h_n,0}(x-X_i)I\{ Y_i^* - \hat q_{\tau}^*(X_i) \leq y\}\\
\vdots \\
\sum_i h_n^p K_{h_n,p}(x-X_i)I\{ Y_i^* - \hat q_{\tau}^*(X_i) \leq y\}
\end{array}
\right),
\\
\hat F_{e,L,U}^*(y|x) &:=& F_e(y|x) + u_1^t\mathcal{M}(K)^{-1}
\Big( T_{e,n,0,L}^*(x,y),\dots, T_{e,n,p,L}^*(x,y) \Big)^t
\eean

where the weights $W_i(x)$ are the same as in equation (\ref{def-LL}). At the end of the proof, we will establish the following assertions uniformly in $y \in \Y, x \in [4h_n,1-4h_n]$
\begin{eqnarray}
\label{eq:help2}
\quad \quad T_{e,n,0,L}^*(x,y) &=& \tilde T_{e,n,0,L}^*(x,y) + o_P(n^{-1/2}).
\\
\label{eq:help1}
\quad \quad \hat F_{e,U}^*(y|x) &=& \hat F_{e,L,U}^*(y|x) + o_P(n^{-1/2}),
\\ 
\label{eq:help6}
\partial_x^m T_{e,n,k,L}^*(x,y) &=&  
O_P\Big(\frac{\log n}{nh_n^{2m+1}} \Big)^{1/2}, \quad m = 0,1,2.
\end{eqnarray}
Now assertions (\ref{eq:iv1}), (\ref{eq:iv2}) follows from (\ref{eq:help00a}), (\ref{eq:help1}) and (\ref{eq:help2}) since 
\begin{eqnarray*}
(\hat F_{e,U}^*(\cdot|x)*\frac{1}{d_n}\omega(\cdot/d_n))(y) &=& \hat F_{e}^*(y|x),
\\
(F_e(\cdot|x)*\frac{1}{d_n}\omega(\cdot/d_n))(y) &=& F_e(y|x) + O(d_n^{p_\omega}) = F_e(y|x) + o(1/\sqrt n),
\end{eqnarray*}
uniformly in $x\in[4h_n,1-4h_n],y\in\Y$.

On the other hand we have
\[
\partial_x^m \hat F_{e,L,U}^*(y|x) := \partial_x^m F_e(y|x) + u_1^t\mathcal{M}(K)^{-1}
\Big( \partial_x^m T_{e,n,0,L}^*(x,y),\dots, \partial_x^m T_{e,n,p,L}^*(x,y) \Big)^t,
\]
and thus (\ref{eq:help6}) implies, uniformly in $y \in \Y, x \in [4h_n,1-4h_n]$, 
\[
\partial_x^m \hat F_{e,L,U}^*(y|x) = \partial_x^m F_e(y|x) + O_P\Big(\frac{\log n}{nh_n^{2m+1}} \Big)^{1/2}.
\]
This entails (\ref{eq:iv'}) since 
\bean
\partial_x^k \partial_y^l \Big(\hat F_{e,L,S}^*(y|x)-F_e(y|x) \Big) 
&=& \frac{1}{d_n^{l}}\Big[\Big(\partial_x^k\hat F_{e,L,U}^*(\cdot|x)-\partial_x^kF_e(\cdot|x)\Big) *\Big(\frac{1}{d_n}\omega^{(l)}\Big(\frac{\cdot}{d_n}\Big)\Big) \Big](y)
\\&& + \Big((\partial_x^k \partial_y^l F_e(\cdot|x))*\Big(\frac{1}{d_n}\omega\Big(\frac{\cdot}{d_n}\Big)\Big)\Big)(y) - \partial_x^k \partial_y^l F_e(y|x).
\eean
Now, since by assumption $\partial_x^k F_e(y|x)$ is $r$ times continuously differentiable with respect to $y$, the second summand is of order $d_n^{r-l} = O\Big(\frac{\log n}{nh_n^{2k+1}d_n^{2l}} \Big)^{1/2}$. The first summand can be bounded by $\frac{1}{d_n^l} O_P\Big(\frac{\log n}{nh_n^{2k+1}} \Big)^{1/2}$.\\ 

The proof will thus be complete after we establish (\ref{eq:help2})-(\ref{eq:help6}). In order to do so, observe that there exists a set $D_n$ such that the probability of $D_n$ tends to one and such that on $D_n$ we have, for any sequence $c_n$ such that $c_n/r_n \to \infty$ [this is a consequence of (\ref{le2n09}) and the uniform rates of convergence for $\hat s_L,\hat q_{\tau,L}, \hat q_{\tau,L}^*$ which follow from parts (i)-(iii) of Lemma \ref{lem:quantlin} and Lemma \ref{lem:unifratesquant}]
\bean
&&\Big|\tilde F_\eps\Big(\frac{y}{\hat s_L(X_i)} + \frac{\hat q_{\tau,L}^*(X_i)-\hat q_{\tau,L}(X_i)}{\hat s_L(X_i)} \Big) - F_\eps\Big( \frac{y}{s(X_i)}\Big)\Big|
\\
&\leq& \sup_{y \in (1+ \Y/c_s)}|\tilde F_\eps(y) - F_\eps(y)| + 0.5c_n \sup_{y \in (1+ \Y/c_s)} |y f_\eps(y)| \leq c_n
\eean 
where the last bound follows from (\ref{eq:c72}) in Lemma \ref{lem:propfepsdach}. In particular, on $D_n$ we have
\bea \label{eq:help4}
I\Big\{U_i \leq  F_\eps\Big( \frac{y}{s(X_i)}\Big) - c_n\Big\} 
&\leq& I\Big\{U_i \leq \tilde F_\eps\Big(\frac{y}{\hat s_L(X_i)} + \frac{\hat q_{\tau,L}^*(X_i)-\hat q_{\tau,L}(X_i)}{\hat s_L(X_i)} \Big)\Big\} 
\\ \label{eq:help5}
&\leq&
I\Big\{U_i \leq  F_\eps\Big( \frac{y}{s(X_i)}\Big) + c_n\Big\}.
\eea


\textbf{Proof of (\ref{eq:help2})}

Recall that $Y_i^* = \hat q_\tau(X_i) + \hat s(X_i) \eps_i^*$ and $\eps_i^* = \tilde F_\eps^{-1}(U_i)$. Observe the identity
\[
I\{Y_i^* \leq y + \hat q_{\tau,L}^*(X_i)\} = I\Big\{U_i \leq \tilde F_\eps\Big(\frac{y}{\hat s(X_i)} + \frac{\hat q_{\tau,L}^*(X_i)-\hat q(X_i)}{\hat s(X_i)} \Big)\Big\}.
\]
Moreover, a Taylor expansion shows that, with probability tending to one, 
\bean
&&\Big|
I\Big\{U_i \leq \tilde F_\eps\Big(\frac{y}{\hat s(X_i)} + \frac{\hat q_{\tau,L}^*(X_i)-\hat q(X_i)}{\hat s(X_i)} \Big)\Big\}-
I\Big\{U_i \leq \tilde F_\eps\Big(\frac{y}{\hat s_L(X_i)} + \frac{\hat q_{\tau,L}^*(X_i)-\hat q_{\tau,L}(X_i)}{\hat s_L(X_i)} \Big)\Big\} 
\Big| 
\\
&\leq& I\Big\{ \Big|U_i - \tilde F_\eps\Big(\frac{y}{\hat s_L(X_i)} + \frac{\hat q_{\tau,L}^*(X_i)-\hat q_{\tau,L}(X_i)}{\hat s_L(X_i)} \Big)\Big| \leq C\gamma_n \sup_{y \in 2\Y/c_s }|y\tilde f_\eps(y)|\Big\} 
\eean
where $\gamma_n = o(1/\sqrt{n})$, and thus arguments similar to those in the proof of Lemma \ref{lem:proclin} yield
\bean
&&\frac{1}{n}\sum_{i=1}^n\frac{K_{h_n,k}(x-u)}{h_n}\Big(\frac{1}{f_X(u)} - \frac{1}{f_X(x)} \Big)
\Big( I\Big\{U_i \leq \tilde F_\eps\Big(\frac{y}{\hat s(X_i)} + \frac{\hat q_{\tau,L}^*(X_i)-\hat q(X_i)}{\hat s(X_i)} \Big)\Big\}\Big)
\\
&=& \frac{1}{n}\sum_{i=1}^n\frac{K_{h_n,k}(x-u)}{h_n}\Big(\frac{1}{f_X(u)} - \frac{1}{f_X(x)} \Big)
\Big( I\Big\{U_i \leq \tilde F_\eps\Big(\frac{y}{\hat s_L(X_i)} + \frac{\hat q_{\tau,L}^*(X_i)-\hat q_{\tau,L}(X_i)}{\hat s_L(X_i)} \Big)\Big\}\Big)
\\
&& + o_P(1/\sqrt n).
\eean
Next, observe that  by part (i)-(iii) of Lemma \ref{lem:quantlin}, Lemma \ref{lem:unifratesquant} and by (\ref{le2n09}) there exists a set $D_n$ whose probability tends to one such that on $D_n$ we have for some $\delta > 0$
\bean
&&\hat s_L \in \tilde C_C^{1+\delta}([3h_n,1-3h_n]), \quad
\hat q_{\tau,L}^*,\hat q_{\tau,L} \in C_C^{1+\delta}([3h_n,1-3h_n]),
\\
\\
&&\sup_{u \in [3h_n,1-3h_n], y \in \Y} \Big| \tilde F_\eps\Big(\frac{y}{\hat s_L(u)} + \frac{\hat q_{\tau,L}^*(u)-\hat q_{\tau,L}(u)}{\hat s_L(u)} \Big) - F_{e}(y|u) \Big| \leq r_nh_n^{-1/4},
\eean
and $\tilde F_\eps \in \mathcal{D}$ defined in (\ref{cal-D}). 
Additionally, (\ref{cal-D-cov}) and the arguments from Proposition 3 in Neumeyer (2009a) show that for the class of functions
\bean
\G_{n} &:=& \Big\{ (u,v) \mapsto I\Big\{ u \leq F\Big(\frac{y}{a_1(v)} + \frac{ a_2(v)}{a_1(v)}\Big)\Big\}
\\
&& \quad\quad\quad  \Big|\  y\in \Y, F \in \mathcal{D}, a_1 \in \tilde C_C^{1+\delta}([3h_n,1-3h_n]), a_2 \in C_C^{1+\delta}([3h_n,1-3h_n])  \Big\}
\eean
we have, denoting by $P$ the product measure of the uniform random variable $U_1$ and the covariate $X_1$, $\sup_n \log N_{[\ ]}(\eps,\G_{n}, L^2(P)) \leq C\eps^{-2\alpha}$ for some $\alpha < 1$. Next, define the class of functions
\bean
\mathcal F_n &:=& \Big\{ (u,v) \mapsto \frac{K_{h_n,k}(x-u)}{h_n}\Big(\frac{1}{f_X(u)} - \frac{1}{f_X(x)} \Big)\times
\\
&& \qquad
\times\Big(I\Big\{v \leq \tilde F_\eps\Big(\frac{y}{\hat s_L(u)} + \frac{\hat q_{\tau,L}^*(u)-\hat q_{\tau,L}(u)}{\hat s_L(u)} \Big)\Big\} - F_e(y|u)\Big) \Big|x\in [4h_n,1-4h_n], y \in \mathcal Y \Big\}.
\eean
In particular, observe that,  due to the continuous differentiability of $f_X$ and the compact support of $K$, the functions in $\mathcal F_n$ are bounded uniformly over $n$. Additionally, combining the bound on $\sup_n \log N_{[\ ]}(\eps,\G_{n}, L^2(P))$ with parts 1, 3 and 4 of Lemma \ref{lem:entrnum}, we find that on $D_n$
\[
\sup_n \log N_{[\ ]}(\eps,\mathcal F_n, L^2(P)) \leq \tilde C\eps^{-2\tilde \alpha}
\]
for some $\tilde \alpha < 1$ and finite $\tilde C$. Moreover, again on $D_n$, we find that for each $f \in \mathcal F_n$
\[
\E f(X_i,U_i) = O(h_n^{3/4}r_n) = o(1/\sqrt n),\quad \E f^2(X_i,U_i) = O(h_n). 
\]
To see the second statement, observe that every $f \in \mathcal F_n$ satisfies
\[
|f(X_i,U_i)| \leq 2 \Big| \frac{K_{h_n,k}(x-X_i)}{h_n}\Big(\frac{1}{f_X(X_i)} - \frac{1}{f_X(x)}\Big) \Big|,
\]
the assertion now follows from a Taylor expansion of $f_X$. For the bound on $\E f(X_i,U_i) $, observe that
\bean
| \E [f(X_i,U_i)]| &\leq& \int \Big| \frac{K_{h_n,k}(x-u)}{h_n}\Big(\frac{1}{f_X(u)} - \frac{1}{f_X(x)}\Big) \Big| r_n h_n^{-1/4} f_X(u) du,
\eean
the claimed bound now follows from a Taylor expansion of $1/f_X(u)$ around $x$. 
Thus by Lemma \ref{lem:base} $\sup_{f \in \mathcal F_n}|\sum_i f(X_i,U_i)| = o_P(1/\sqrt n)$ and (\ref{eq:help2}) follows.

\newpage

\textbf{Proof of (\ref{eq:help1})} 
Define $\mathcal{H} := \mbox{diag}(1,h_n,...,h_n^p)$ and observe that by (\ref{eq:help2}) we have uniformly in $x\in [4h_n,1-4h_n], y \in \Y$ 
\bean
\hat F_{e,U}^*(y|x) - \hat F_{e,L,U}^*(y|x)
&=& \sum_{i} W_i(x)(I\{Y_i^* - \hat q_\tau^*(X_i) \leq y\}-I\{Y_i^* - \hat q_{\tau,L}^*(X_i) \leq y\})
\\ 
&& + \frac{u_1^t(\XX^t\WW\XX)^{-1} }{nh_n}\mathcal{H}
\left(
\begin{array}{c}
\sum_i K_{h_n,0}(x-X_i)F_e(y|X_i)\\
\vdots \\
\sum_i K_{h_n,p}(x-X_i)F_e(y|X_i)
\end{array}
\right)
- F_e(y|x)
\\
&&+ \Big(u_1^t(\XX^t\WW\XX)^{-1}\mathcal{H} - \frac{u_1^t\mathcal{M}(K)^{-1}}{f_X(x)}\Big) 
\left(
\begin{array}{c}
f_X(x)\tilde T_{e,n,0,L}^*(x,y)\\
\vdots \\
f_X(x)\tilde T_{e,n,p,L}^*(x,y)
\end{array}
\right)
\\
&& + u_1^t\mathcal{M}(K)^{-1}\left(
\begin{array}{c}
\tilde T_{e,n,0,L}^*(x,y) - T_{e,n,0,L}^*(x,y)\\
\vdots \\
\tilde T_{e,n,p,L}^*(x,y) - T_{e,n,p,L}^*(x,y)
\end{array}
\right)
\\ \vspace{.2cm}
&=:& R_{n,1}(x,y) + R_{n,2}(x,y) + R_{n,3}(x,y) + R_{n,4}(x,y).
\eean

Note that a Taylor expansion of $F_e(y|X_i)$ with respect to $X_i$ around the point $x$ combined with the fact that
\[
\frac{1}{nh_n}u_1^t(\XX^t\WW\XX)^{-1}
\left(
\begin{array}{c}
h_n^k\sum_i K_{h_n,k}(x-X_i)\\
\vdots \\
h_n^{p+k}\sum_i K_{h_n,p+k}(x-X_i)
\end{array}
\right) 
= I\{k=0\}
\]
for $k=0,...,p$ yields the representation
\[
\frac{u_1^t(\XX^t\WW\XX)^{-1} }{nh_n}
\left(
\begin{array}{c}
\sum_i K_{h_n,0}(x-X_i)F_e(y|X_i)\\
\vdots \\
\sum_i h_n^p K_{h_n,p}(x-X_i)F_e(y|X_i)
\end{array}
\right)
= F_e(y|x) + O_P(h_n^{p+1}) = F_e(y|x) + o_P(n^{-1/2})
\]
uniformly in $x\in[4h_n,1-4h_n],y\in\mathcal{Y}$, so that $R_{n,2}$ is small.\\ 

Next, consider $R_{n,3}$. By Lemma \ref{lem:propweight} and observing that $u_1^t \mathcal{H}^{-1} = u_1^t$ we find
\[
\Big(\frac{u_1^t(\XX^t\WW\XX)^{-1} }{nh_n}\mathcal{H} - \frac{u_1^t\mathcal{M}(K)^{-1}}{f_X(x)}\Big) = O_P(h_n),
\]
and together with the fact that
\[
\sup_{x [4h_n,1-4h_n], y\in \Y}\sup_{k=0,...,p} |T_{e,n,k,L}^*(x,y)| = O_P\Big(\frac{\log n}{nh_n} \Big)^{1/2}
\] 
which follows by similar arguments as the proof of (\ref{eq:help6}), this shows that $R_{n,3}$ is small.\\ 

The negligibility of $R_{n,4}$ follows from (\ref{eq:help2}).\\ 

Finally, consider $R_{n,1}$. Observe that, by similar arguments as in the proof of (\ref{eq:help2}), there exists a deterministic sequence $\xi_n = o(n^{-1/2})$ such that, with probability tending to one, we have for any $X_i \in [3h_n,1-3h_n]$ 
\[
\Big| I\{Y_i^* - \hat q_\tau^*(X_i) \leq y\}-I\{Y_i^* - \hat q_{\tau,L}^*(X_i) \leq y\} \Big|\leq I\Big\{\Big|U_i - \tilde F_\eps\Big(\frac{y}{\hat s_L(X_i)} + \frac{\hat q_{\tau,L}^*(X_i)-\hat q_{\tau,L}(X_i)}{\hat s_L(X_i)} \Big)\Big| \leq \xi_n \Big\}.
\]
Now arguments similar to those in the proof of Lemma \ref{lem:proclin} yield for every $k = 0,...,p$
\[
d_{n,k} := \frac{1}{n}\sum_{i=1}^n \frac{|K_{h_n,k}(x-u)|}{h_n} \frac{1}{f_X(x)}I\Big\{\Big|U_i - \tilde F_\eps\Big(\frac{y}{\hat s_L(X_i)} + \frac{\hat q_{\tau,L}^*(X_i)-\hat q_{\tau,L}(X_i)}{\hat s_L(X_i)} \Big)\Big| \leq \xi_n \Big\}
= o_P(n^{-1/2})
\]
uniformly over $x \in [4h_n,1-4h_n], y \in \Y$. Moreover, by Lemma \ref{lem:propweight} we have
\[
|R_{n,1}(x,y)| \leq (p+1) \Big(\max_{k=0,...,p} (u_1^t (\XX^t\WW\XX)^{-1}\mathcal H)_k\Big)\Big(\max_{k=0,...,p} |d_{n,k}(x,y)|\Big)
\]
This shows that $R_{n,1}$ is negligible and completes the proof of (\ref{eq:help1}).

\textbf{Proof of (\ref{eq:help6})}
Consider the decomposition 
\[
\partial_x^m T_{e,n,k,L}^*(x,y) = A_{n,k,m}^{+}(x,y) + A_{n,k,m}^{-}(x,y)
\]
where
\[
A_{n,k,m}^{+}(x,y) := \frac{1}{nh_n}\frac{1}{h_n^m}\sum_{i=1}^n \frac{K_{h_n,k}^{(m)}(x-X_i)}{f_X(X_i)}  I\Big\{K_{h_n,k}^{(m)}(x-X_i)>0\Big\} \Big(I\{Y_i^* \leq y + \hat q_{\tau,L}^*(X_i)\} - F_e(y|X_i)\Big)
\]
and $A_{n,k,m}^{-}$ is defined analogously. On the set $D_n$ (defined in the beginning of this proof) we have
\bean
A_{n,k,m}^{+}(x,y) 
&\leq& \frac{1}{nh_n^{m+1}}\sum_{i=1}^n \frac{K_{h_n,k}^{(m)}(x-X_i)}{f_X(X_i)}  I\Big\{K_{h_n,k}^{(m)}(x-X_i)>0\Big\}\times
\\
&&\qquad\qquad\qquad\qquad\qquad\qquad
\times \Big(I\Big\{U_i \leq  F_\eps\Big( \frac{y}{s(X_i)}\Big) + c_n\Big\}  - F_\eps\Big( \frac{y}{s(X_i)}\Big)\Big)
\\
&=:& \frac{1}{nh_n^{m+1}}\sum_{i=1}^n g^{(n,m,+)}_{x,y}(X_i,U_i,c_n).
\eean
The expectation of each summand $g^{(n,m,+)}_{x,y}(X_i,U_i,c_n)$ in the above sum is of the order $O(h_nc_n)$. Moreover, the class of functions
\[
\Big\{ (u,v) \mapsto g^{(n,m,+)}_{x,y}(u,v,c_n) \Big| x\in [4h_n,1-4h_n], y \in \mathcal{Y}\Big\}
\]
is with probability tending to one contained in a class that satisfies the assumptions of part 2 of Lemma \ref{lem:base} with $\delta_n = h_n$, this follows from a combination of assumption \ref{as:k2} with parts 1,2,4,6 of Lemma \ref{lem:entrnum} where part 6 is applied with the class of functions $\G := \{v \mapsto F_\eps(y/s(v))+z|y\in \Y,z \in [0,1]\}$. This yields the bound
\[
\frac{1}{nh_n^{m+1}}\sum_{i=1}^n g^{(n,m,+)}_{x,y}(X_i,U_i,c_n) = o\Big(\frac{c_n h_n}{h_n^{m+1}}\Big) + O_P\Big(\frac{\log n}{nh_n^{2m+1}} \Big)^{1/2}
\] 
uniformly in $x\in [4h_n,1-4h_n], y \in \mathcal{Y}$. Since $c_n/r_n$ can tend to infinity arbitrarily slowly, the above result implies
\[
\frac{1}{nh_n^{m+1}}\sum_{i=1}^n g^{(n,m,+)}_{x,y}(X_i,U_i,c_n) = O_P\Big(\frac{\log n}{nh_n^{2m+1}} \Big)^{1/2}.
\]
Summarizing, we have obtained the bound
$
A_{n,k,m}^{+}(x,y) 
\leq  O_P\Big(\frac{\log n}{nh_n^{2m+1}} \Big)^{1/2},
$
and a corresponding lower bound can be obtained by similar arguments. Analogous reasoning yields a bound for $A_{n,k,m}^{-}(x,y)$ and altogether this implies (\ref{eq:help6}).

Thus we have established (\ref{eq:help2})-(\ref{eq:help6}) and the proof of the Lemma is complete.  \hfill $\Box$

\begin{lemma}\label{lem:propweight}
Under assumptions \ref{as:k1} and \ref{as:fx} if additionally $(nh_n)^{-1} = o(h_n\sqrt{\log n})$ we have the decomposition (holding uniformly in $x \in [h_n,1-h_n]$)
\[
nh_n (\XX^t\WW\XX)^{-1} = \frac{1}{f_X(x)}\mathcal{H}^{-1}\mathcal{M}(K)^{-1}\mathcal{H}^{-1} + \mathcal{H}^{-1}1_{(p+1)\times (p+1)}O_P(h)\mathcal{H}^{-1}
\]
where 
$\mathcal{H} = diag(1,h_n,...,h_n^p)$, and $1_{(p+1)\times (p+1)}$ is a matrix with 1 in every entry. 

\end{lemma}
\textbf{Proof} The elements of the matrix $\XX^t\WW\XX$ are of the form 
\bean
\frac{1}{nh_n}(\XX^t\WW\XX)_{k,l} = \frac{1}{nh_n}\sum_i K_{h_n,0}(x-X_i)(x-X_i)^{m} = \frac{h_n^{m}}{nh_n^d}\sum_i K_{h_n,m}(x-X_i)
\eean
where $m = k+l-2$. In particular, continuous differentiability of $f_X$ together with an application of Lemma \ref{lem:base} and Lemma \ref{lem:entrnum} implies that
\[
\frac{1}{nh_n}\sum_i K_{h_n,k}(x-X_i) = \mu_k f_X(x) + O_P(\Big(\frac{\log n}{nh_n} \Big)^{1/2} + h_n) 
\]
uniformly in $x$. Thus we obtain a representation of the form 
\[
\frac{1}{nh_n} \XX^t\WW\XX = \mathcal{H}\Big(\mathcal{M}(K) f_X(x) + 1_{N\times N}O_P(h_n)\Big)\mathcal{H} 
\]
where $M_0 = \mathcal{M}(K)$ is invertible and $\mathcal{H}$ is a diagonal matrix with entries $1,h_n,...,h_n^{p}$. Thus for $h_n$ sufficiently small an application of the Neumann series yields the assertion with probability tending to one.
\hfill $\Box$

\newpage

\newpage
\subsection{Additional technical results} \label{sec:tec}

\begin{lemma} \label{lem:propfepsdach}
Let $n\alpha_n^4 = o(1)$ and assume that the conditions of (i), (i)', (ii), (ii)' of Lemma \ref{lem:propfdach} hold. Then for any bounded $\Y \subset \R$ and any $\delta_n \to 0$ we have 
\bea \label{eq:c71}
\sup_{a,b\in\Y,|a-b|\leq \delta_n} \Big|\tilde F_\eps(a) - \tilde F_\eps(b) - \Big(\bar F_\eps(a) - \bar F_\eps(b)\Big) \Big| &=& o_P(1/\sqrt{n}),
\\ \label{eq:c72}
\sup_{y \in \Y} \Big|\tilde F_\eps(y) - F_\eps(y)\Big| &=& O_P\Big(\Big(\frac{\log n}{n h_n}\Big)^{1/2}\Big),
\eea
where 
\[
\bar F_\eps(a) := \frac{\sum_k I_{[2h_n,1-2h_n]}(X_k)F_Y(\hat q_{\tau,L}(X_k) + a\hat s_L(X_k)|X_k)}{\sum_l I_{[2h_n,1-2h_n]}(X_l)}.
\]
\end{lemma}
\textbf{Proof of Lemma \ref{lem:propfepsdach}} 
Recalling the definition of $\tilde F_\eps,$ it is easy to see that $\tilde F_\eps(y) = \frac{1}{\alpha_n}\Big(\hat F_\eps(\cdot) * \phi( \cdot/\alpha_n)\Big)(y)$ where
\[
\hat F_\eps(y) := \frac{\sum_k I_{[2h_n,1-2h_n]}(X_k) I\{Y_k - \hat q(X_k) \leq y \hat s(X_k) \}}{\sum_l I_{[2h_n,1-2h_n]}(X_l)}.
\]
Standard calculations show that 
\[
\frac{1}{\alpha_n}\Big(\bar F_\eps(\cdot) * \phi( \cdot/\alpha_n)\Big)(y) = \bar F_\eps(y) + o_P(1/\sqrt n)
\]
uniformly in $y \in\Y$. Thus it suffices to establish that, for any bounded $\tilde \Y$ 
\bea \label{eq:c71'}
\sup_{a,b \in \tilde \Y, \ |a-b|\leq \delta_n} \Big|\hat F_\eps(a) - \hat F_\eps(b) - \Big(\bar F_\eps(a) - \bar F_\eps(b)\Big) \Big| &=& o_P(1/\sqrt{n})
\\ \label{eq:c72'}
\sup_{y \in \tilde\Y} \Big|\hat F_\eps(y) - F_\eps(y)\Big| &=& O_P\Big(\Big(\frac{\log n}{n h_n}\Big)^{1/2}\Big).
\eea
To simplify the notation, write $\Y$ for $\tilde \Y$.\\

\textbf{Proof of (\ref{eq:c71'})}
Since $\frac{1}{n}\sum_l I_{[2h_n,1-2h_n]}(X_l) = 1 + o_P(1)$, we only need to consider the enumerator. Since $\Y$ is bounded we have, with probability tending to one, uniformly in $y\in\Y$
\bean
&&\Big| I\{Y_k - \hat q_\tau(X_k) \leq y \hat s(X_k) \} -  I\{Y_k - \hat q_{\tau,L}(X_k) \leq y \hat s_L(X_k) \}\Big|
\\
&\leq& I\{ Y_k - \hat q_{\tau,L}(X_k) - y\hat s_L(X_k) \leq \gamma_n  \} - I\{ Y_k - \hat q_{\tau,L}(X_k) - y\hat s_L(X_k) \leq -\gamma_n  \}
\eean
for some $\gamma_n = o(1/\sqrt n)$. Moreover an application of parts 1 and 6 of Lemma \ref{lem:entrnum} combined with Theorem 2.7.1 in van der Vaart, Wellner (1996) shows that the functions
\[
(u,v) \mapsto I\{ v - \hat q_{\tau,L}(u) - y\hat s_L(u) \leq \gamma_n  \} - I\{ v - \hat q_{\tau,L}(u) - y\hat s_L(u) \leq -\gamma_n  \}
\]
are, with probability tending to one, contained in a class of functions satisfying the assumptions of the first part of Lemma \ref{lem:base} with the additional property that each element has expectation of order $o(1/\sqrt n)$. Combined with parts 1 and 4 of Lemma \ref{lem:entrnum}, this implies 
\[
\sup_{y \in \Y} \Big|\sum_k I_{[2h_n,1-2h_n]}(X_k)\Big( I\{Y_k - \hat q_\tau(X_k) \leq y \hat s(X_k) \} - I\{Y_k - \hat q_{\tau,L}(X_k) \leq y \hat s_L(X_k) \} \Big) \Big| = o_P(1/\sqrt n),
\]
and thus it remains to consider
\bean
&&\sup_{a,b\in\Y, |a-b| \leq \delta_n}\frac{1}{n}\sum_k I_{[2h_n,1-2h_n]}(X_i)\Big(I\{Y_k \leq \hat q_{\tau,L}(X_k) + a\hat s_L(X_k)\} - I\{Y_k \leq \hat q_{\tau,L}(X_k) + b\hat s_L(X_k)\}
\\
&&\qquad\qquad\qquad\qquad\qquad\qquad    
- F_Y(\hat q_{\tau,L}(X_k) + a\hat s_L(X_k)|X_k) + F_Y(\hat q_{\tau,L}(X_k) + b\hat s_L(X_k)|X_k) \Big)
\eean
By arguments similar to those given above, it is easily seen that this quantity is of order $o_P(1/\sqrt n)$ if one notes that the smoothness assumptions on $F_Y$ imply that with $\hat q_{\tau,L},\hat s_L \in C_C^{1+\delta}$ with probability tending to one the same holds for the function $u\mapsto F_Y(\hat q_{\tau,L}(u) + y\hat s_L(u)|u)$ uniformly in $y\in\Y$. This completes the proof of (\ref{eq:c71'}).\\

\textbf{Proof of (\ref{eq:c72'})}
Write
\[
\hat F_\eps(y) - F_\eps(y) = \frac{n^{-1}\sum_k I_{[2h_n,1-2h_n]}(X_k) \Big(I\{Y_k - \hat q_\tau(X_k) \leq y \hat s(X_k) \} - F_\eps(y)\Big)}{n^{-1}\sum_l I_{[2h_n,1-2h_n]}(X_l)}.
\]
Since $n^{-1}\sum_l I_{[2h_n,1-2h_n]}(X_l) = 1 + o_P(1)$, it suffices to consider the enumerator. Observe that 
\[
I\{Y_k - \hat q_\tau(X_k) \leq y \hat s(X_k) \} = I\Big\{\eps_k \leq y\frac{\hat s(X_k)}{s(X_k)} + \frac{\hat q_{\tau}(X_k) - q_\tau(X_k)}{s(X_k)} \Big\}
\]
and thus, for any $c_n/r_n \to \infty$ we have with probability tending to one, uniforly over $y\in \Y$ 
\bean
\Big|I\{Y_k - \hat q_\tau(X_k) \leq y \hat s(X_k) \} - I\{\eps_k < y\}\Big| \leq I\{|\eps_k -y| \leq c_n\}.
\eean
Thus standard 
\bean
&&\sup_{y \in \Y} \Big|n^{-1}\sum_k I_{[2h_n,1-2h_n]}(X_k) \Big(I\{Y_k - \hat q_\tau(X_k) \leq y \hat s(X_k) \} - I\{\eps_k < y\}\Big)\Big|
\\
&\leq& \sup_{y \in \Y} n^{-1}\sum_k I_{[2h_n,1-2h_n]}(X_k)I\{|\eps_k -y| \leq c_n\} = O_P(c_n),
\eean
where the last equality follows by standard empirical process arguments. This shows that, uniformly in $y \in \Y$,
\[
\hat F_\eps(y) - F_\eps(y) = \frac{n^{-1}\sum_k I_{[2h_n,1-2h_n]}(X_k) \Big(I\{\eps_k \leq y \} - F_\eps(y)\Big)}{n^{-1}\sum_l I_{[2h_n,1-2h_n]}(X_l)} + O_P(c_n) = O_P(c_n).
\]

Since $c_n$ was arbitrary, this completes the proof of (\ref{eq:c72'}) and hence also of the Lemma. \hfill $\Box$

\newpage
\begin{lemma} \label{lem:genlin} 
Assume that $\kappa$ is a symmetric, uniformly bounded density with support $[-1,1]$ and let $b_n = o(1)$.\\ 
(a) If the function $F: [0,1]\rightarrow \R$ is strictly increasing and $F^{-1}$ is $k$ times continuously differentiable in a neighborhood of the point $\tau$, we have for $b_n$ small enough
\[
H_{id,\kappa, \tau,b_n}(F) = F^{-1}(\tau) + \sum_{i=1}^k \frac{b_n^i}{i!}(F^{-1})^{(i)}(\tau)\mu_{i+1}(\kappa) + R_n(\tau) 
\] 
with $|R_n(\tau)| \leq C_k(\kappa) b_n^k \sup_{|s-\tau|\leq b_n}|(F^{-1})^{(k)}(\tau)-(F^{-1})^{(k)}(s)|$, $\mu_i(\kappa) := \int u^i\kappa(u)du$ and a constant $C_k$ depending only on $k$ and $\kappa$. In particular, if we assume that $F: \R\rightarrow [0,1]$ is strictly increasing and $F^{-1}$ is two times continuously differentiable in a neighborhood of $\tau$ and $G:\R\to(0,1)$ is two times continuously differentiable in a neighborhood of $F^{-1}(\tau)$ with $G'(F^{-1}(\tau))>0$ we have 
\[
|F^{-1}(\tau) - Q_{G,\kappa,\tau,b_n}(F)| \leq Cb_n^2\sup_{|s-G\circ F^{-1}(\tau)|\leq R_{n,1}}|(G^{-1})'(s)|\sup_{|s-\tau|\leq b_n }|(G\circ F^{-1})''(s)| =: R_{n,2}
\]
for some constant $C$ that depends only on $\kappa$ where $R_{n,1} := Cb_n^2\sup_{|s-\tau|\leq b_n }|(G\circ F^{-1})''(s)|$.\\
\\
(b) Assume that $\kappa$ is additionally differentiable with Lipschitz-continuous derivative and that the functions $G,G^{-1}$ have derivatives that are uniformly bounded on any compact subset of $\R$ [the bound is allowed to depend on the interval].
Then for any increasing function $F$ with uniformly bounded first derivative we have $|H(F_1)-H(F_2)|\leq R_{n,3} + R_{n,4}$ and
\[
|Q_{G,\kappa,\tau,b_n}(F_1) - Q_{G,\kappa,\tau,b_n}(F_2)| \leq \sup_{u\in \mathcal{U}(H(F_1),H(F_2))}|(G^{-1})'(u)|(R_{n,3} + R_{n,4})
\]
where $C$ is a constant that depends only on $\kappa$, $\mathcal{U}(a,b) := [a \wedge b, a \vee b]$, and
\[
R_{n,3} := \frac{Cc_n}{b_n}\|F_1-F_2\|_\infty \sup_{|v-\tau|\leq c_n }|(G\circ F^{-1})'(v)|, \quad R_{n,4} := R_{n,3} \frac{\|F_1-F\|_\infty + \|F_1-F_2\|_\infty}{b_n}
\]
with $c_n := b_n + 2\|F_1-F_2\|_\infty + \|F_1-F\|_\infty$.\\
\\
(c) If additionally to the assumptions made in (b), the function $F_1$ is two times continuously differentiable in a neighborhood of $F^{-1}(\tau)$ with $F_1'(F_1^{-1}(\tau )) > 0$  and $G$ is two times continuously differentiable in a neighborhood of $F_1^{-1}(\tau)$ with $G'(F^{-1}(\tau))>0$, we have
\bean
Q_{G,\kappa,\tau,b_n}(F_1) - Q_{G,\kappa,\tau,b_n}(F_2)
&=& -\frac{1}{F_1'(F_1^{-1}(\tau ))}\int_{-1}^1 \kappa(v)\Big(F_2(F_1^{-1}(\tau+vb_n)) - F_1(F_1^{-1}(\tau+vb_n)) \Big)dv
\\
&&+ R_n,
\eean
where
\bean
|R_n| &\leq& R_{n,5} + R_{n,6} + \frac{Cb_n \sup_{|s-\tau|\leq b_n}(G\circ F^{-1})''(s)\|F_1-F_2\|_\infty + R_{n,4}}{G'(F_1^{-1}(\tau))}
\eean
with a constant $C$ depending only on $\kappa$ and
\bean
R_{n,5} &:=& \frac{1}{2}\sup_{u\in \mathcal{U}(H(F_1),H(F_2))}|(G^{-1})''(u)|(H(F_1)-H(F_2))^2
\\
R_{n,6} &:=& \sup_{u\in \mathcal{U}(H(F_1),G(F_1^{-1})(\tau))}|(G^{-1})''(u)|\cdot |H(F_1) - G(F_1^{-1})(\tau)|\cdot |H(F_1)-H(F_2)|.
\eean
\end{lemma}

\noindent\textbf{Proof}
See Volgushev et al.\ (2013).

\begin{lemma}[Basic Lemma]\label{lem:base}\
\begin{enumerate}
\item Assume that the classes of functions $\ef_n$ consist of uniformly bounded functions (with the bound, say $D$, not depending on $n$) with $N_{[ ]}(\ef_n,\eps,L^2(P))\leq C\exp(-c\eps^{-a})$ for every $\eps\leq\delta_n$ for some $a<2$ and constants $C,c$ not depending on $n$. Then we have
\[
\sqrt{n} \sup_{f\in \ef_n, \|f\|_{P,2} \leq \delta_n} \Big(\int f dP_n - \int f dP\Big) = o_P^*(1)
\]
where the $^*$ denotes outer probability, see van der Vaart and Wellner (1996) for a more detailed discussion.
\item If under the assumptions of part one we have $N_{[ ]}(\ef_n,\eps,L^2(P))\leq C\eps^{-a}$ for every $\eps\leq \delta_n$, some $a>0$ and $C$ not depending on $n$, it holds that for any $\delta_n \sim n^{-b}$ with $b<1/2$
\[
\sqrt{n} \sup_{f\in \ef_n, \|f\|_{P,2} \leq \delta_n} \Big(\int f dP_n - \int f dP\Big) = O_P^*\Big( \delta_n |\log \delta_n| \Big)
\]  
\end{enumerate}
\end{lemma}
\textbf{Proof}
See Volgushev et al.\ (2013).

\begin{lemma}\label{lem:entrnum} \ 
\begin{enumerate}
\item \label{entrnum0} Define $\ef+\G := \{f+g|f\in\ef,g\in\G\}, \ef\G := \{fg|f\in\ef,g\in\G\}$. Then 
\[
N_{[]}(\ef+\G, \eps, \rho)\leq N_{[]}(\ef, \eps/2, \rho)N_{[]}(\G, \eps/2, \rho)
\]

If additionally the classes $\ef,\G$ are uniformly bounded by the constant $C$, we have
\[
N_{[]}(\ef\G, \eps, \|.\|) \leq N_{[]}^2(\ef, \eps/4C, \|.\|)N_{[]}^2(\G, \eps/4C, \|.\|)
\] 
for any seminorm $\|.\|$ with the additional property that $|f_1|\leq|f_2|$ implies $\|f_1\|\leq \|f_2\|$.
\item \label{entrnum2} Assume that the Kernel $K$ has compact support $[-1,1]$, that $K_{1,k}^{(m)}$ is uniformly bounded and Lipschitz-continuous, and that $f_X$ is uniformly bounded. Then the $L^2(P_X)$ bracketing numbers $N_{[]}(\mathcal F_n, \eps, L^2(P_X))$ of the set 
\[
\mathcal F_n := \Big\{ u \mapsto K_{h_n,k}^{(m)}(x-u)\Big| x \in [h_n,1-h_n] \Big\}
\]
are bounded by $C\eps^{-3}$ for some constant $C$ independent of $n$.
\item \label{entrnum3} Assume that the Kernel $K$ has compact support $[-1,1]$, that $K$ is uniformly bounded and Lipschitz continuous, and that $f_X$ is uniformly bounded away from zero on $[0,1]$ and Lipschitz-continuous. Then for the set of function
\[
\mathcal F_n := \Big\{ u \mapsto \frac{1}{h_n}\Big(\frac{1}{f_X(x)} - \frac{1}{f_X(u)}\Big)K_{h_n,k}(x-u)\Big| x \in [h_n,1-h_n] \Big\}
\]
we have $N_{[]}(\mathcal F_n, \eps, L^2(P)) \leq C\eps^{-5}$ for some constant $C$ independent of $n$.
\item \label{entrnum4} For any measure $P$ on the unit interval with uniformly bounded density $f$, the class of functions 
\[
\ef := \Big\{ u \mapsto I\{u\leq s\} \Big| s \in [0,1] \Big\} \cup \Big\{ u \mapsto I\{u < s\} \Big| s \in [0,1] \Big\}
\] 
can be covered by $C\eps^{-(2)}$ brackets of $L^2(P)$ length $\eps$. 

\item \label{entrnum5} Consider the class of distribution functions $\ef:=\Big\{ u  \mapsto F(y|u) \Big| y \in \R \Big\}$ with densities $f(y|u)$ and assume that $\sup_{u,y}|y|^\alpha(F(y|u)\wedge(1-F(y|u))\leq D$ for some $\alpha>0$ and additionally $\sup_{u,y} f(y|u) \leq D$. Then we have $N_{[]}(\mathcal F, \eps, \|\ \|_\infty) \leq C\eps^{-\frac{\alpha+1}{\alpha}}$ for some constant $C$ independent of $\alpha$.
\item \label{entrnum7} For any measure $P$ on $\R\times\R^k$ with uniformly bounded conditional density $f_{V|U}$ the class of functions
\[
\G := \Big\{(u,v) \mapsto I\{v \leq f(u)\} \Big| f \in \ef \Big\}
\]
satisfies $N_{[]}(\G, \eps, \|.\|_{P,2}) \leq N_{[]}(\ef, C\eps^2, \|.\|_\infty)$ for some constant $C$ independent of $\eps$. 
\end{enumerate}
\end{lemma}
\textbf{Proof}\\
\textbf{Part \ref{entrnum0}} 
The first assertion is obvious from the definition of bracketing numbers. For the second assertion, note that $\ef\G = (\ef+C)(\G+C) - C\ef - C\G + C^2.$ Moreover, all elements of the classes $\ef+C,\G+C$ are by construction non-negative and thus it also is possible to cover them with brackets consisting of non-negative functions and amounts equal to the brackets of $\ef,\G$, respectively. Finally, observe that if $0 \leq f_l \leq f \leq f_u$ and $0 \leq g_l \leq g \leq g_u$, we also have $f_lg_l \leq fg \leq f_ug_u$. Moreover $\|f_lg_l - f_ug_u\| \leq C\|f_u-f_l\| + C\|g_u-g_l\|$. Thus the class $(\ef+C)(\G+C)$ can be covered by at most $\leq N_{[]}(\ef, \eps, \|.\|)N_{[]}(\G, \eps, \|.\|)$ brackets of length $2C\eps$. Finding brackets for the classes $C\ef, C\G$ is trivial, and applying the first assertion of the Lemma completes the proof. 
\\
\textbf{Part \ref{entrnum2}+\ref{entrnum3}} 
Without loss of generality, assume that $h=h_n<1$. Note that the class of functions $\mathcal{F}_n$ from part \ref{entrnum2} can be represented as $\mathcal{F}_n = \{u \mapsto g_x(u) | x \in [h_n,1-h_n]\}$ where the functions $g_x$ satisfy $\sup_{x \in [h_n,1-h_n]} \|g_x\|_\infty \leq C$, $\sup_{u \in \R } |g_x(u) - g_y(u)| \leq \tilde C|x-y| h_n^{-1}$ for some constants $C,\tilde C$ independent of $n,x,y$. To see the latter inequality, observe that by assumption $u \mapsto K^{(m)}_{1,k}(u)$ is uniformly bounded and Lipschitz continuous. Additionally, the support of the functions $g_x$ is contained in $[x-h_n,x+h_n]$.

Similarly, $\mathcal{F}_n$ from part \ref{entrnum3} can be represented as $\mathcal{F}_n = \{u \mapsto g_x(u) | x \in [h_n,1-h_n]\}$ where the functions $g_x$ satisfy $\sup_{x \in [h_n,1-h_n]} \|g_x\|_\infty \leq C$, $\sup_{u \in \R } |g_x(u) - g_y(u)| \leq \tilde C|x-y| h_n^{-2}$ for some constants $C,\tilde C$ independent of $n,x,y$ (and possibly different from those for part \ref{entrnum2}), and the support of the functions $g_x$ is contained in $[x-h_n,x+h_n]$.

Thus it suffices to establish that for any class of functions $\mathcal F$ of the form $\mathcal{F} = \{u \mapsto g_x(u) | x \in [h,1-h]\}$ with $0 \leq h \leq 1/2$ with elements $g_x$ that have support contained in $[x-h,x+h]$ and satisfy $\sup_{x \in [h,1-h]} \|g_x\|_\infty \leq C$, $\sup_{u \in \R } |g_x(u) - g_y(u)| \leq \tilde C|x-y| h^{-L}$ for some constants $C,\tilde C$ independent of $h,x,y$ we have we have $N_{[\ ]}(\mathcal F, \eps, L^2(P_X)) \leq c\eps^{-(2L+1)}$ for some $c$ that does not depend on $h$. 

To prove this statement, consider two cases.
\begin{itemize}
\item[1] $\eps> 4h^{1/2}$\\ 
Divide $[0,1]$ into $N := 2/\eps^2$ subintervals of length $2\alpha := \eps^2$ with centers $r\alpha$ for $r=1,...,N$ and call the intervals $I_1,...,I_N$. Note that two adjunct intervals overlap by $\alpha > 2h$. This construction ensures that every set of the form $[x-h,x+h]$ with $x\in[h,1-h]$ is completely contained in at least one of the intervals defined above. Then a collection of $N$ brackets of $L^2$-length $D\eps$ for some $D>0$ independent of $h$ is given by $(-CI\{u \in I_j\}, CI\{u \in I_j\})$.
\item[2] $\eps \leq 4h^{1/2}$\\
Consider the points $t_i := i/(N+1), i=1,...,N$ with $N := 4^{2L+2}\tilde C/\eps^{2L+1}$. By construction, to every $x\in[h,1-h]$ there exists $i(x)$ with $|t_{i(x)}-x|\leq \eps^{2L+1}/(4^{2L+2}\tilde C)$. This implies
\[
\sup_u |g_x(u) - g_{t_{i(x)}}(u)| \leq \tilde C  \eps^{2L+1} h^{-L}/(4^{2L+2}\tilde C) < \eps/2
\]
Then $N$ $\|\cdot\|_\infty-$brackets of length $\eps$ covering $\ef$ are given by $(g_{t_i}(\cdot)-\eps/2,g_{t_i}(\cdot)+\eps/2)$, $i=1,...,N$. From those one can easily construct $L^2(P_X)$-brackets. 
\end{itemize}
\textbf{Part \ref{entrnum4}} Follows by standard arguments.
\\
\textbf{Part \ref{entrnum5}} For any $\eps>0$, set $y_\eps := \eps^{-1/\alpha}D^{1/\alpha}$ and define $t_i:= -y_\eps + i\eps/D$ for $i=1,...,N$ with $N$ such that $1 + y_\eps \geq t_N \geq y_\eps$. Note that $N \leq C\eps^{-\frac{\alpha+1}{\alpha}}$ for some fixed, finite constant $C$ which can depend on $D$ but not on $\eps$. The collection of brackets $(f\equiv 0,f\equiv\eps),(f\equiv 1-\eps,f\equiv 1),(F(y_{t_i}|.)-\eps/2,F(y_{t_i}|.)+\eps/2)$ with $i=1,...,N$ covers the class $\ef$. To see that, let $f \in \mathcal{F}$. Then there exists $y \in \R$ such that $f(\cdot) = F(y|\cdot)$. If $y < - y_\eps$ we have 
\[
0 \leq F(y|u) \leq \sup_u F(-y_\eps|u) \leq y_\eps^{-\alpha}\sup_u y_\eps^{\alpha}F(-y_\eps|u) \leq D (\eps^{-1/\alpha}D^{1/\alpha})^{-\alpha} = \eps.
\]
Similarly, $y > y_\eps$ implies $1- \eps \leq F(y|u) \leq 1$
Finally, if $-y_\eps \leq y \leq y_\eps$, there exists $i \in \{1,...,N\}$ such that $|y-t_i| \leq \eps/(2D)$. In that case 
\[
F(t_i|u)-\eps/2 \leq |F(t_i|u) - F(y|u)| + F(y|u) - \eps/2 \leq F(y|u) \leq F(t_i|u) + \eps/2
\]
since $|F(t_i|u) - F(y|u)| \leq D |t_i-y| \leq \eps/2$ by the assumption $\sup_{u,y} f(y|u) \leq D$.
\\ 
\textbf{Part \ref{entrnum7}} Follows from $|I\{v \leq g_1(u)\} - I\{v \leq g_2(u)\}| \leq I\{|v-g_1(u)|\leq 2 \|g_1-g_2\|_\infty\}$.

\hfill $\Box$

\newpage

\subsection{Main results for proofs}
\label{app-tech}

Define $\hat \eps_{i,L}$ as the estimated residuals based on linearized versions $\hat q_{\tau,L}, \hat s_L$ [see Appendix \ref{app-quant} for their definition], i.e. $\hat\eps_{i,L} := (Y_i - \hat q_{\tau,L}(X_i))/\hat s_ L(X_i)$, and $\hat\eps_{i,L}^*$ as the corresponding quantities in the bootstrap setting, that is
\[
\hat\eps_{i,L}^*=\frac{\hat s_L(X_i)\eps_i^*+\hat q_{\tau,L}(X_i)-\hat q_{\tau,L}^*(X_i)}{\hat s_L^*(X_i)}
\]
The following Lemma demonstrates, that the sequential empirical process based on the residuals $\hat \eps_i = (Y_i - \hat q_\tau(X_i))/\hat s(X_i)$ computed from the initial estimators $\hat q_\tau, \hat s$ and the sequential empirical process of residuals based on $\eps_{i,L}$ have the same first order expansion.  

\begin{lemma}\label{lem:proclin}
Assume that \ref{as:k1}-\ref{as:Gs}, \ref{as:fx}-\ref{as:fbound}, \ref{as:bw} hold. Then
\[
\sup_{t \in [2h_n,1-2h_n],y\in \R} \Big|\frac{1}{\sqrt n}\sum_i I\{2h_n \leq X_i\leq t\}(I\{\hat \eps_{i} \leq y\} - I\{\hat\eps_{i,L} \leq y\})\Big| = o_P(1).
\]
If additionally \ref{as:b0}-\ref{as:b2} hold we also have
\[
\sup_{t \in [4h_n,1-4h_n],y\in \R} \Big|\frac{1}{\sqrt n}\sum_i I\{4h_n\leq X_i\leq t\}(I\{\hat \eps_{i}^* \leq y\} - I\{\hat\eps_{i,L}^* \leq y\})\Big| = o_P(1).
\]
\end{lemma}
\textbf{Proof} We only proof the second assertion since the first one follows by similar but easier arguments. Start by observing that under the assumptions of the Lemma there exists a set $D_n$ whose probability tends to one such that on $D_n$ we have
\bean
(i) & \sup_{x \in [4h_n,1-4h_n]} \max\Big(|\hat q_\tau(x) - \hat q_{\tau,L}(x)|,|\hat q_\tau^*(x) - \hat q_{\tau,L}^*(x)|,|\hat s(x) - \hat s_L(x)|,|\hat s^*(x) - \hat s_L^*(x)|\Big) \leq \gamma_n&
\\
(ii) &\inf_{x \in [4h_n,1-4h_n]} \min(\hat s_L(x), \hat s_L^*(x)) \geq c > 0&
\\
(iii) & \sup_{y \in \R} |y \tilde f_\eps(y)| \leq C&
\eean 
for some deterministic sequence $\gamma_n = o(1/\sqrt{n})$ and finite constants $C,c>0$. Here (i) and (ii) follow from Lemma \ref{lem:quantlin} and Lemma \ref{lem:unifratesquant} together Assumption \ref{as:s}, while (iii) is a consequence of (\ref{le2n09}) in the main body of the paper.

A standard Taylor expansion shows that on $D_n$ 
\bean
\Big| I\{\hat \eps_{i}^* \leq y\} - I\{\hat \eps_{i,L}^* \leq y\} \Big|
&\leq& I\Big\{ \Big| U_i - \tilde F_\eps\Big(y\frac{\hat s_L^*(X_i)}{\hat s_L(X_i)} + \frac{\hat q_{\tau,L}^*(X_i) - \hat q_{\tau,L}(X_i)}{\hat s_L(X_i)}\Big)\Big| \leq C\gamma_n \Big\}
\\
&=:& g_{n,y,C\gamma_n}(U_i,X_i),
\eean
this follows from the representations
\bean
I\{\hat \eps_{i}^* \leq y\} &=& I\Big\{  U_i \leq \tilde F_\eps\Big(y\frac{\hat s^*(X_i)}{\hat s(X_i)} + \frac{\hat q_{\tau}^*(X_i) - \hat q_{\tau}(X_i)}{\hat s(X_i)}\Big)\Big\},
\\
I\{\hat \eps_{i,L}^* \leq y\} &=& I\Big\{  U_i \leq \tilde F_\eps\Big(y\frac{\hat s_L^*(X_i)}{\hat s_L(X_i)} + \frac{\hat q_{\tau,L}^*(X_i) - \hat q_{\tau,L}(X_i)}{\hat s_L(X_i)}\Big)\Big\},
\eean
a Taylor expansion of $\tilde F_\eps$ and (i)-(iii). In the same manner as the proof of Proposition 3 in Neumeyer (2009a) it follows from assumptions \ref{as:b1} and \ref{as:b2} that, with probability tending to one
\beq \label{eq:class1}
\quad\quad \G_{n} := \Big\{ (u,v) \mapsto I\Big\{ u \leq z + \tilde F_\eps\Big(y\frac{\hat s_L^*(v)}{\hat s_L(v)} + \frac{\hat q_{\tau,L}^*(v) - \hat q_{\tau,L}(v)}{\hat s_L(v)}\Big)\Big\}  \Big|\  y\in \R, z \in [-2,2] \Big\} 
\eeq
is contained in the class
\begin{eqnarray*}
 \tilde \G_n =\Big\{(u,v) \mapsto I\Big\{ u \leq z + F\Big(y\frac{a_3(v)}{a_1(v)} + \frac{a_2(v)}{a_1(v)}\Big)\Big\}  \Big|\  
 F\in\mathcal{D}, a_1,a_3\in \tilde C_C^{1+\delta}([4h_n,1-4h_n]),\\
 a_2\in C_C^{1+\delta}([4h_n,1-4h_n]), y\in \R, z \in [-2,2]
 \Big\},
 \end{eqnarray*}
where $\mathcal{D}$ is defined in (\ref{cal-D}).
Now, denoting by $P$ the product measure of the uniform random variable $U_1$ and the covariate $X_1$, 
\beq \label{eq:class2} 
\log N_{[\ ]}(\eps,\tilde \G, L^2(P)) \leq C\eps^{-2\alpha}
\eeq
for some $\alpha < 1$ , this can be shown by similar arguments as in the proof of Proposition 3 in Neumeyer (2009a). Next, since $I\{|U_1-a|\leq b\} = I\{U_1 \leq a + b\} - I\{U_1 \leq a-b\}$ a.s., we find that, with probability tending to one
\bean
\mathcal{F}_n &:=&\Big\{(u,v) \mapsto I\{s\leq v \leq t\} g_{n,y,C\gamma_n}(v,u) \Big|s,t \in [4h_n,1-4h_n], y \in \R \Big\} 
\\
&\subseteq& 
\Big\{(u,v) \mapsto I\{s\leq v \leq t\} (g_1(v,u) - g_2(v,u)) \Big| s,t \in [4h_n,1-4h_n], g_1,g_2 \in \tilde \G_n \Big\} =: \G_{n,1}.
\eean
Combining parts (1) and (4) of Lemma \ref{lem:entrnum} thus yields that $\log N_{[\ ]}(\eps,\mathcal{F}_n,L^2(P)) \leq \tilde C\eps^{-2\alpha}$ for some constant $\tilde C$. Moreover, standard arguments (employing Taylor expansions and the bounds in (\ref{le2n09}) from the main body of the paper) show that 
$
\sup_{g \in \mathcal F_n} \int gdP = o(1/\sqrt n)
$
and 
$
\sup_{g \in \mathcal F_n} \int g^2 dP = o(1).
$ 
Here, $P$ denotes the probability distribution of $(X_i,U_i)$ and $g^2=g$ for all $g\in\mathcal{F}_n$. Finally observe that, with probability tending to one, 
\bean
&&\sup_{t \in [4h_n,1-4h_n], y \in \R} \frac{1}{\sqrt n}\sum_i \Big(I\{h_n\leq X_i\leq t\}g_{n,y,C\gamma_n}(U_i,X_i) - \int_{h_n}^t \int  g_{n,y,C\gamma_n}(v,u)f_X(u)dvdu ]\Big)
\\
&\leq& \sqrt{n}\sup_{g\in \mathcal{F}_n} (\int g dP_n - \int g dP),
\eean
and the right-hand side of the inequality is of order $o_P(1)$ by part one of Lemma \ref{lem:base} . Moreover, standard arguments yield
\[
\int_{h_n}^t \int g_{n,y,C\gamma_n}(v,u)f_X(u)dvdu = o_P(1/\sqrt{n}). 
\] 
Summarizing, we have obtained the estimate
\[
\sup_{t \in [4h_n,1-4h_n], y \in \R} \frac{1}{\sqrt n}\sum_i I\{4h_n\leq X_i\leq t\}g_{n,y,C\gamma_n}(U_i,X_i) = o_P(1).  
\]
and thus the proof is complete.
\hfill $\Box$

\newpage

\begin{lemma}\label{lem:B2boot}
Assume that the conditions \ref{as:k1}-\ref{as:Gs}, \ref{as:fx}-\ref{as:fbound}, \ref{as:bw} hold. Then
\[
\int_{h_n}^{t} \frac{\hat q_{\tau,L}(x) - q_\tau(x)}{s(x)}f_X(x)f_\eps(0)dx = - \frac{1}{n}\sum_{i=1}^n (I\{\eps_i \leq 0\} - \tau)I_{[h_n,t]}(X_i) + o_P(1/\sqrt{n})
\]
uniformly in $t\in [h_n,1-h_n]$ and
\bean
&&\int_{2h_n}^t\frac{\hat s_L(x)-s(x)}{\hat s(x)}f_X(x)\,dx\\
&=& -\frac{1}{n}\sum_{i=1}^n \frac{I_{[2h_n,t]}(X_i)}{f_{|\eps|}(1)}\Big(
I\{|\eps_i|\leq 1\} - \frac{1}{2} - \frac{(I\{\eps_i\leq 0\}-\tau)(f_\eps(1) - f_\eps(-1))}{f_\eps(0)}\Big)
+ o_P(\frac{1}{\sqrt{n}})
\eean
uniformly in $t\in [2h_n,1-2h_n]$.\\
If additionally \ref{as:b0}-\ref{as:b2} hold
\[
\int_{3h_n}^{t} \frac{\hat q_\tau^*(x) - \hat q_{\tau,L}(x)}{\hat s_L(x)}f_X(x)dx = - \frac{1}{n}\sum_{i=1}^n \frac{I\{\eps_i^* \leq 0\} - \tau}{f_\eps(0)}I_{[3h_n,t]}(X_i) + o_P(1/\sqrt{n})
\]
uniformly in $t\in [3h_n,1-3h_n]$ and
\bean
&&\int_{4h_n}^t\frac{\hat s^*(x)-\hat s(x)}{\hat s(x)}f_X(x)\,dx\\
&=& -\frac{1}{n}\sum_{i=1}^n \frac{I_{[4h_n,t]}(X_i)}{f_{|\eps|}(1)}\Big(
I\{|\eps_i^*|\leq 1\} - \frac{1}{2} - \frac{(I\{\eps_i^*\leq 0\}-\tau)(f_\eps(1) - f_\eps(-1))}{f_\eps(0)}\Big)
+ o_P(\frac{1}{\sqrt{n}})
\eean
uniformly in $t\in [4h_n,1-4h_n]$.\\
\end{lemma}
\textbf{Proof} We will only prove the representation for $\int_{3h_n}^{t} \frac{\hat q^*(x) - \hat q_{\tau,L}(x)}{\hat s_L(x)}f_X(x)dx$ since all other results can be derived by analogous arguments.\\
Observe the decomposition $\hat q_\tau^*(x) - \hat q_{\tau,L}(x) = \hat q_\tau^*(x) - q_\tau(x) + q_\tau(x) - \hat q_{\tau,L}(x)$. By Lemma \ref{lem:unifratesquant} and Lemma \ref{lem:quantlin} we have 
\[
\hat q_\tau^*(x) - \hat q_{\tau,L}^*(x) = o_P(1/\sqrt n),\ \  \hat q_{\tau,L}^*(x) - q_\tau(x) = O_P(r_n),\ \  \hat s_L(x) - s(x) = O_P(r_n),
\] 
uniformly in $x\in[3h_n,1-3h_n]$. It thus suffices to establish 
\bean
\int_{3h_n}^{t} \frac{\hat q_{\tau,L}^*(x) - q_\tau(x)}{s(x)}f_X(x)dx &=& \int_{3h_n}^{t} \frac{\hat q_{\tau,L}(x) - q_\tau(x)}{s(x)}f_X(x)dx - \frac{1}{n}\sum_{i=1}^n \frac{I\{\eps_i^* \leq 0\} - \tau}{f_\eps(0)}I_{[3h_n,t]}(X_i)
\\
&& + o_P(1/\sqrt{n})
\eean
uniformly in $t\in [3h_n,1-3h_n]$. By definition of $\hat q_{\tau,L}^*$, by part (iii)' of Lemma \ref{lem:propfdach}, and since $f_{e}(0|x) = s(x)f_\eps(0)$ we have
\bean
&&\frac{f_X(x)(\hat q_{\tau,L}^*(x) - q_\tau(x))}{s(x)} 
\\ 
&=& -\frac{f_X(x)u_1^t\mathcal{M}(K)^{-1}}{f_\eps(0)} \int_{-1}^1 \kappa(v)\Big(\tilde T_{n,0,L,S}^*(x,q_{\tau+vb_n}(x)),\dots, \tilde T_{n,p,L,S}^*(x,q_{\tau+vb_n}(x)) \Big)^t dv + o_P(1\sqrt n)
\eean
where
\bean
\tilde T_{n,k,L,S}^*(x,y) &=& \frac{1}{nh_n}\frac{1}{f_X(x)}\sum_{i=1}^n K_{h_n,k}(x-X_i) \Big(\Omega\Big(\frac{Y_i^*-y}{d_n} \Big) - F_Y(y|X_i)\Big).
\eean

The remaining proof is based on the following intermediate results which we will establish later on. First of all, uniformly in $t \in [3h_n, 1-3h_n]$, we have 

\bea \label{eq:c2h0}
&& \int_{3h_n}^t \tilde T_{n,k,L,S}^*(x,q_{\tau+vb_n}(x))f_X(x)dx
\\ \nonumber
&=& \frac{1}{n}\sum_i I_{[3h_n,t-h_n]}(X_i)\int_{-1}^1 K_{1,k}(u) 
\Big(\Omega\Big(\frac{Y_i^*-q_{\tau+vb_n}(X_i+uh_n)}{d_n}\Big)
\\ \nonumber
&& \qquad\qquad\qquad\qquad\qquad\qquad\qquad\qquad
- F_Y(q_{\tau+vb_n}(X_i+uh_n)|X_i)\Big)du
+ o_P(1/\sqrt n).
\eea

Moreover we have uniformly in $u \in [-1,1], t \in [3h_n, 1-3h_n]$ 
\bea \nonumber
\quad\quad &&\frac{1}{n}\sum_i I_{[3h_n,t-h_n]}(X_i) \Omega\Big(\frac{Y_i^*-q_{\tau+vb_n}(X_i+uh_n)}{d_n}\Big)
\\ \label{eq:c2h1}
&=& \frac{1}{n}\sum_i I_{[3h_n,t-h_n]}(X_i) \Big(\Omega\Big(\frac{\eps_i^*\hat s_L(X_i)}{d_n}\Big)
+ vb_n\gamma_n(X_i) + \sum_{j=1}^p \xi_j(X_i,v,n)(uh_n)^j \Big)
\\ \nonumber
&& + \frac{f_\eps(0)}{n}\sum_i I_{[3h_n,t-h_n]}(X_i) \Big(\frac{q_\tau(X_i) - \hat q_{\tau,L}(X_i)}{\hat s_L(X_i)} \Big) + o_P(n^{-1/2})
\\ \label{eq:c2h2}
&=& \frac{1}{n}\sum_i I_{[3h_n,t-h_n]}(X_i)
\Big( vb_n\gamma_n(X_i) + \sum_{j=1}^p \xi_j(X_i,v,n)(uh_n)^j \Big)
\\ \nonumber
&& + \frac{1}{n}\sum_i I_{[3h_n,t]}(X_i)I\{\eps_i^* \leq 0\} + f_\eps(0)\int_{3h_n}^{t} \frac{q_{\tau}(x) - \hat q_{\tau,L}(x)}{s(x)} f_X(x) dx + o_P(n^{-1/2}), 
\eea
where $\xi_j,\gamma_n$ denote some functions that do not depend on $u$. Additionally, a Taylor expansion of $(u,v) \mapsto F_Y(q_{\tau+vb_n}(X_i+uh_n)|X_i)$ shows that
\bea \nonumber
&&\frac{1}{n}\sum_i I_{[3h_n,t-h_n]}(X_i) F_Y(q_{\tau+vb_n}(X_i+uh_n)|X_i) 
\\ \label{eq:c2h3}
&=& \frac{1}{n}\sum_i I_{[3h_n,t-h_n]}(X_i)\Big(\tau + vb_n + \sum_{j=1}^p \zeta_j(X_i,v,n)(uh_n)^j \Big) + o_P(n^{-1/2}),
\eea
where the remainder holds uniformly in $u \in [-1,1], t \in [3h_n, 1-3h_n]$ and the functions $\zeta_j$ are again independent of $u$. Plugging (\ref{eq:c2h2}) and (\ref{eq:c2h3}) into (\ref{eq:c2h0}) we find that
\bean
&&\int_{-1}^1 \kappa(v) \int_{3h_n}^t \tilde T_{n,k,L,S}^*(x,q_{\tau+vb_n}(x))dx dv = \sum_{j=0}^p \mu_{k+j}(K)w_j(t) + o_P(1/\sqrt n)
\eean
where
\bean
w_0(t) &:=& \Big( f_\eps(0)\int_{3h_n}^{t} \frac{q_{\tau}(u) - \hat q_{\tau,L}(u)}{s(u)} f_X(u) du + \frac{1}{n} \sum_{i=1}^n I_{[3h_n,t]}(X_i)(I\{\eps_i^* \leq 0\} - \tau) \Big),
\\
w_j(t) &:=& \frac{h_n^j}{n}\sum_{i=1}^n I_{[3h_n,t-h_n]}(X_i)\int_{-1}^1 \kappa(v)(\xi_j(X_i,v,n) - \zeta_j(X_i,v,n)) dv, \quad j=1,...,p.
\eean
Thus, uniformly in $t \in [3h_n, 1-3h_n]$,
\bean
&& f_X(x)\int_{-1}^1 \kappa(v)\Big(\tilde T_{n,0,L,S}^*(x,q_{\tau+vb_n}(x)),\dots, \tilde T_{n,p,L,S}^*(x,q_{\tau+vb_n}(x)) \Big)^t dv
\\
&=& \mathcal{M}(K)(w_0(t),...,w_p(t))^t + o_P(1/\sqrt n).
\eean
Hence the proof will be complete once we establish (\ref{eq:c2h0})-(\ref{eq:c2h2}).\\
\\

\textbf{Proof of (\ref{eq:c2h0})}

Recalling that $K$ has support $[-1,1]$, we obtain for any $t\in [3h_n,1-3h_n]$ the decomposition
\[
K_{h_n,k}(x-X_i)I_{[3h_n,t]}(x) = K_{h_n,k}(x-X_i)I_{[3h_n,t]}(x)\Big(I_{(t-h_n,t+h_n]}(X_i) + I_{[2h_n,3h_n)}(X_i) + I_{[3h_n,t-h_n]}(X_i)\Big).
\]
We will now show that the contributions corresponding to the summands containing $I_{[2h_n,3h_n)}(X_i)$ and $I_{(t-h_n,t+h_n]}(X_i)$ are negligible. Since both expressions can be treated analogously, we only provide the arguments for $I_{(t-h_n,t+h_n]}(X_i)$. By similar arguments as in the proof of Lemma \ref{lem:propfdach} it is easy to show that 
\bean
&&\sup_{t,x \in [3h_n,1-3h_n], y \in \mathcal{Y}}\Big| \frac{1}{nh_n}\sum_{i=1}^n \frac{K_{h_n,k}(x-X_i)}{f_X(x)}I_{(t-h_n,t+h_n]}(X_i) \Big(\Omega\Big(\frac{Y_i^*-y}{d_n} \Big) - F_Y(y|X_i)\Big)\Big|
\\
&&=: A_n(\Y) = O_P(r_n)
\eean 
for any bounded $\Y\subset\R$. Observe that $K_{h_n,k}$ vanishes outside $[-h_n,h_n]$, and since
\[
I\{|x-X_i|\leq h_n\}I_{[3h_n,t]}(x)I_{(t-h_n,t+h_n]}(X_i) \leq I_{[t-2h_n,t+2h_n]}(x)I_{[t-h_n,t+h_n]}(X_i)
\]
we obtain, for a suitably chosen $\Y$,
\bean
&&\Big|\int_{3h_n}^{t}\frac{1}{nh_n}\sum_{i=1}^n \frac{K_{h_n,k}(x-X_i)}{f_X(x)}I_{[t-h_n,t+h_n]}(X_i) \Big(\Omega\Big(\frac{Y_i^*-q_{\tau+vb_n}(x)}{d_n}\Big) - F_Y(q_{\tau+vb_n}(x)|X_i)\Big)dx\Big|
\\
&\leq& \int_{t-2h_n}^{t+2h_n}A_n(\Y) dx = O_P(h_nr_n) = o_P(1/\sqrt n)
\eean
uniformly in $t \in [3h_n,1-3h_n], v\in [-1,1]$. This completes the proof of (\ref{eq:c2h0}).

\newpage

\textbf{Proof of (\ref{eq:c2h1})}
Throughout this part of the proof, let $\Y \subset \R$ denote a fixed, bounded set containing the interval $[-d_n,d_n]$ for sufficiently large $n$. The following statement will be proved later 
\bea \label{eq:c2h4}
&& \frac{1}{n}\sum_i I_{[3h_n,t-h_n]}(X_i) \Big(I\{Y_i^* \leq q_{\tau+vb_n}(X_i+uh_n)+y\} - I\{\eps_i^* \leq y/\hat s_L(X_i)\} \Big)
\\ \nonumber
&=& \frac{1}{n}\sum_i I_{[3h_n,t-h_n]}(X_i)\Big(\bar F_\eps\Big(\frac{q_{\tau+vb_n}(X_i+uh_n) - \hat q_{\tau,L}(X_i)+y}{\hat s_L(X_i)} \Big) 
\\ \nonumber
&& -\bar F_\eps\Big(\frac{q_{\tau}(X_i) - \hat q_{\tau,L}(X_i)+y}{\hat s_L(X_i)} \Big) + f_\eps\Big(\frac{y}{\hat s_L(X_i)} \Big)\frac{q_{\tau}(X_i) - \hat q_{\tau,L}(X_i)}{\hat s_L(X_i)} \Big) + o_P(1/\sqrt n)
\eea

uniformly in $t \in [3h_n,1-3h_n],u,v\in [-1,1], y\in\Y$ where $\bar F_\eps$ is defined in Lemma \ref{lem:propfepsdach}. Now convolving both sides of (\ref{eq:c2h4}) [with respect to the argument $y$] with $\frac{1}{d_n}\omega(\cdot/d_n)$ and evaluating the result in $0$ yields the identity 
\bean
&&\frac{1}{n}\sum_i I_{[3h_n,t-h_n]}(X_i) \Big(\Omega\Big(\frac{Y_i^*-q_{\tau+vb_n}(X_i+uh_n)}{d_n}\Big) - \Omega\Big(\frac{\hat s_L(X_i)\eps_i^*}{d_n}\Big) \Big)
\\
&=& \frac{1}{n}\sum_i I_{[3h_n,t-h_n]}(X_i)\Big(\bar F_\eps\Big(\frac{q_{\tau+vb_n}(X_i+uh_n) - \hat q_{\tau,L}(X_i)}{\hat s_L(X_i)} \Big) -\bar F_\eps\Big(\frac{q_{\tau}(X_i) - \hat q_{\tau,L}(X_i)}{\hat s_L(X_i)} \Big)
\\
&&+ f_\eps(0)\frac{q_{\tau}(X_i) - \hat q_{\tau,L}(X_i)}{\hat s_L(X_i)} \Big) + o_P(1/\sqrt n).
\eean
Observe that the smoothness properties of $\bar F_\eps$ (defined in Lemma \ref{lem:propfepsdach}) yield the representation
\bean
&&\bar F_\eps\Big(\frac{q_{\tau}(X_i) - \hat q_{\tau,L}(X_i)}{\hat s_L(X_i)} \Big) - 
\bar F_\eps\Big(\frac{q_{\tau+vb_n}(X_i+uh_n) - \hat q_{\tau,L}(X_i)}{\hat s_L(X_i)} \Big)
\\
&& = vb_n\gamma_n(X_i)+ \sum_{j=1}^p \xi_j(X_i,v,n)(uh_n)^j + r_{n,1}
\eean
where the remainder terms $r_{n,1}$ is of order $O(b_n^2+h_n^{p+1}) = o(1/\sqrt n)$ uniformly in $u,v$ and $\xi_j,\gamma_n$ denote some functions that do not depend on $u$. Thus the proof of (\ref{eq:c2h1}) will be complete once we establish (\ref{eq:c2h4}). To this end, observe that 
\[
I\Big\{ Y_i^* \leq q_{\tau+vb_n}(X_i+uh_n)+y \Big\} 
= I\Big\{ \eps_i^* \leq \frac{q_{\tau+vb_n}(X_i+uh_n) - \hat q_\tau(X_i)+y}{\hat s(X_i)}\Big\}
\]
and
\bean
&&\frac{1}{n}\sum_i I_{[3h_n,t-3h_n]}(X_i) \Big(I\Big\{ \eps_i^* \leq \frac{q_{\tau+vb_n}(X_i+uh_n) - \hat q(X_i) +y}{\hat s(X_i)}\Big\} - I\Big\{\eps_i^* \leq \frac{y}{\hat s_L(X_i)}\Big\}\Big)
\\
&=& \frac{1}{n}\sum_i I_{[3h_n,t-3h_n]}(X_i) \Big(I\Big\{ \eps_i^* \leq \frac{q_{\tau+vb_n}(X_i+uh_n) - \hat q_{\tau,L}(X_i)+y}{\hat s_L(X_i)}\Big\} - I\Big\{\eps_i^* \leq \frac{y}{\hat s_L(X_i)}\Big\}\Big) + o_P(1/\sqrt n)
\\
&=& \frac{1}{n}\sum_i I_{[3h_n,t-3h_n]}(X_i) \Big(\tilde F_\eps\Big(  \frac{q_{\tau+vb_n}(X_i+uh_n) - \hat q_{\tau,L}(X_i)+y}{\hat s_L(X_i)} \Big) - \tilde F_\eps(y/\hat s_L(X_i)) \Big) + o_P(1/\sqrt n)
\eean
uniformly in $t,v,u$, which follows by arguments similar to those used in the proof of Lemma \ref{lem:proclin}. 
Consider the decomposition 
\bean
&&\tilde F_\eps\Big(  \frac{q_{\tau+vb_n}(X_i+uh_n) - \hat q_{\tau,L}(X_i)+y}{\hat s_L(X_i)} \Big) - \tilde F_\eps\Big(\frac{y}{\hat s_L(X_i)} \Big)
\\
&=& \tilde F_\eps\Big(  \frac{q_{\tau+vb_n}(X_i+uh_n) - \hat q_{\tau,L}(X_i)+y}{\hat s_L(X_i)} \Big) - \tilde F_\eps\Big(  \frac{q_{\tau}(X_i) - \hat q_{\tau,L}(X_i)+y}{\hat s_L(X_i)} \Big)
\\
&& + \tilde F_\eps\Big(  \frac{q_{\tau}(X_i) - \hat q_{\tau,L}(X_i)+y}{\hat s_L(X_i)} \Big)  - \tilde F_\eps\Big(\frac{y}{\hat s_L(X_i)} \Big) .
\eean
For the first term in this decomposition, an application of Lemma \ref{lem:propfepsdach} yields
\bean
&&\frac{1}{n}\sum_i I_{[3h_n,t-3h_n]}(X_i) \Big[ \tilde F_\eps\Big( \frac{q_{\tau+vb_n}(X_i+uh_n) - \hat q_{\tau,L}(X_i) + y}{\hat s_L(X_i)} \Big) - \tilde F_\eps\Big( \frac{q_{\tau}(X_i) - \hat q_{\tau,L}(X_i) + y}{\hat s_L(X_i)} \Big)\Big]
\\
&=& \frac{1}{n}\sum_i I_{[3h_n,t-3h_n]}(X_i) \Big[ \bar F_\eps\Big( \frac{q_{\tau+vb_n}(X_i+uh_n) - \hat q_{\tau,L}(X_i) + y}{\hat s_L(X_i)} \Big) - \bar F_\eps\Big( \frac{q_{\tau}(X_i) - \hat q_{\tau,L}(X_i) + y}{\hat s_L(X_i)} \Big)\Big]
\\
&& +  o_P(1/\sqrt n),
\eean
where $\bar F_\eps$ is defined in Lemma \ref{lem:propfepsdach}. Noting that
\[
\tilde F_\eps\Big(\frac{q_{\tau}(X_i) - \hat q_{\tau,L}(X_i)+y}{\hat s_L(X_i)} \Big) - \tilde F_\eps\Big(\frac{y}{\hat s_L(X_i)} \Big) = \tilde f_\eps\Big(\frac{y}{\hat s_L(X_i)} \Big)\frac{q_{\tau}(X_i) - \hat q_{\tau,L}(X_i)}{\hat s_L(X_i)} + o_P(1/\sqrt n),
\]
and recalling that $\tilde f_\eps$ converges to $f_\eps$ uniformly with rate $o_P((h_n/\log n)^{1/2})$ [see (\ref{le2n09})] combined with $r_n(h_n/\log n)^{1/2} = o(1)$ yields
\[
\tilde F_\eps\Big(  \frac{q_{\tau}(X_i) - \hat q_{\tau,L}(X_i)+y}{\hat s_L(X_i)} \Big)  - \tilde F_\eps\Big(\frac{y}{\hat s_L(X_i)} \Big) = f_\eps\Big(\frac{y}{\hat s_L(X_i)} \Big)\frac{q_{\tau}(X_i) - \hat q_{\tau,L}(X_i)}{\hat s_L(X_i)} + o_P(1/\sqrt n)
\] 
which completes the proof of (\ref{eq:c2h4}) and thus (\ref{eq:c2h1}) is also established.

\textbf{Proof of (\ref{eq:c2h2})} It suffices to show that, uniformly in $t \in [3h_n,1-3h_n]$
\bea \label{eq:c2h5}
 && \quad\quad \frac{1}{n}\sum_i I_{[3h_n,t-h_n]}(X_i)\Big( \Omega\Big(\frac{\hat s_L(X_i)\eps_i^*}{d_n} \Big) - I\{\eps_i^* \leq 0\} \Big) 
= o_P(1/\sqrt n),
\\ \label{eq:c2h6}
&& \quad\quad\frac{1}{n}\sum_i I_{[3h_n,t-h_n]}(X_i)(I\{\eps_i^* \leq 0\} - \tau) = \frac{1}{n}\sum_i I_{[3h_n,t]}(X_i)(I\{\eps_i^* \leq 0\} - \tau) + o_P(1/\sqrt n),
\\ \label{eq:c2h7}
&& \quad\quad\frac{1}{n}\sum_i I_{[3h_n,t-h_n]}(X_i)\frac{q_{\tau}(X_i) - \hat q_{\tau,L}(X_i)}{\hat s_L(X_i)}
= \int_{3h_n}^{t} \frac{q_{\tau}(u) - \hat q_{\tau,L}(u)}{s(u)} f_X(u) du + o_P(1/\sqrt n).
\eea

The statement in (\ref{eq:c2h7}) follows since, for $t \in [4h_n,1-3h_n]$,
\bean
\frac{1}{n}\sum_i I_{[3h_n,t-h_n]}(X_i)\frac{q_{\tau}(X_i) - \hat q_{\tau,L}(X_i)}{\hat s_L(X_i)}
&=& \frac{1}{n}\sum_i I_{[3h_n,t-h_n]}(X_i)\frac{q_{\tau}(X_i) - \hat q_{\tau,L}(X_i)}{s(X_i)} + o_P(1/\sqrt n)
\\
&=& \int_{3h_n}^{t-h_n} \frac{q_{\tau}(u) - \hat q_{\tau,L}(u)}{s(u)} f_X(u) du + o_P(1/\sqrt n)
\\
&=& \int_{3h_n}^{t} \frac{q_{\tau}(u) - \hat q_{\tau,L}(u)}{s(u)} f_X(u) du + o_P(1/\sqrt n),
\eean
where the first equality follows from the rates of convergence for $\hat q_{\tau,L} - q_\tau, \hat s_L - s$ [see Lemma \ref{lem:unifratesquant} and Lemma \ref{lem:quantlin}], the second equality is a consequence of  the fact that $\hat q_{\tau,L} \in C_C^{\delta}$ with probability tending to one [see Lemma \ref{lem:unifratesquant}] combined with Lemma \ref{lem:base}. For $t < 4h_n$, the left-hand side of (\ref{eq:c2h6}) is zero and the right-hand side of order $o_P(n^{-1/2})$ by Lemma \ref{lem:unifratesquant} and Lemma \ref{lem:quantlin}.\\
For a proof of (\ref{eq:c2h5}), observe that
\[
\Omega\Big(\frac{\hat s_L(X_i)\eps_i^*}{d_n}\Big) - I\{\eps_i^* \leq 0\} = \frac{1}{d_n}\int_{-d_n}^{d_n} \Big(I\{\eps_i^* \leq a/\hat s_L(X_i)\} - I\{\eps_i^* \leq 0\}\Big) \omega\Big(\frac{a}{d_n}\Big) da.
\]
Define the sequence of sets 
\[
S(\delta_n) := \{(t,y_n,z_n)|t\in[3h_n,1-3h_n],y_n,z_n\in \Y, |y_n-z_n| \leq \delta_n\}
\] 
for some $\delta_n = o(1)$. Observe that, with probability tending to one, 
\bean
&&\sup_{(t,y_n,z_n)\in S(\delta_n)} \Big|\frac{1}{n}\sum_{i=1}^n I_{[3h_n,t-3h_n]}(X_i)\Big(I\{\eps_i^*\leq y_n\} - I\{\eps_i^*\leq z_n\} + \tilde F_\eps(z_n) - \tilde F_\eps(y_n)\Big)\Big|
\\
&=& \sup_{(t,y_n,z_n)\in S(\delta_n)} \Big|\frac{1}{n}\sum_{i=1}^n I_{[3h_n,t-3h_n]}(X_i)\Big(I\{U_i\leq \tilde F_\eps(y_n)\} - I\{U_i\leq \tilde F_\eps(z_n)\} + \tilde F_\eps(z_n) - \tilde F_\eps(y_n)\Big)\Big|
\\
&\leq& \sup_{(t,y_n,z_n)\in S(C\delta_n)}\Big|\frac{1}{n}\sum_{i=1}^n I_{[3h_n,t-3h_n]}(X_i) \Big(I\{U_i \leq y_n\} - I\{U_i \leq z_n\} + z_n - y_n\Big)\Big| 
\\
&=& o_P(1/\sqrt n).
\eean
Here, for the first inequality we made use of (B.1). This implies that, with probability tending to one, $\tilde F_\eps$ has a uniformly bounded derivative which shows that, with probability tending to one, $|y_n-z_n|\leq \delta_n$ implies $|\tilde F_\eps(y_n) - \tilde F_\eps(z_n)| \leq C\delta_n$ for some finite constant $C$. The last bound above follows by standard empirical process arguments provided that $\delta_n = o(1)$. Thus
\bean
&& \frac{1}{n}\sum_i I_{[3h_n,t-h_n]}(X_i)\Big( \Omega\Big(\frac{\hat s_L(X_i)\eps_i^*}{d_n}\Big) - I\{\eps_i^* \leq 0\}\Big)
\\ 
&=& \frac{1}{n}\sum_i I_{[3h_n,t-h_n]}(X_i) \frac{1}{d_n}\int_{-d_n}^{d_n} \Big(\tilde F_\eps(a/\hat s_L(X_i)) - \tilde F_\eps(0)\Big) \omega\Big(\frac{a}{d_n}\Big) da
 + o_P(n^{-1/2})
\\
&=& \frac{1}{n}\sum_i I_{[3h_n,t-h_n]}(X_i) \frac{1}{d_n}\int_{-d_n}^{d_n} \Big(\bar F_\eps(a/\hat s_L(X_i)) - \bar F_\eps(0)\Big) \omega\Big(\frac{a}{d_n}\Big) da
 + o_P(n^{-1/2})
\\
&=& o_P(n^{-1/2})
\eean
where the second to last line follows by Lemma \ref{lem:propfepsdach} and the last line is a consequence of the smoothness properties of $\bar F_\eps$. \\
Thus (\ref{eq:c2h5}) follows and it remains to establish (\ref{eq:c2h6}). To this end, observe that it suffices to establish
\[
\sup_{t \in [3h_n,1-3h_n]} \Big|\frac{1}{n}\sum_{i=1}^n I_{[t-h_n,t]}(X_i)(I\{\eps_i^* \leq 0\} - \tau)\Big| = o_P(n^{-1/2}).
\]
Now
\[
\frac{1}{n}\sum_{i=1}^n I_{[t-h_n,t]}(X_i)(I\{\eps_i^* \leq 0\} - \tau) = \frac{1}{n}\sum_{i=1}^n I_{[t-h_n,t]}(X_i)(I\{U_i \leq \tilde F_\eps(0)\} - \tau),
\]
and by (\ref{eq:c72}) in Lemma \ref{lem:propfepsdach} we have $\tilde F_\eps(0) - \tau = \tilde F_\eps(0) - F_\eps(0) = O_P(r_n)$. Thus we have with probability tending to one $|\tilde F_\eps(0) - \tau|\leq r_n h_n^{-1/4}$ and in particular
\bean
&&\sup_{t \in [3h_n,1-3h_n]} \Big|\frac{1}{n}\sum_{i=1}^n I_{[t-h_n,t]}(X_i)(I\{\eps_i^* \leq 0\} - \tau)\Big| \\
&\leq& \sup_{t \in [3h_n,1-3h_n]} \sup_{|y| \leq r_n h_n^{-1/4}} \Big|\frac{1}{n}\sum_{i=1}^n I_{[t-h_n,t]}(X_i)(I\{U_i \leq y\} - \tau)\Big| = o_P(n^{-1/2})
\eean
where the first inequality holds with probability tending to one and the equality follows by standard empirical process arguments. Thus (\ref{eq:c2h5}) follows. This completes the proof of Lemma \ref{lem:B2boot}.

\hfill $\Box$

\end{appendix}

\end{document}